\documentclass[reqno, 11pt, german]{amsart}

\usepackage{amssymb,bbm,amsmath,enumitem}             
\usepackage{float}                
\usepackage{graphicx, marvosym}
\usepackage{german}
\usepackage{hyperref,color}
\definecolor{darkblue}{rgb}{0,0,0.7}
\definecolor{darkgreen}{rgb}{0,0.4,0}
\definecolor{darkred}{rgb}{0.7,0,0}
\hypersetup{breaklinks=true, colorlinks=true, urlcolor=darkblue, linkcolor=darkblue, citecolor=darkblue}

\numberwithin{equation}{section}
\begin{document}
\title[Ein optimiertes Gl\"attungsverfahren]{Ein optimiertes Gl\"attungsverfahren motiviert durch eine technische Fragestellung}
\author{Frank Klinker \and G\"unter Skoruppa}
\thanks{{\Large\Letter}: Fakult\"at f\"ur Mathematik, TU Dortmund, 44221 Dortmund, Germany}
\thanks{{\large\Email}\ : \href{mailto:frank.klinker@math.tu-dortmund.de}{frank.klinker@math.tu-dortmund.de} }
\thanks{{\large\Email}\ : \href{mailto:guenter.skoruppa@math.tu-dortmund.de}{guenter.skoruppa@math.tu-dortmund.de} }
\thanks{\vspace*{1ex} Math.~Semesterber.~{\bf 59} (2012), no.1, 29-55.
\href{http://dx.doi.org/10.1007/s00591-012-0098-1}{DOI 10.1007/s00591-012-0098-1}} 

\begin{abstract} Ausgehend von einer konkreten technischen Fragestellung diskutieren wir in dieser Notiz die Anwendung verschiedener Gl\"attungsverfahren auf Datens\"atze mit vorgegebener Struktur. Wir stellen die Verfahren im Detail vor und besprechen die Vor- und Nachteile. Insbesondere stellen wir hier die symmetrisierte exponentielle Gl\"attung vor, die ein sehr gutes Gl\"attungsverhalten mit einem hohen Ma\ss~an Symmetrieerhaltung kombiniert.
\end{abstract}
\maketitle
\setcounter{tocdepth}{2}

\section{Ausgangssituation und Fragestellung}\label{einl}

\begin{figure}[h]\caption{Keilsicherungsscheibe DIN 25201}\label{din}
 \vspace*{1.5ex}\includegraphics[width=5.2cm]{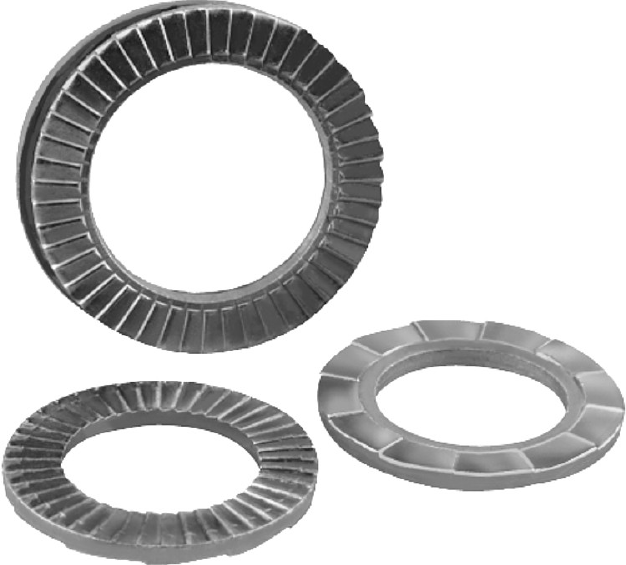}
\end{figure}

\begin{figure}[h]\caption{Schematischer Aufbau der Messapparatur}\label{laser}
\includegraphics[width=0.5\textwidth]{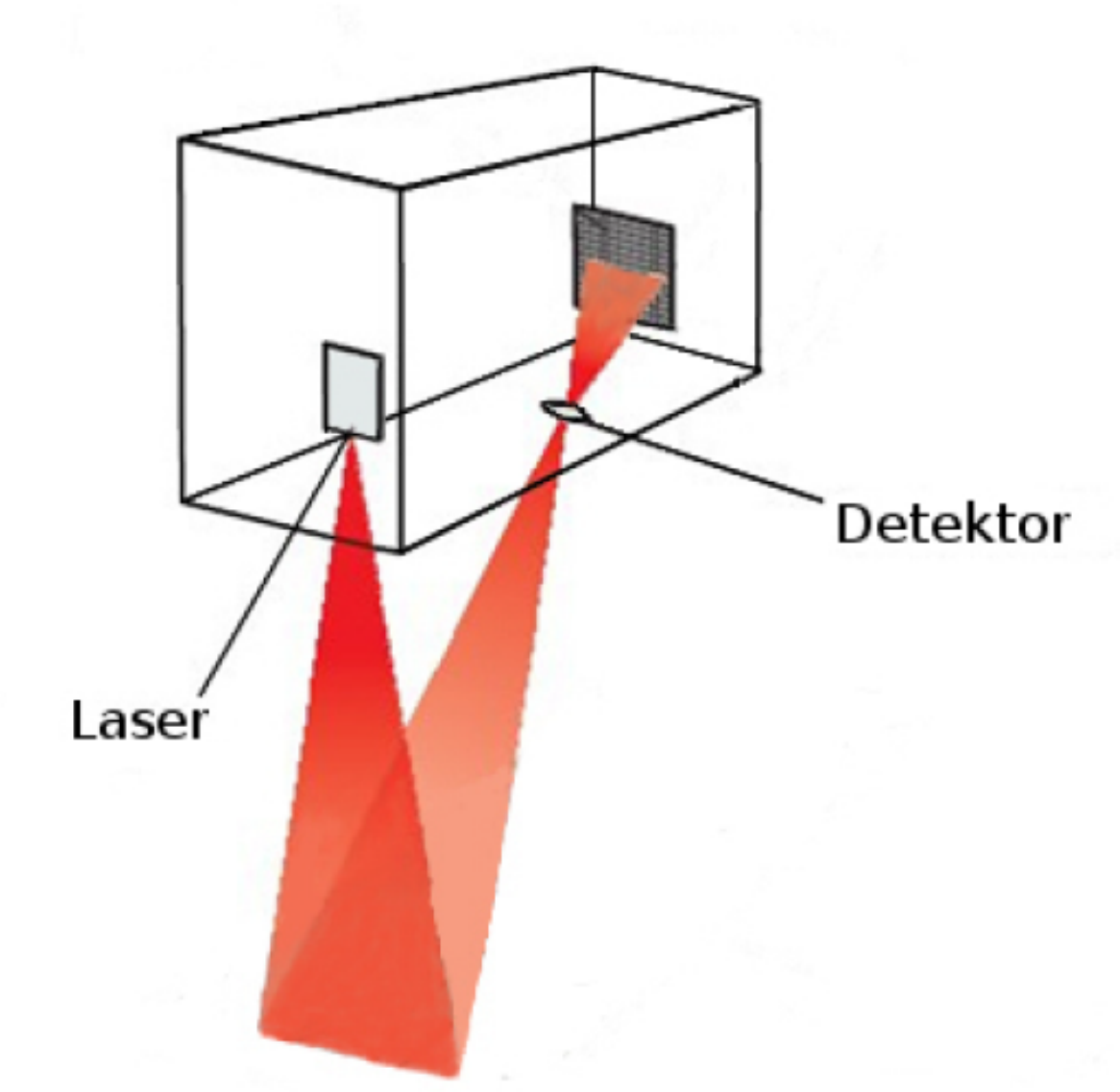} \hspace{2em} 
\includegraphics[width=0.4\textwidth]{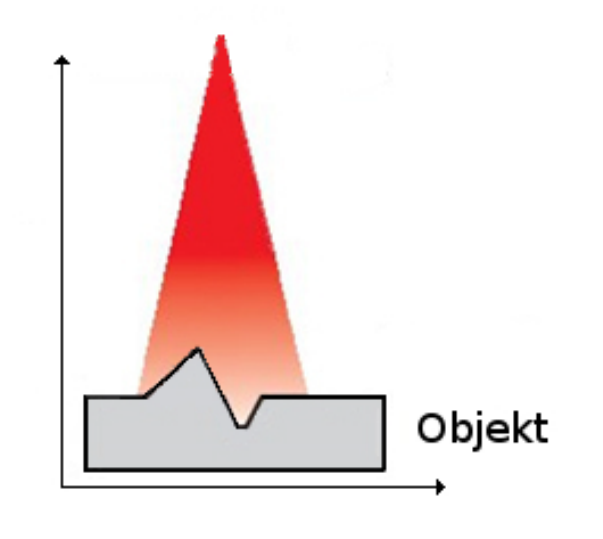} 
\end{figure}

Wir pr\"asentieren in dieser Notiz eine L\"osung zur folgenden technischen Problemstellung. 
Gegeben ist das H\"ohenprofil einer Seite einer Keilsicherungsscheibe (DIN 25201). Diese Scheibe ist beidseitig gezahnt mit jeweils einer grob und einer fein gezahnten Seite, siehe Abbildung \ref{din}.
Die Ermittlung des Profils erfolgt mit Hilfe eines 2D/3D-Lasersensors. Dieser liefert  \"uber $N\leq 750$ \"aquidistanten St\"utzstellen einer Sekante der Scheibe eine Reihe von Messpunkten, wobei jeder einzelne Messpunkt der H\"ohe des Objekts \"uber der jeweiligen St\"utzstelle entspricht. Zum schematischen Aufbau der Messapparatur vergleiche Abbildung \ref{laser}.
Das Ziel ist, die Zahl der Z\"ahne der ausgemessenen Scheibenseite zu ermitteln. Dies geschieht durch Ausz\"ahlen der Maxima des Datensatzes. Dieses Ausz\"ahlen und die daraus resultierende Entscheidung, ob sich die Scheibe in der richtigen Lage befindet oder gedreht werden muss, soll in einem  automatisierten Prozess geschehen.

Dabei tritt das Problem auf, dass die Messwerte durch Verunreinigungen der Scheibe vom Soll abweichen k\"onnen, zum Beispiel durch Verschmutzung w\"ahrend des Arbeitsprozesses oder auch durch Farbauftr\"age oder eine nicht ganz fehlerfreie Produktion des Werkst\"ucks. 
Geringe Abweichungen, insbesondere solche durch Verunreinigungen, sollen toleriert werden und die Zahnz\"ahlung nicht behindern.

\begin{figure}[h!]\caption{Typische und pathologische Messreihe}\label{fig1}
 \includegraphics[width=0.48\textwidth]{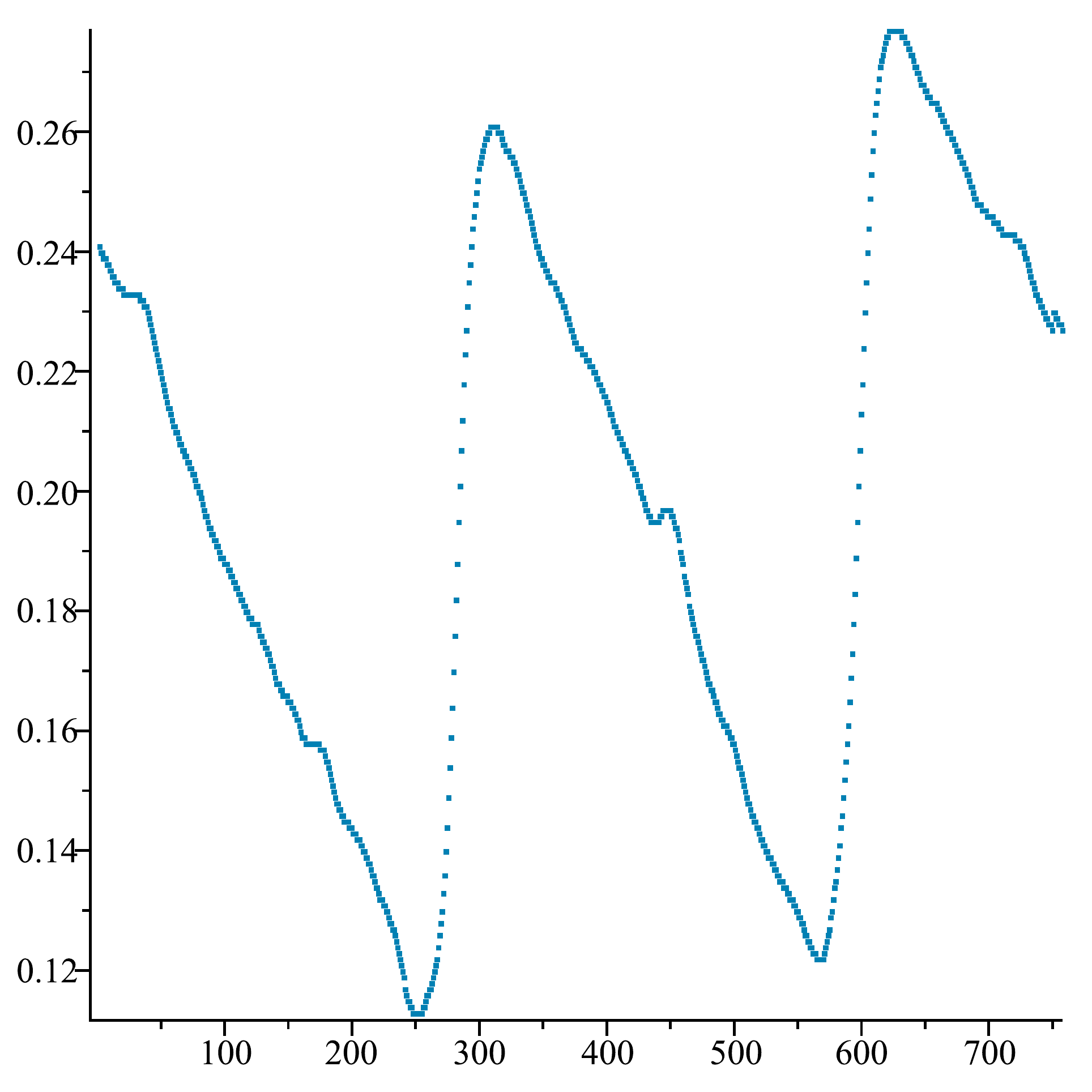}\
 \includegraphics[width=0.48\textwidth]{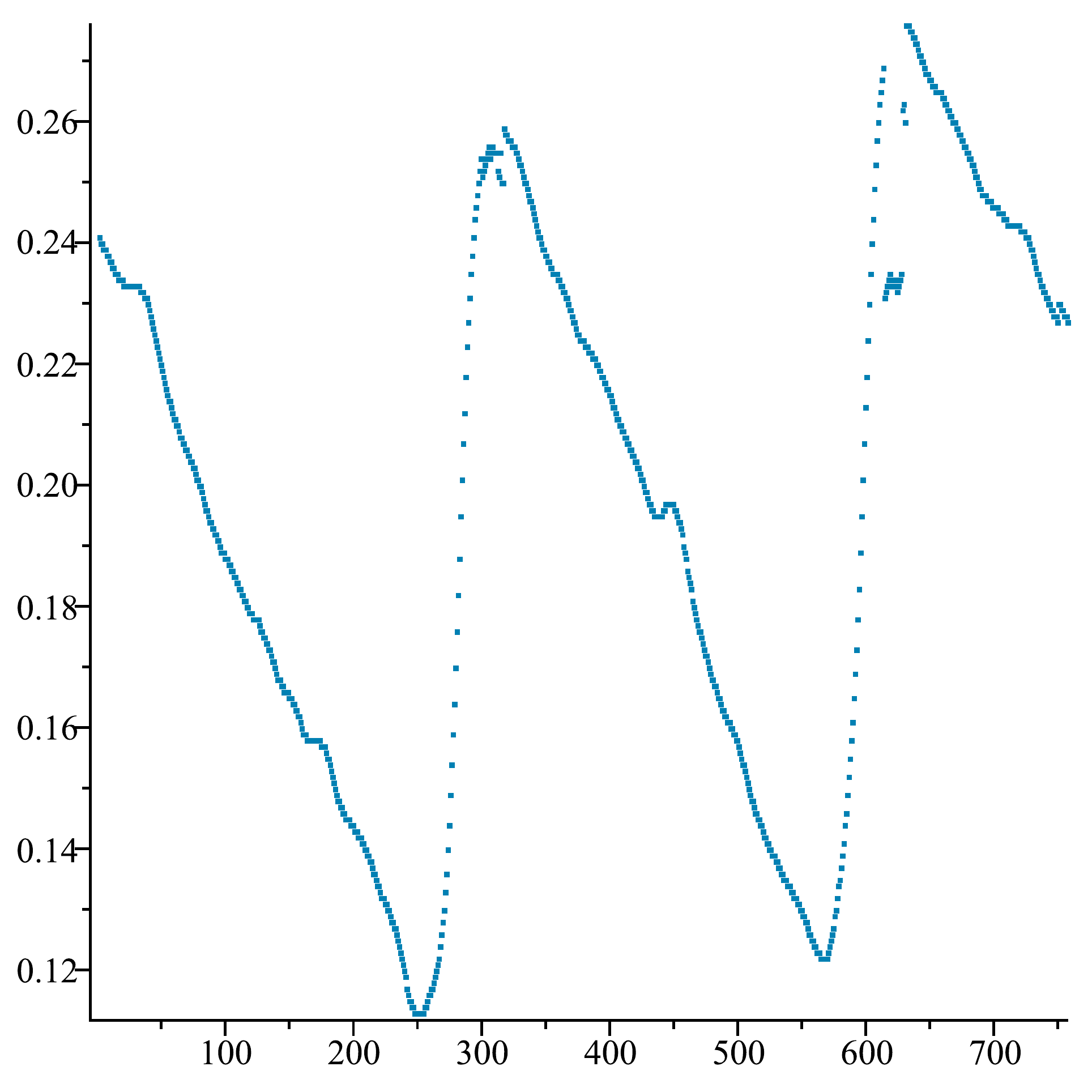} 
\end{figure}

Die Grafiken in Abbildung \ref{fig1} zeigen eine typische und eine hieraus erzeugte pathologische Messreihe.
Gut erkennbar sind in Abbildung \ref{fig1} die Imperfektionen an den nach rechts absteigenden Flanken. Diese werden als zus\"atzliche Z\"ahne die Z\"ahlung verf\"alschen. 
Der rechte Teil der Abbildung \ref{fig1} wurde aus dem linken durch manuell eingef\"ugte, lokale St\"orungen in den Messwerten erzeugt. Solche St\"orungen treten auch in realen Beispielen auf, wie die Messreihe {\bf I} zeigt.

\subsection*{Fragestellung}
Lassen sich  unwesentliche St\"orungen derart beseitigen, dass das bereinigte Zahnprofil durch Ausz\"ahlung der Spitzen unmittelbar die Entscheidung der technischen Fragestellung erm\"oglicht?
Dabei darf die urspr\"ungliche Messreihe ver\"andert werden mit der Einschr\"ankung, dass die Zahnzahl erhalten bleiben soll. Die Modifikation der Daten wird dabei umso besser sein, je mehr der "`Charakter"' der urspr\"unglichen Daten erhalten bleibt.

\section{L\"osungsans\"atze.} \label{glaettung}

Gl\"attungsoperationen, wie sie etwa in der Statistik eingesetzt werden, sind ein geeignetes Mittel zur St\"orungsbeseitigung in den Daten. 
Wir werden uns hier zun\"achst auf die Vorstellung zweier g\"angiger Verfahren beschr\"anken und ihre Vor- und Nachteile  diskutieren. Im Einzelnen sind das der {\em gleitende arithmetische Durchschnitt} ({\em moving arithmetic mean}, MA), und die {\em exponentielle Gl\"attung} ({\em exponential average}, EA). Zur Diskussion allgemeiner Durchschnitte und ihrer gleitenden Varianten verweisen wir gerne auf \cite{Kl}.
Anschliessend bieten wir ein Gl\"attungsverfahren an, das jeweils die Vorteile beider vereint. Dieses Verfahren nennen wir {\em symmetrisierte exponentielle Gl\"attung}, ({\em symmetrized exponential average}, SEA).

\subsection{Der gleitende arithmetische Durchschnitt (MA)}\label{subMA}

Gegeben sei ein Datensatz\footnote{Ist der Datensatz endlich, so sei im Folgenden immer $N$ die Anzahl der Elemente.} $(y_1,y_2,\ldots,y_N,\ldots)$. Die Idee des gleitenden arithmetischen Mittels ist es nun, statt des Wertes $y_i$ ein arithmetisches Mittel $\bar{y}_i$ eines Teils des Datensatzes zu ersetzen.
Die Berechnung des gleitenden arithmetischen Durchschnitts ben\"otigt zwei Parameter. Zum einen die {\em Gleitl\"ange}  $n+1\leq N$ zum anderen das {\em Gewicht} $\ell\in\{0,\ldots, n\}$. Die Spezialf\"alle $\ell=0$ und $\ell=n$ nennt man {\em rechts-} und {\em linkseitigen gleitenden arithmetischen Durchschnitt der L\"ange $n+1$} und den Spezialfall $n$ gerade und $\ell=\frac{n}{2}$ nennt man {\em zentrierten gleitenden arithmetischen Durchschnitt der L\"ange $n+1$}.
Das oben angesprochene arithmetische Mittel wird dann aus den $n+1$ Werten der Menge $\{y_{i-\ell},\ldots,y_{i+n-\ell}\}$ berechnet, also

\begin{equation}\label{MA}\begin{aligned}
 \bar{y}_i 
& :=  \frac{1}{n+1}\left(y_{i-\ell} + y_{i-\ell+1} + \ldots + y_i + \ldots  + y_{i+n-\ell-1} + y_{i+n-\ell} \right) \\
& =   \frac{1}{n+1} \sum_{k=-\ell}^{n-\ell} y_{i+k}\,,
\end{aligned}\end{equation}
f\"ur $i=\ell+1, \ldots, N-n+\ell$. Der neue Datensatz ist insbesondere um $n$ Werte k\"urzer als der alte, und zwar entfallen am Anfang $\ell$ und am Ende $n-\ell$ Werte. Ist der Datensatz nicht endlich, so entfallen selbstverst\"andlich nur die Werte am Anfang.

Wie man die Parameter w\"ahlt, h\"angt von der Anwendung ab. Da in unserer Anwendung alle Werte der Messreihe gleichberechtigt sind, haben wir uns f\"ur die zentrierte Variante entschieden und MA bezeichnet im Folgenden genau diese.\footnote{Die Beschr\"ankung auf die zentrierte Variante ist insofern nicht zwingend notwendig, da einerseits die Wahl von $\ell\in\{0,\ldots,n\}$ unerheblich f\"ur das Ergebnis der Gl\"attung ist, und andererseits in dem hier vorliegenden Fall $n\ll N$ auch die Symmetrie der Originalmessreihe im Wesentlichen erhalten bleibt.} Hier ist der Verlust an Daten am Anfang und am Ende des Datensatzes gleich. In diesem Fall ist dann $n=2\ell$ und (\ref{MA}) wird zu
\begin{equation}\label{ZMA}
 \bar{y}_i  = \frac{1}{2\ell+1} \sum_{k=-\ell}^{\ell} y_{i+k}\,.
\end{equation}
Zur Berechnung des MA muss man nicht in jedem Schritt die gesamte Summe  (\ref{MA}) oder (\ref{ZMA}) berechnen. Die Folge der Durchschnitte $\bar y_i$ erf\"ullt vielmehr die folgende rekursive Gleichung:
\begin{equation}\label{MArek}
\bar{y}_i =  \bar{y}_{i-1}+  \frac{1}{2\ell+1} \left(  y_{i+\ell} - y_{i-1-\ell} \right)\,.
 \end{equation}
Daher sieht ein den Aufwand minimierender Berechnungsalgorithmus f\"ur den zentrierten MA wie folgt aus.

\begin{enumerate}[leftmargin=2em]
\item W\"ahle einen positiven ganzzahligen Parameter $\ell$, so dass $2\ell+1\leq N$.
\item Berechne den Startwert 
\[ 
\bar{y}_{\ell+1} :=  \frac{1}{2\ell+1}\left( y_{1} + y_{2} + \ldots + y_{2\ell+1} \right)\,.
\]
\item Berechne f\"ur $i = \ell+2, \ldots, N-\ell$ rekursiv 
\[
 \bar{y}_{i} :=  \bar{y}_{i-1}+  \frac{1}{2\ell+1} \left(  y_{i+\ell} - y_{i-1-\ell} \right) \,.
\]
\end{enumerate}

\subsection*{Bemerkungen zum MA}
\begin{itemize}[leftmargin=1.5em]
\item Je gr\"o\ss er die Gleitl\"ange $n+1\in\mathbbm{N}$ gew\"ahlt wird, desto st\"arker wird die Gl\"attung ausfallen. Zu gro\ss e Gleitl\"angen $n+1$ mitteln die Werte so stark, dass gegebenenfalls wichtige charakteristische Eigenschaften der Originalmessreihe verloren gehen. Mit anderen Worten, der Einfluss des einzelnen Messwertes $y_i$ geht verloren, weil alle Werte mit dem gleichen Gewicht $\frac{1}{n+1}$ in den Durchschnitt eingehen.
\item Neben dem eventuellen Verlust an charakteristischen Eigenschaften hat der MA insbesondere bei gro\ss en Gleitl\"angen den Nachteil, dass sie den Datensatz stark verk\"urzen. 
\item F\"ur den Grenzfall $2\ell+1=N$ erh\"alt man genau einen Wert $\bar y_{\ell+1}$, der dem gew\"ohnlichen arithmetischen Mittel aller Messpunkte entspricht, und f\"ur $\ell=0$ bekommt man als Ergebnis die urspr\"ungliche Messreihe zur\"uck.
\item Sieht man einmal von der einmaligen Startwertberechnung ab, ben\"otigt die Berechnung eines neuen Datenwertes  3 Rechenoperationen, vgl.~(\ref{MArek}).
\item Der MA ber\"ucksichtigt die Symmetrie eines Systems: Hat der beschr\"ankte Ausgangsdatensatz die Symmetrie $y_i=y_{N-i}$, so hat der gegl\"attete Datensatz diese ebenfalls. 
\end{itemize}

\subsection{Die exponentielle Gl\"attung (EA)}\label{subEA}

Bei dem in Abschnitt \ref{subMA} vorgestellten Verfahren der Gl\"attung wird nur ein Ausschnitt des gesamten Datensatzes um einen ausgezeichneten Wert $y_i$ betrachtet und alle Werte dieses Abschnitts tragen mit dem gleichen Gewicht $\frac{1}{n+1}$ zur Gl\"attung, also zum Wert $\bar y_{i}$, bei, siehe (\ref{MA}).
Bei der exponentiellen Gl\"attung werden im Gegensatz dazu bei der Berechnung alle Werte der Testreihe bis zur betrachteten Stelle $y_i$ mit einbezogen. Dabei tragen allerdings die "{}\"alteren Werte{}" -- das sind die Werte im Datensatz, die sich vor der betrachteten Stelle $y_i$ befinden -- mit fallendem Gewicht zur Berechnung der Gl\"attung bei.  
Die Definition des EA erfolgt rekursiv, analog zur Berechnung des MA gem\"a\ss\ (\ref{MArek}).

Gegeben sei ein endlicher oder unendlicher Datensatz $(y_1,y_2, \ldots,y_N,\ldots)$. Der EA dieses Datensatzes ist gegeben durch
\begin{equation}\label{EArek}\begin{aligned}
\hat y_1 	 &:=y_1\,, \\
\hat y_i &:= (1-\alpha)\hat{y}_{i-1} + \alpha y_{i},\quad\text{f\"ur }i=2,\ldots,N \,.
\end{aligned}\end{equation}
Die Berechnung des EA h\"angt hier von einem Parameter $0\leq\alpha\leq 1$ ab. F\"ur $\alpha\neq 0$ hei\ss t das Inverse dieses Parameters, $\frac{1}{\alpha}$, die {\em Gleitl\"ange} der exponentiellen Gl\"attung. Wie auch im Fall des MA gibt es hier eine geschlossene Berechnungsformel. In dieser sieht man auch das Fallen der Gewichte, mit dem die \"alteren Datenwerte zum Durchschnittswert beitragen. Es gilt
\begin{equation}\begin{aligned}\label{EA}
\hat y_i &= (1-\alpha)^{i-1}y_1+ \alpha (1-\alpha)^{i-2} y_{2}+\ldots +\alpha(1-\alpha)y_{i-1}+\alpha y_i\\
&= (1-\alpha)^{i-1}y_1+  \alpha \sum_{r=2}^{i}(1-\alpha)^{i-r} y_{r}\,.
\end{aligned}\end{equation}

Die algorithmische Beschreibung der Berechnung des EA erfolgt wieder \"uber die rekursive Formel (\ref{EArek}).
\begin{enumerate}
\item W\"ahle Parameter $\alpha \in [0,1]$ und setze $\hat{y}_1 := y_1$.
\item Berechne rekursiv f\"ur\footnote{Ist der Datensatz unendlich, so entf\"allt selbstverst\"andlich die obere Beschr\"ankung hier, genauso wie schon in (\ref{EArek}).} $i = 2,3, \ldots, N$:
\begin{equation*}
\hat{y}_{i} = (1-\alpha)  \hat{y}_{i-1} + \alpha y_{i}  \,.
\end{equation*}
\end{enumerate}

\subsection*{Bemerkungen zum EA}
\begin{itemize}[leftmargin=1.5em]
\item Anders als der zentrierte MA ist das Verfahren einseitig, d.h.~es betrachtet von $y_i$ aus gesehen nur fr\"uhere Daten. Dies macht EA besonders interessant f\"ur Prognosen, etwa bei B\"orsenkursen. F\"ur eine Diskussion dieser Anwendung siehe zum Beispiel \cite{Kl}.
\item Dass die exponentielle Gl\"attung einseitig arbeitet und die komplette Historie reflektiert, kann bei nahezu periodischen oder anderen Symmetrien aufweisenden Datens\"atzen von Nachteil sein, da die Einseitigkeit diese Symmetrie zerst\"ort.  Dieses Ph\"anomen wird insbesondere am Beispiel in Abschnitt \ref{sub4} deutlich werden.
\item Die Einseitigkeit des EA hat auch einen negativen Einfluss auf die Gl\"attung, insbesondere wenn ein Bereich, der eine starke Gl\"attung erfordert, recht fr\"uh im Datensatz zu finden ist. Vergleiche dazu Abbildung \ref{fig23} und beachte, dass dort der rEA gem\"a\ss~(\ref{rEA}) angewendet wurde und somit der fr\"uhe Bereich rechts in den Kurven zu finden ist.
\end{itemize}

\subsection{Die symmetrisierte exponentielle Gl\"attung (SEA)}\label{subSEA}

In diesem Abschnitt liefern wir eine Variante der exponentiellen Gl\"attung, die den Nachteil der Einseitigkeit des EA kompensiert. Dieses Gl\"attungsverfahren kombiniert die guten Gl\"attungseigenschaften des EA und dessen Betonung des Wertes $y_i$ bei der Berechnung der Gl\"attung $\bar y_i$, einerseits, mit der Symmetrieerhaltung des MA, andererseits. Zur Konstruktion ben\"otigen wir einen endlichen Datensatz und wir schauen uns zur Berechnung der Gl\"attung den gesamten Datensatz an und nicht -- wie beim EA -- nur den \glqq linken Teil\grqq~oder -- wie beim MA -- nur einen Ausschnitt.  

Sei also der endliche Datensatz $(y_1,y_2,\ldots ,y_N)$ gegeben. Zu diesem konstruieren wir durch Auslesen des Datensatzes von rechts nach links einen neuen Datensatz $(z_1,\ldots,z_N)$ mit $z_j:=y_{N+1-j}$. Wir definieren dann die {\em exponentielle Gl\"attung von rechts} (rEA) des Datensatzes $(y_1,\ldots ,y_N)$ als die exponentielle Gl\"attung des Datensatzes $(z_1,\ldots,z_N)$. Diesen bezeichnen wir mit $\check{y}_i:=\hat z_{N-i+1}$. Der rEA berechnet sich dann analog zu (\ref{EA}) gem\"a\ss
\begin{equation}\label{rEA}\begin{aligned}
\check{y}_i &=  (1-\alpha)^{N-i}y_N+ \alpha(1-\alpha)^{N-i-1}y_{N-1}+\ldots +\alpha(1-\alpha)y_{i+1} +\alpha y_i  \\
			&= (1-\alpha)^{N-i}y_N +\alpha\sum_{r=i}^{N-1} (1-\alpha)^{r-i}y_{r}\,.
\end{aligned}\end{equation}
Die rekursive Variante berechnet sich dann analog zu (\ref{EArek}) zu
\begin{equation}\begin{aligned}\label{rEArek}
\check{y}_N  &:= y_N\,, \\
\check{y}_i  &:= (1-\alpha)\check{y}_{i+1}+\alpha y_i \quad\text{f\"ur }i=N-1,\ldots, 1 \,.
\end{aligned}\end{equation}

\subsection*{Definition: SEA} 
Die {\em symmetrisierte exponentielle Gl\"attung} (SEA) des Datensatzes ist definiert als das arithmetische Mittel aus EA und rEA:
\begin{equation}\label{SEA}
\bar y_i:=\frac{1}{2}\left(\hat y_i+\check {y}_i\right)\,.
\end{equation}

Die entsprechende algorithmische Beschreibung des SEA lautet
\begin{enumerate}[leftmargin=2em]
\item W\"ahle Parameter $\alpha \in [0,1]$ und setze $\hat{y}_1 := y_1$ und $\check{y}_N=y_N$.
\item Berechne rekursiv f\"ur $i = 2,\ldots, N$ den EA 
\begin{equation*}
\hat{y}_i = (1-\alpha)  \hat{y}_{i-1} + \alpha y_i  \,,
\end{equation*}
und f\"ur $j=N-1,\ldots,1$ den rEA
\begin{equation*}
\check{y}_j = (1-\alpha) \check{y}_{j+1} +\alpha y_j\,.
\end{equation*}
\item Berechne f\"ur $i=1,\ldots, N$ den SEA
\begin{equation*}
\bar y_i=\frac{1}{2}(\hat y_i+\check y_i)\,.
\end{equation*}
\end{enumerate}

\subsection*{Bemerkungen zum EA, rEA und SEA}
\begin{itemize}[leftmargin=1.5em]
\item Bei der Berechnung des EA, rEA und des SEA wird der Datensatz im Gegensatz zum MA nicht verk\"urzt. Das ist insbesondere von Vorteil, wenn zur Gl\"attung hohe Gleitl\"angen ben\"otigt werden.
\item Die Berechnung eines neuen Datenwertes ben\"otigt f\"ur EA und rEA wie beim MA ebenfalls 3 Operationen, vergleiche(\ref{EArek}). Zur Berechnung des SEA ben\"otigt man insgesamt 8 Operationen.
\item Zur Berechnung eines Wertes des SEA werden alle Werte der Ursprungsreihe mit herangezogen und nicht, wie beim MA und EA, nur ein Teil der Werte.
\item Je kleiner $\alpha$, desto st\"arker wird die Gl\"attung ausfallen. Das bedeutet aber gleichzeitig, dass mehr  Charakteristiken des Datensatzes verloren gehen. Die Grenzf\"alle liefern f\"ur  $\alpha=1$ keine Ver\"anderung der Messreihe und f\"ur $\alpha=0$ eine konstante Folge, n\"amlich das arithmetische Mittel aus $y_1$ und $y_N$.
\item Bei den Berechnungen des EA und des SEA geht der Wert $y_i$  mit dem Gewicht $\alpha$ ein. Bei der Berechnung des MA in (\ref{ZMA}) geht $y_i$ mit dem Gewicht $\frac{1}{n+1}$ ein. Das und auch die Diskussion der Grenzf\"alle motiviert die Interpretation und die Definition des Parameters  $\frac{1}{\alpha}$ als Gleitl\"ange.
\item Im Gegensatz zu MA geht der Wert $y_i$ bei der Berechnung des SEA mit dem gr\"o\ss ten Gewicht ein. Dies gilt schon f\"ur den EA und den rEA, vgl.~(\ref{EA}) und (\ref{rEA}).
\end{itemize}

\section{Beispiele}\label{Ex}

In den Abschnitten \ref{sub1} und \ref{sub2} widmen wir uns den vorgestellten Beispielen realer Datens\"atze gem\"ass Abbildung \ref{fig1}.  Daneben  betrachten wir noch weitere Beispiele um die vorgeschlagenen Gl\"attungsverfahren zu testen und ihre  Eigenschaften zu verdeutlichen. Dies sind im Einzelnen in Abschnitt \ref{sub3} eine \"Uberlagerung  zweier Schwingung und in Abschnitt \ref{sub4} ein  fourierentwickeltes Rechtecksignal.

Die in den Abschnitten \ref{sub1} und \ref{sub6} verwendeten Datens\"atze wurden uns freundlicherweise von der Firma \href{http://www.kranz-software-engineering.de/}{\sc Kranz Software Engineering} zur Ver\-f\"u\-gung gestellt.

\subsection{Laserausmessung eines Zahnscheibenprofils}\label{sub1}

Zun\"achst untersuchen wir eine typische Messreihe. Bei dieser ist an einer absteigenden Flanke ein zus\"atzliches Maximum zu erkennen, das bei der Z\"ahlung der Maxima allerdings vernachl\"assigt werden soll.

\begin{figure}[H]\caption{Originaldaten einer Laserausmessung und Gl\"attung mittels MA ($n+1=20$)}\label{fig10}\centering
 \includegraphics[width=0.48\textwidth]{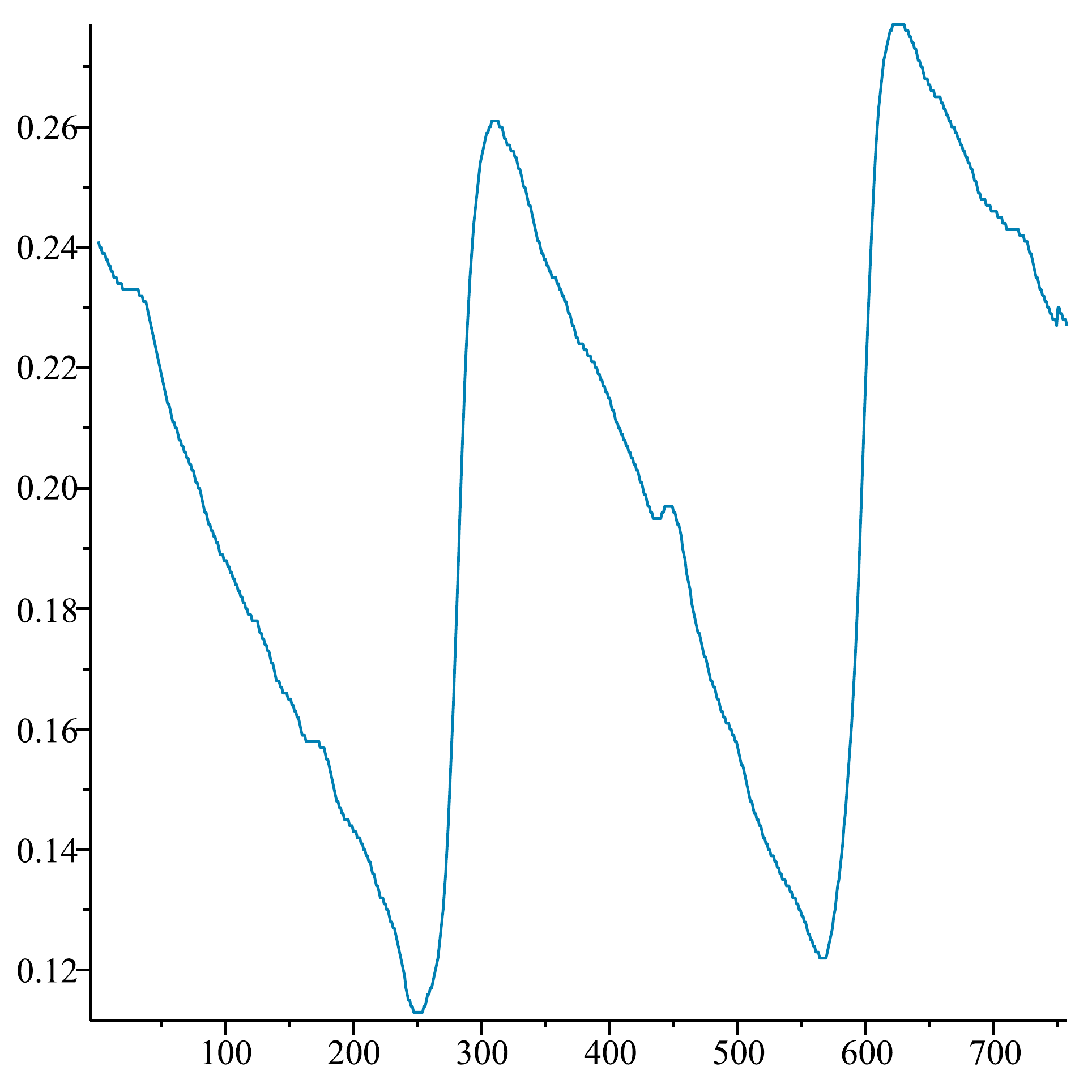}
 \includegraphics[width=0.48\textwidth]{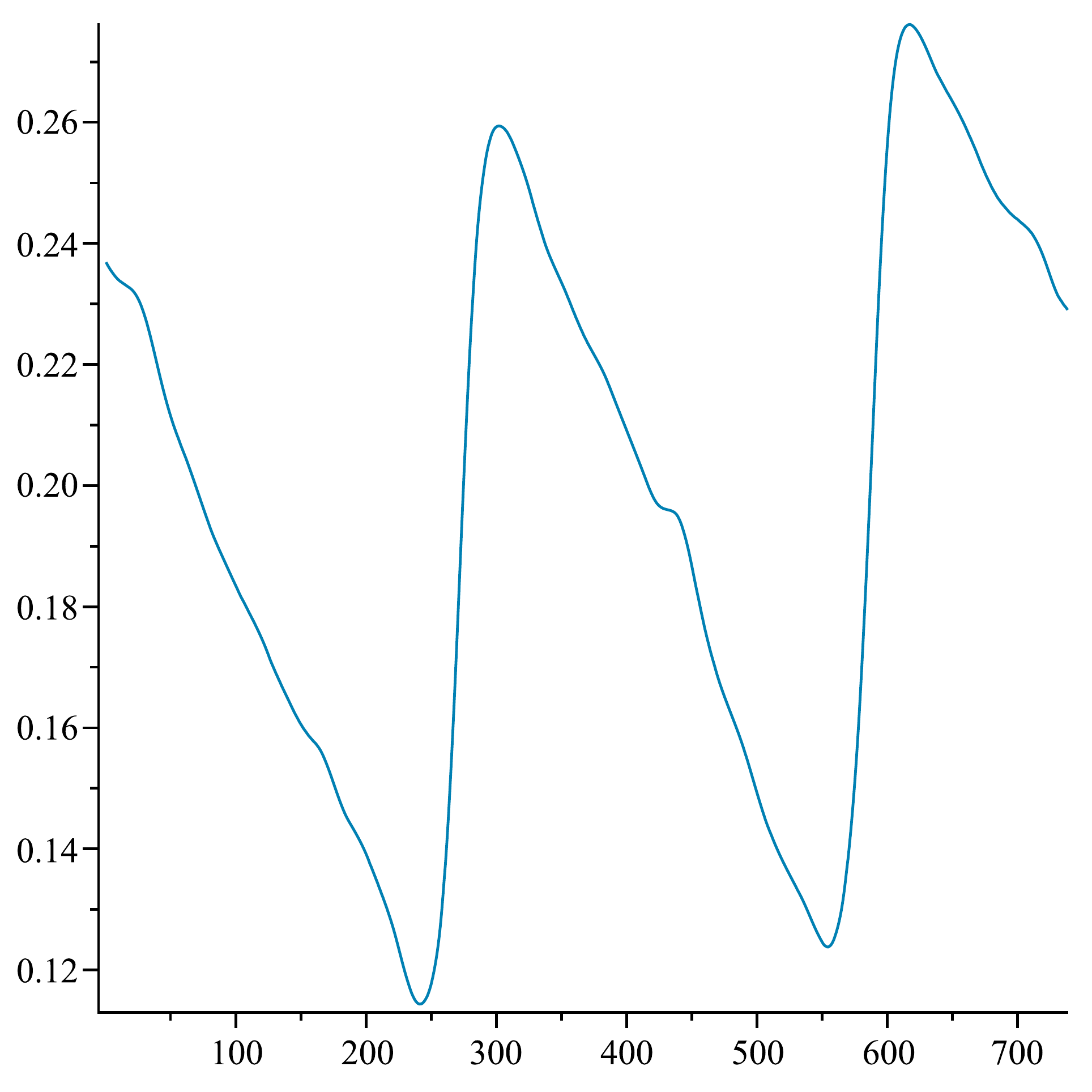}
\end{figure}
\begin{figure}[H]\caption{Gl\"attung mittels EA und rEA ($\alpha^{-1}=20$)}\label{fig12}\centering
 \includegraphics[width=0.48\textwidth]{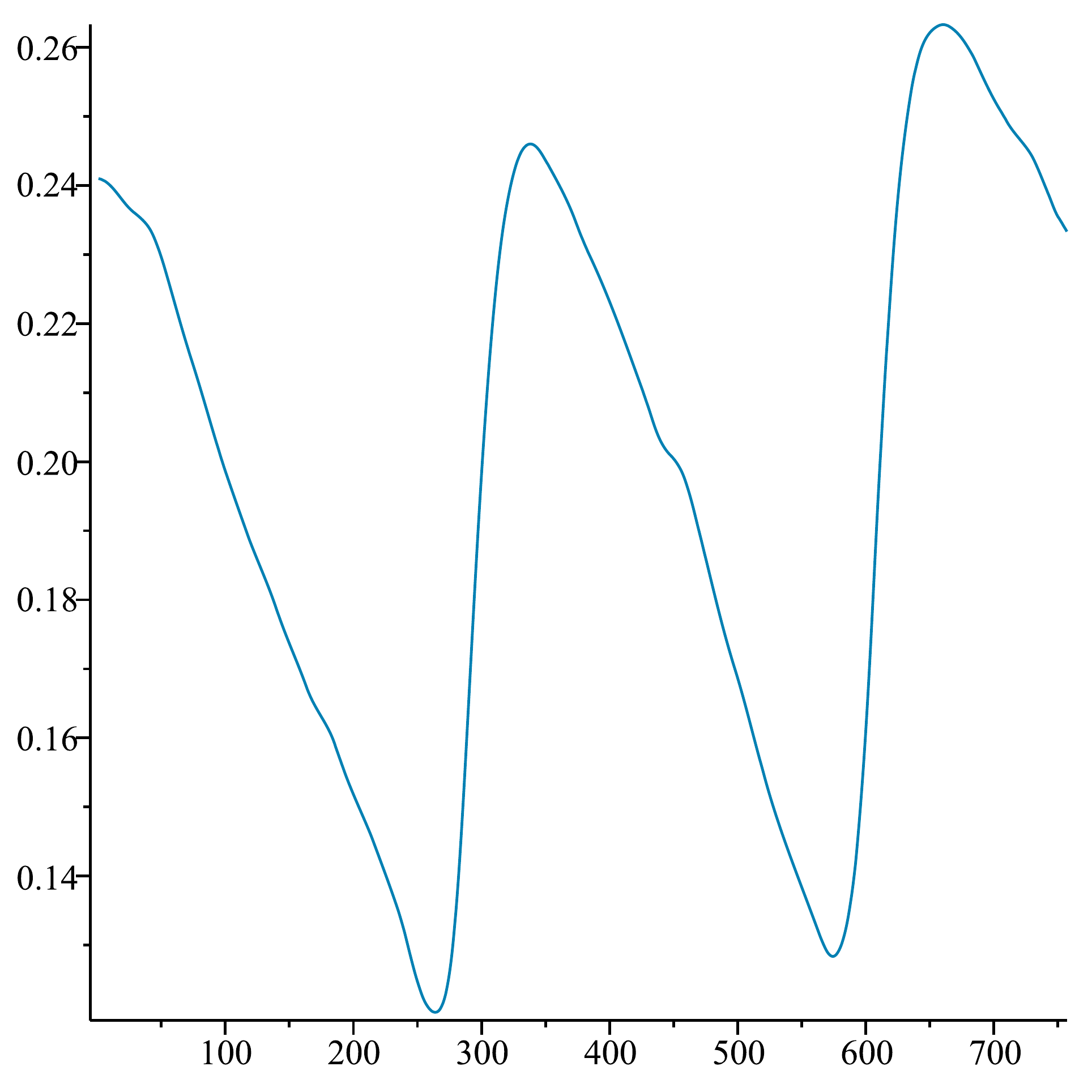}
 \includegraphics[width=0.48\textwidth]{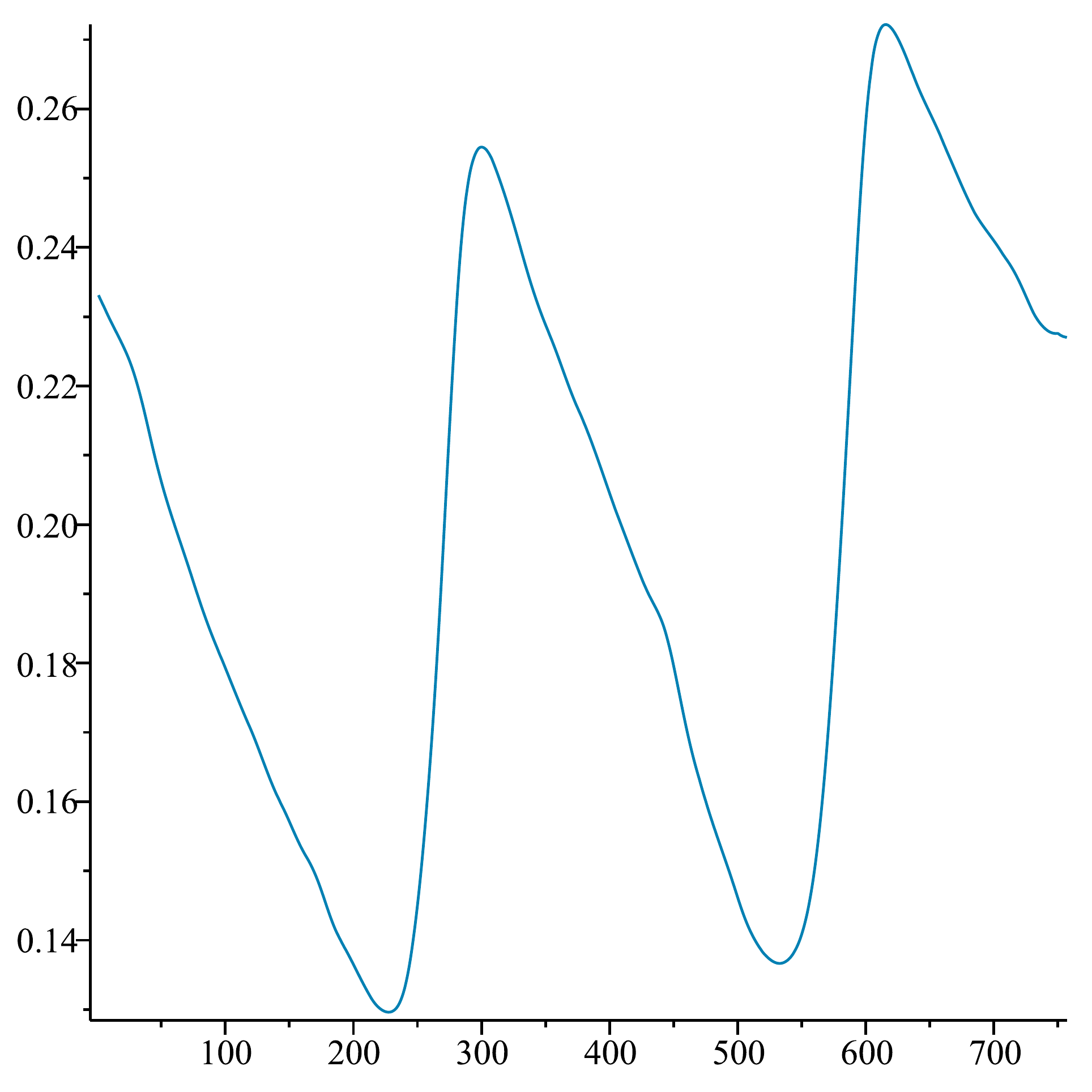}
\end{figure}
\begin{figure}[H]\caption{Gl\"attung mittels SEA ($\alpha^{-1}=20$)}\centering \label{fig13}
 \includegraphics[width=0.48\textwidth]{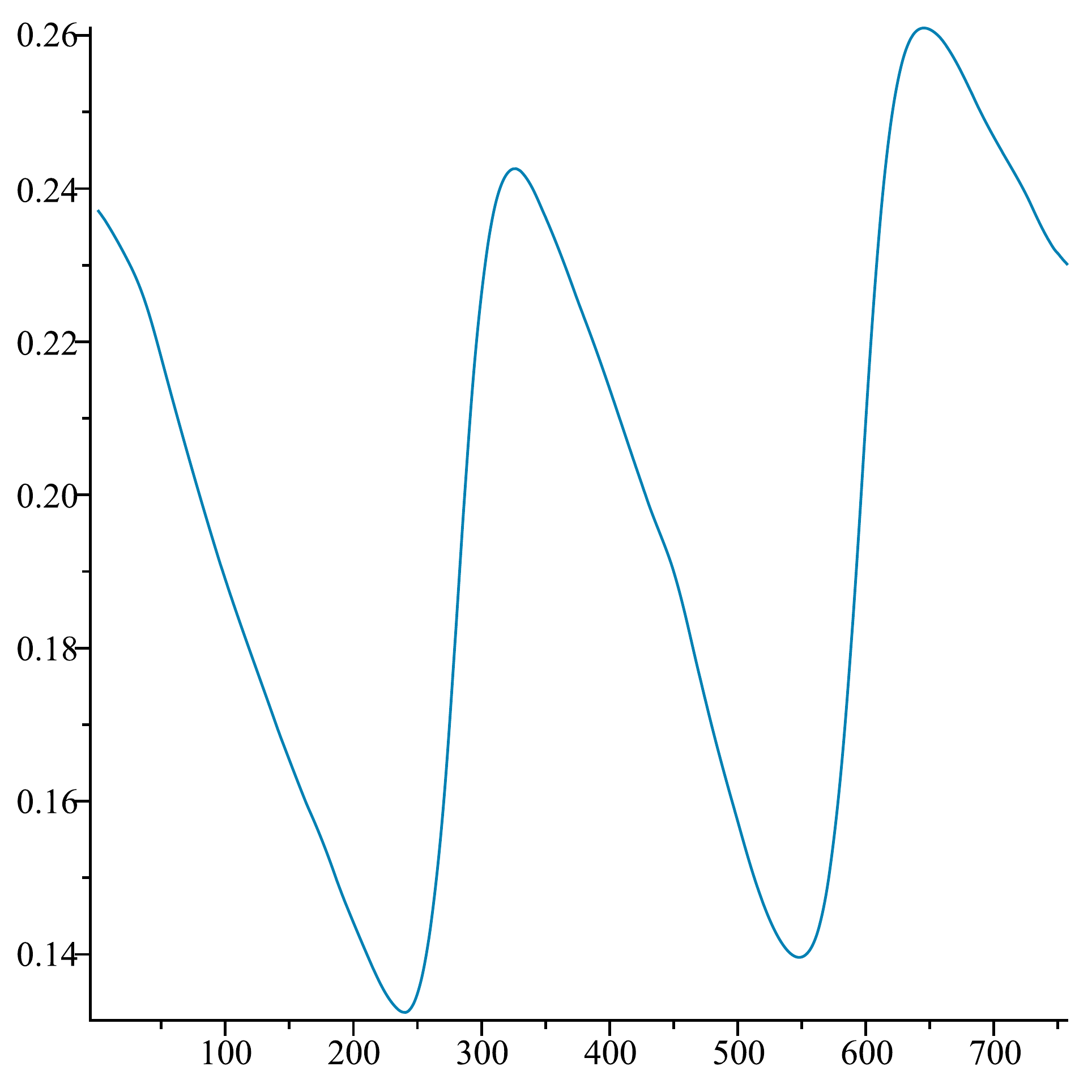}
\end{figure}

Alle Ans\"atze f\"uhren schon bei niedriger Gleitl\"ange (hier 20) zu guten bis sehr guten Ergebnissen und beseitigen die St\"orungen auf den rechten Flanken. Diese Gl\"attung ist bei SEA am besten ausgepr\"agt. 

Mit Beseitigung der St\"orung ist an dieser Stelle eine Gl\"attung zu verstehen, die im Sinne der Fragestellung keine weiteren als die gew\"unschten Extrema aufweist. Insofern, liefert auch die Gl\"attung mittels MA in Abbildung \ref{fig10} trotz der auff\"alligen St\"orung an der Flanke ein zufriedenstellendes Ergebnis.

In Abschnitt \ref{fazit} werden wir kurz auf die Wahl der Gleitl\"ange hinsichtlich der Ausgangsfragestellung eingehen.

\subsection{Daten mit stark destrukturierten Extrema}\label{sub2}

Durch eine starke Destrukturierung des Datensatzes insbesondere in den Extrema kommt es zu einer Vielzahl unerw\"unschter Maxima. Das Ziel ist auch hier die Ausgl\"attung dieser.

\begin{figure}[H]\caption{Originaldaten einer Laserausmessung mit Artefakten}\label{fig21}\centering
 \includegraphics[width=0.48\textwidth]{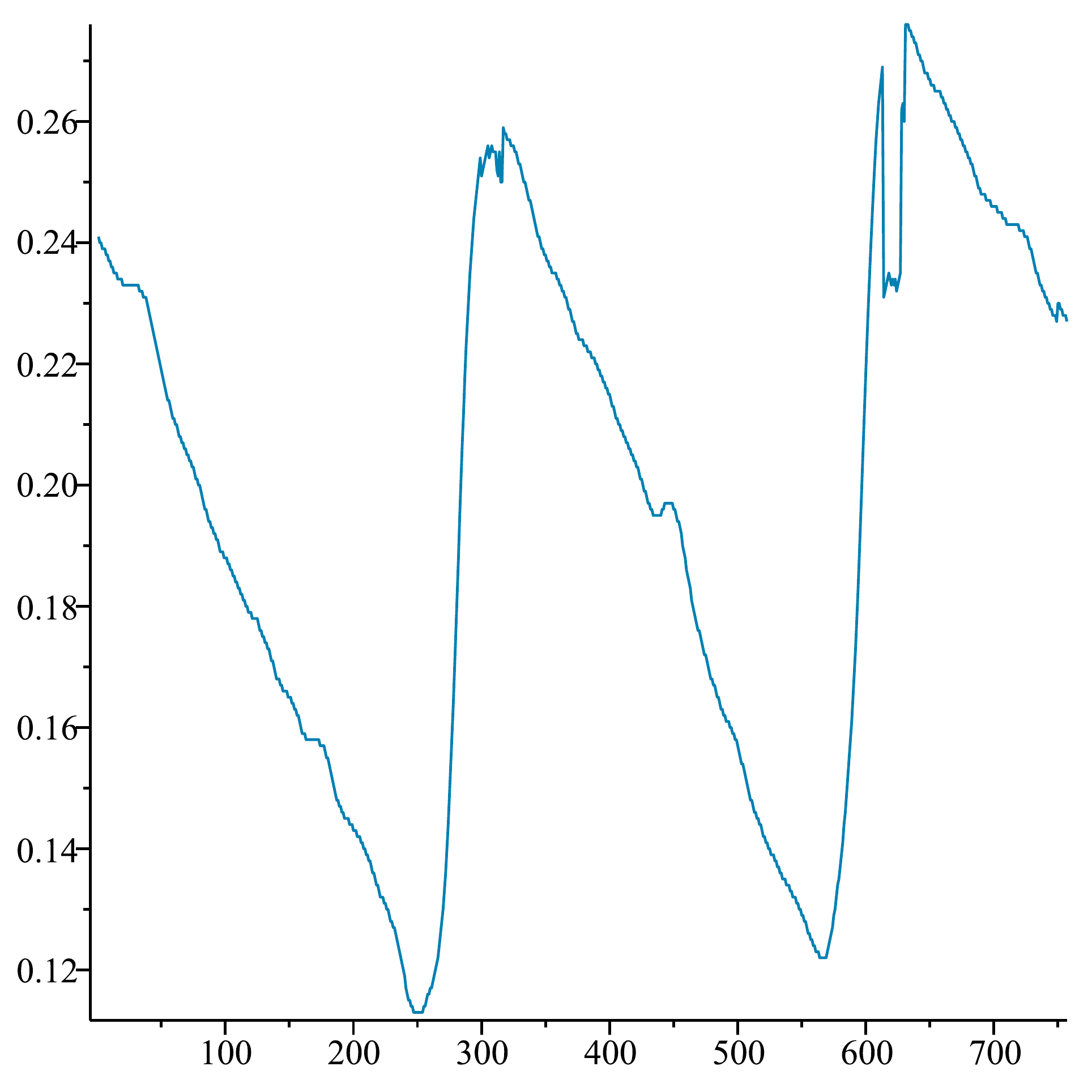}
\end{figure}
\begin{figure}[H]\caption{Gl\"attung mittels MA ($n+1=20$ und $n+1=90$)}\label{fig22}\centering
 \includegraphics[width=0.48\textwidth]{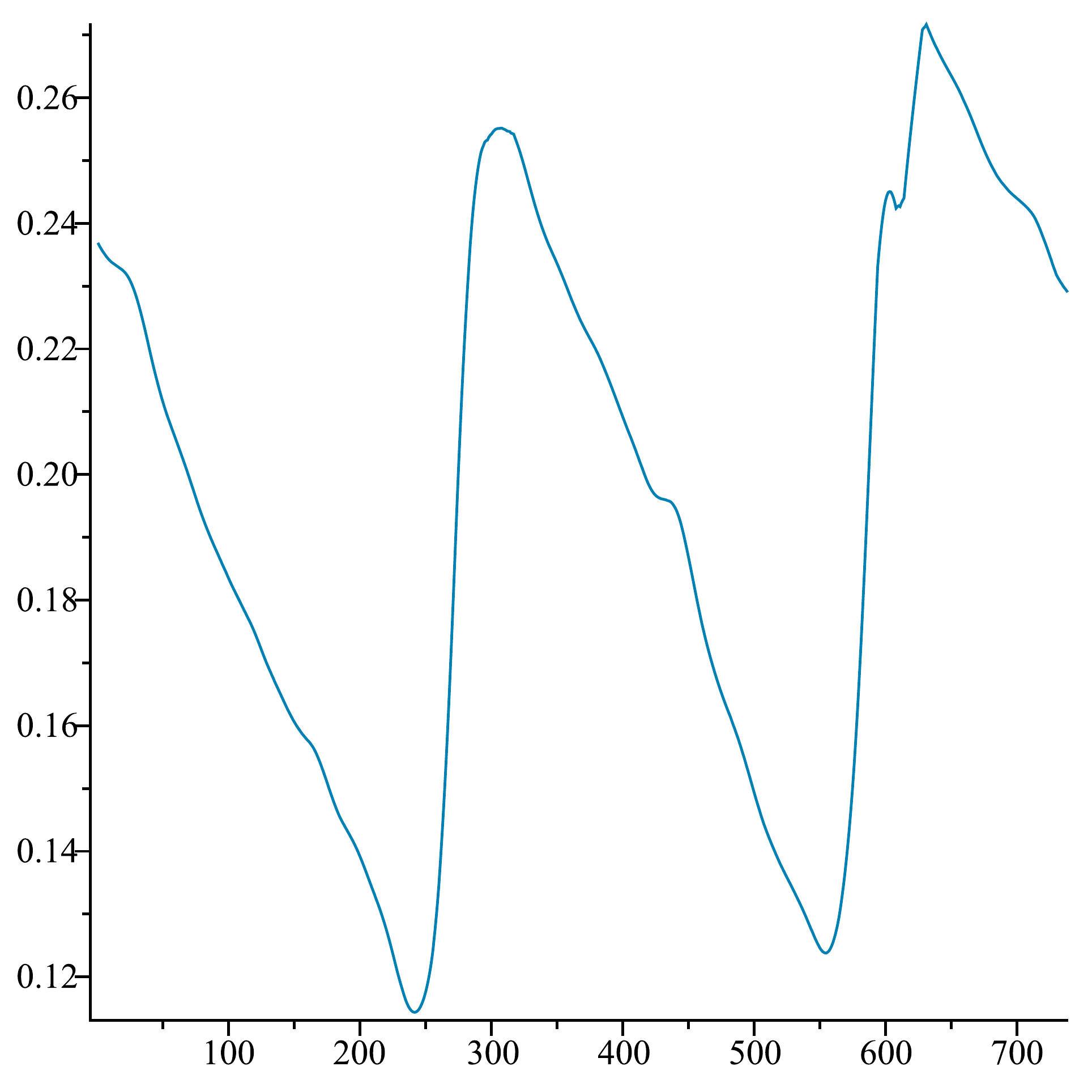}
 \includegraphics[width=0.48\textwidth]{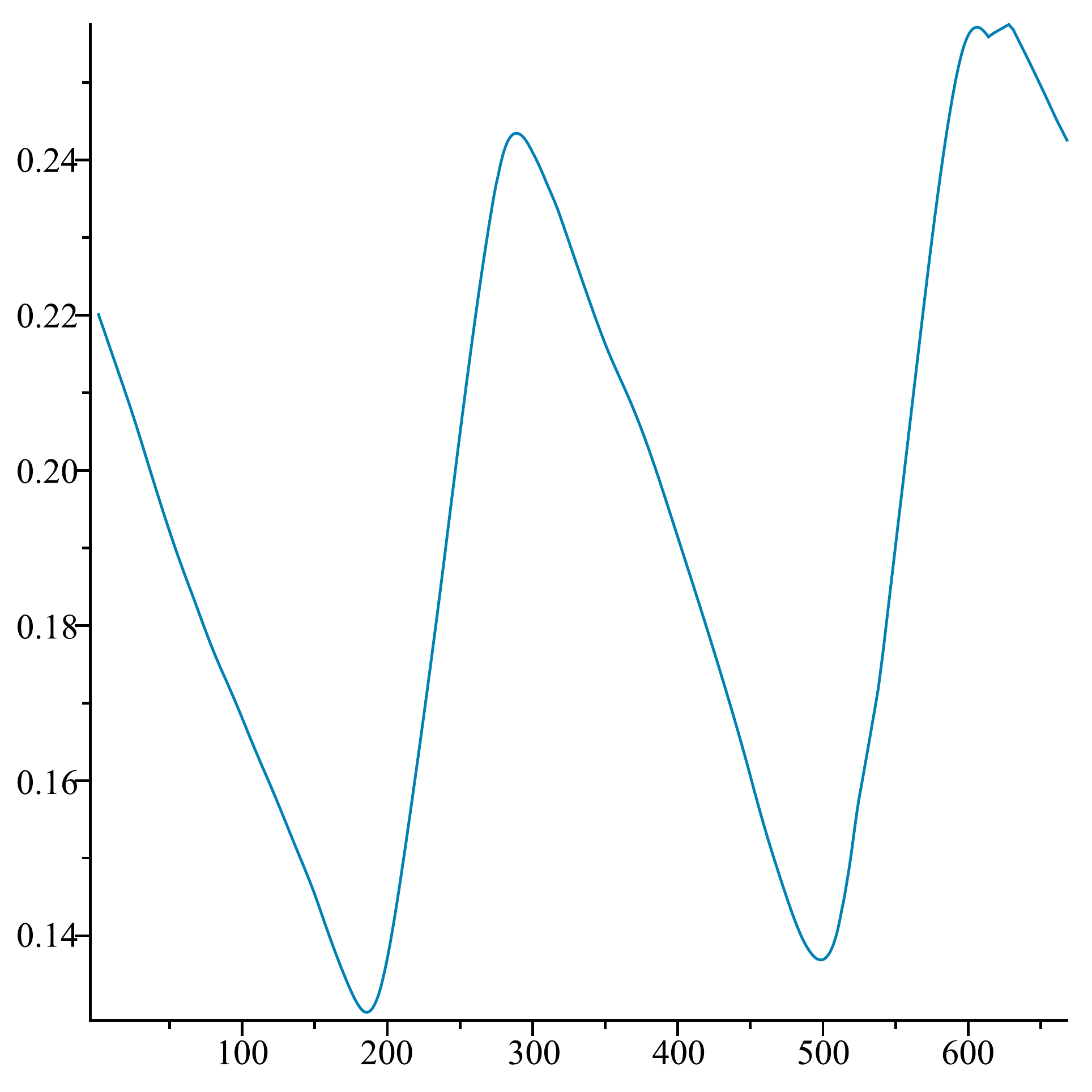}
\end{figure}
\begin{figure}[H]\caption{Gl\"attung mittels rEA ($\alpha^{-1}=20$ und $\alpha^{-1}=90$)}\label{fig23}\centering
 \includegraphics[width=0.48\textwidth]{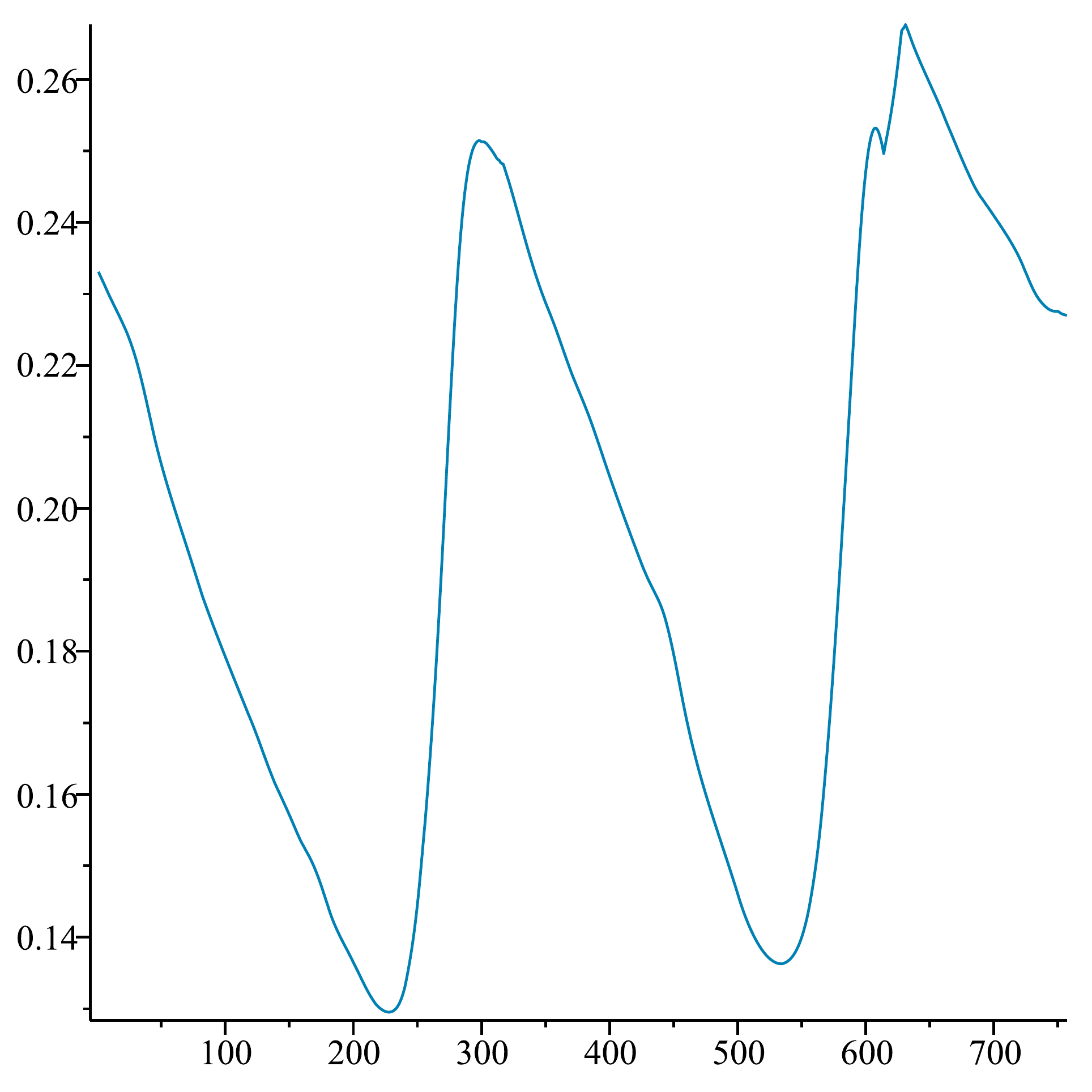}
 \includegraphics[width=0.48\textwidth]{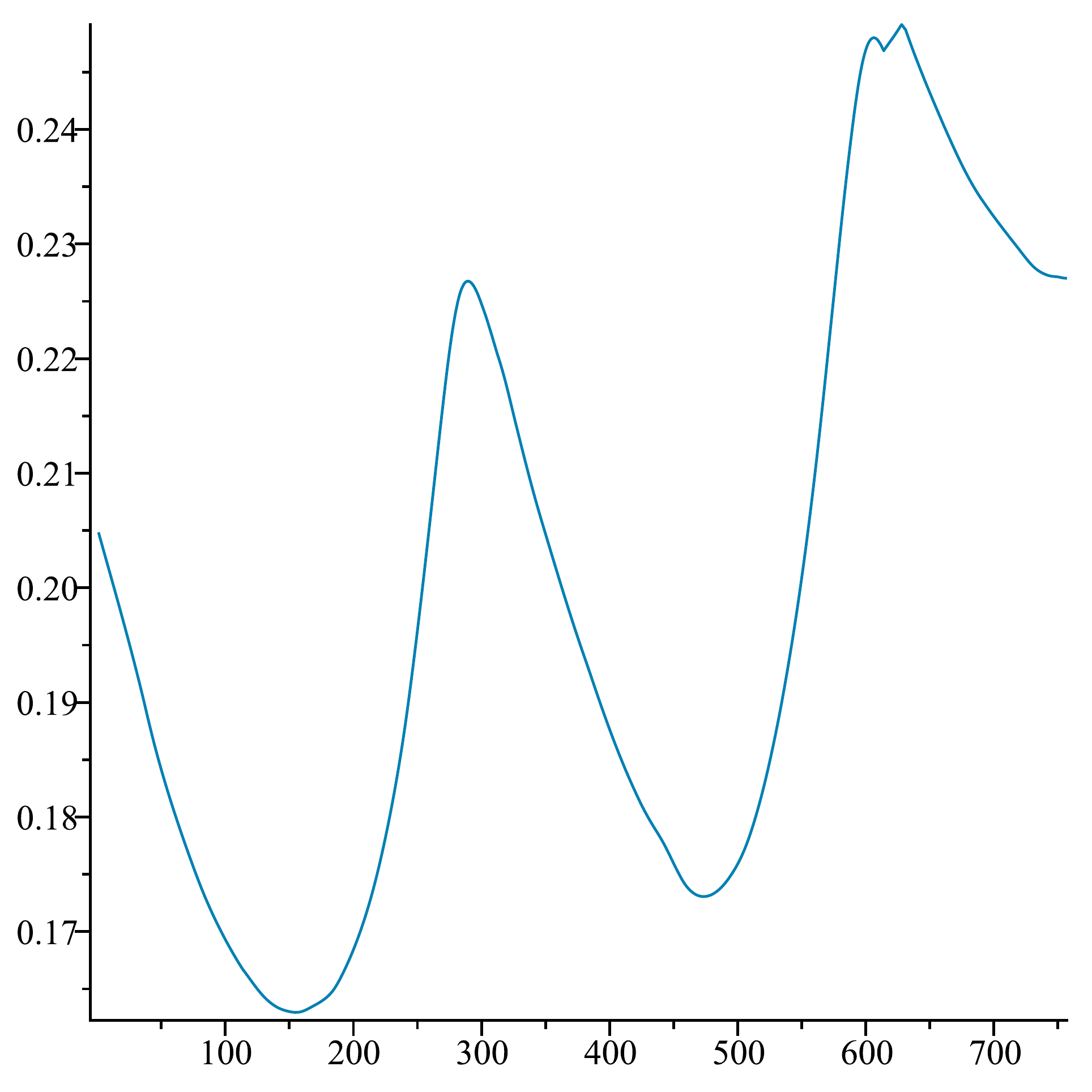}
\end{figure}
\begin{figure}[H]\caption{Gl\"attung mittels EA und SEA ($\alpha^{-1}=20$)} \label{fig24}\centering
 \includegraphics[width=0.48\textwidth]{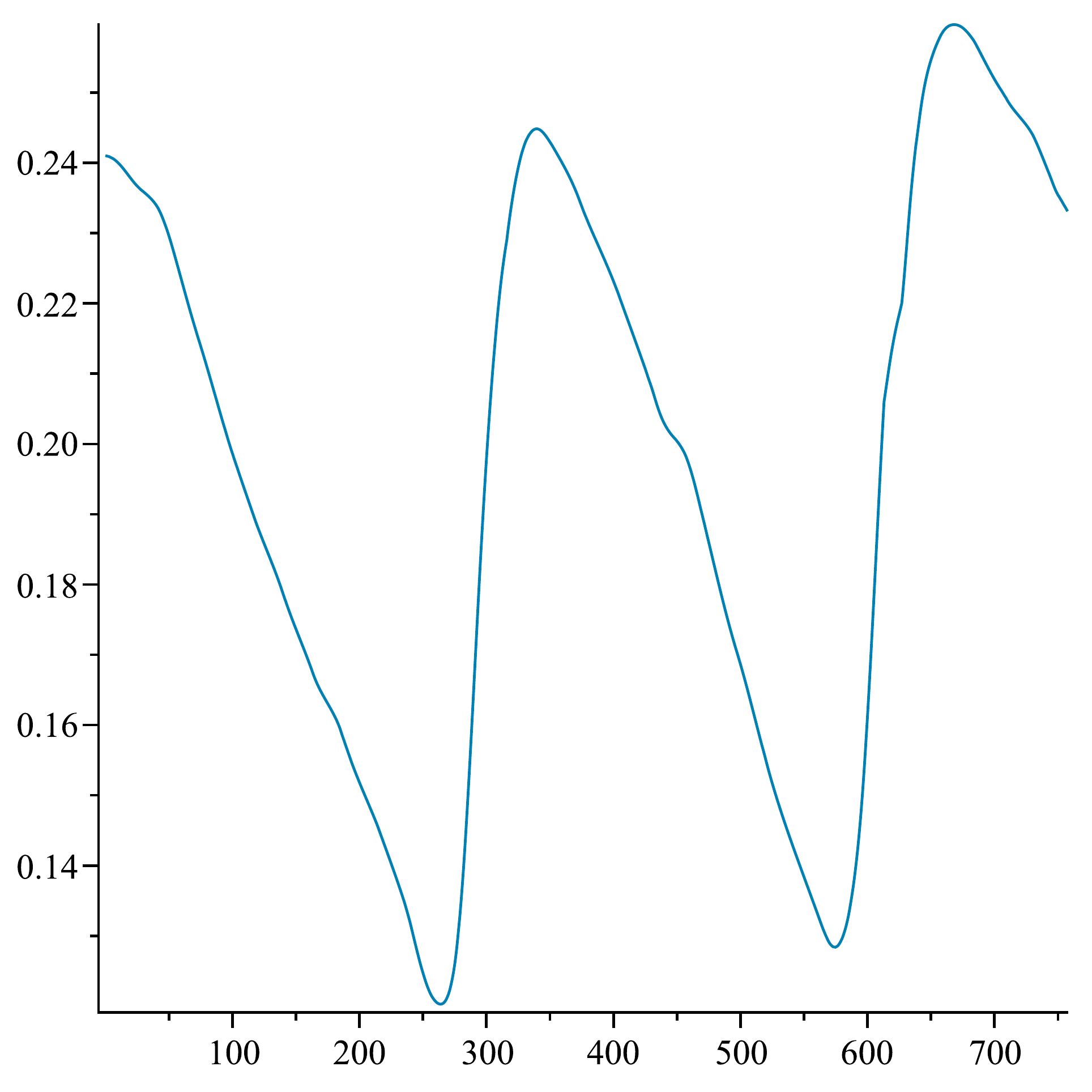}
 \includegraphics[width=0.48\textwidth]{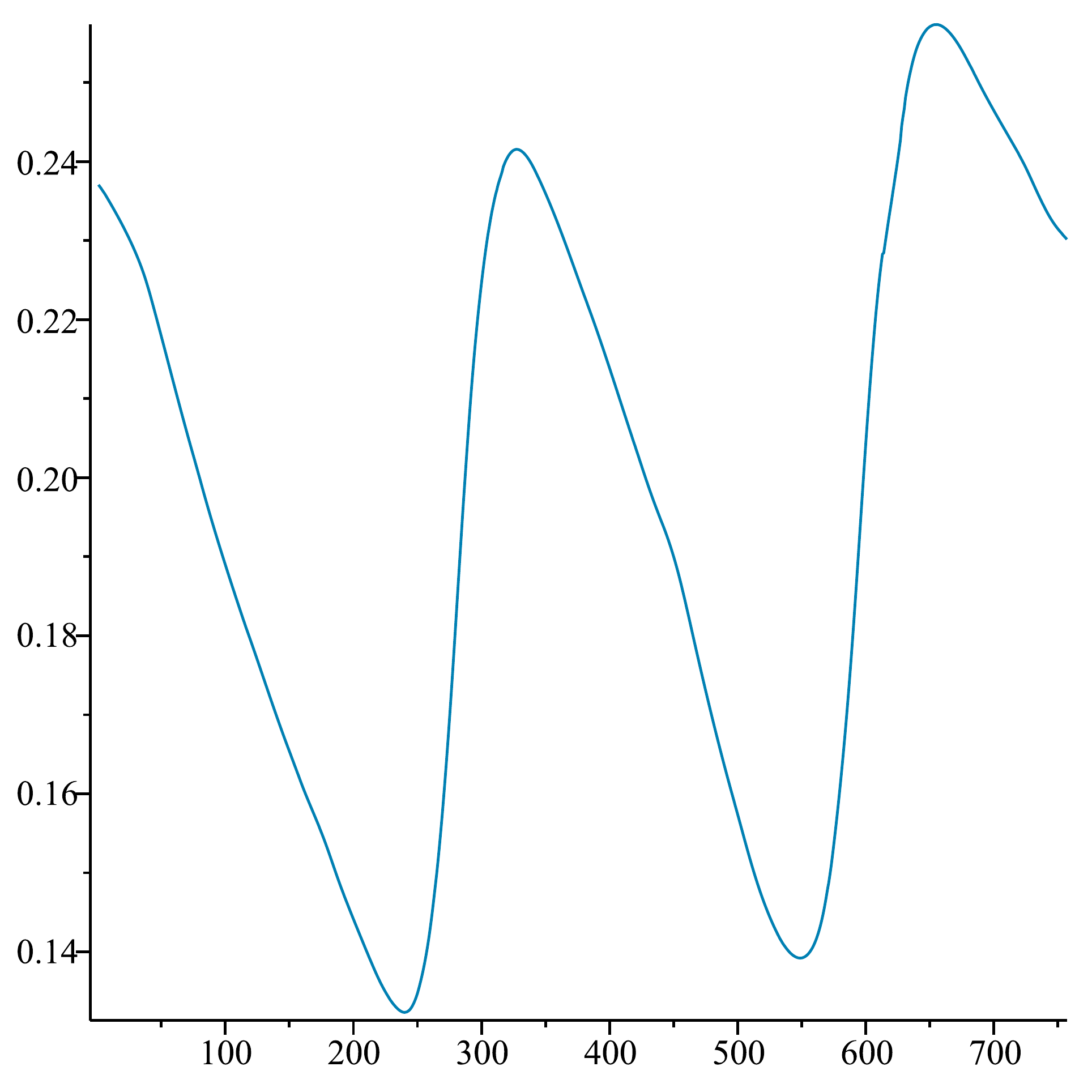}
\end{figure}

In den Abbildungen \ref{fig21} und \ref{fig23} sieht man, dass MA und rEA es nicht schaffen, die starken Artefakte rechts in der Grafik zu gl\"atten. Das Ergebnis kann auch mit steigender Gleitl\"ange nicht verbessert werden. Beim rEA liegt das an der Tatsache, dass die St\"orungen sehr fr\"uh auftreten und daher nur wenig Werte zur Berechnung der Gl\"attung herangezogen werden. Dieser Effekt wird dadurch, dass die St\"orungen  recht stark ausfallen, noch verst\"arkt. Beim MA liegt es an der Tatsache, dass alle Werte mit dem gleichen Gewicht auftreten und deshalb im Allgemeinen sehr starke Schwankungen in der Messreihe erst bei hoher Gleitl\"ange aufgefangen werden k\"onnen.

EA und insbesondere SEA liefern hier auch bei kleiner Gleitl\"ange gute bis sehr gute Ergebnisse. Bei EA liegt es daran, dass die gro\ss en Artefakte sehr sp\"at auftreten und deshalb sehr viele Werte in die Berechnung mit einfliessen (hier ca $80\%$). Bei der Berechnung des SEA spielt die Stelle des Auftretens der Artefakte keine Rolle, da der Mangel des rEA durch die gleichzeite Verwendung des EA aufgehoben wird.

\subsection{\"Uberlagerung zweier periodischer Schwingungen}\label{sub3}

Als erstes k\"unstliches Beispiel zum Testen der Gl\"attungsverfahren betrachten wir eine \"Uberlagerung zweier periodischer Schwingungen mit unterschiedlichen Frequenzen. Dabei ist die Frequenzdifferenz hoch gew\"ahlt und die Amplitude der hochfrequenten St\"orung klein im Vergleich zur Amplitude der niederfrequenten Grundschwingung.

\begin{figure}[H]\caption{Originaldaten der \"Uberlagerung}\label{fig31}
\centering \includegraphics[width=0.48\textwidth]{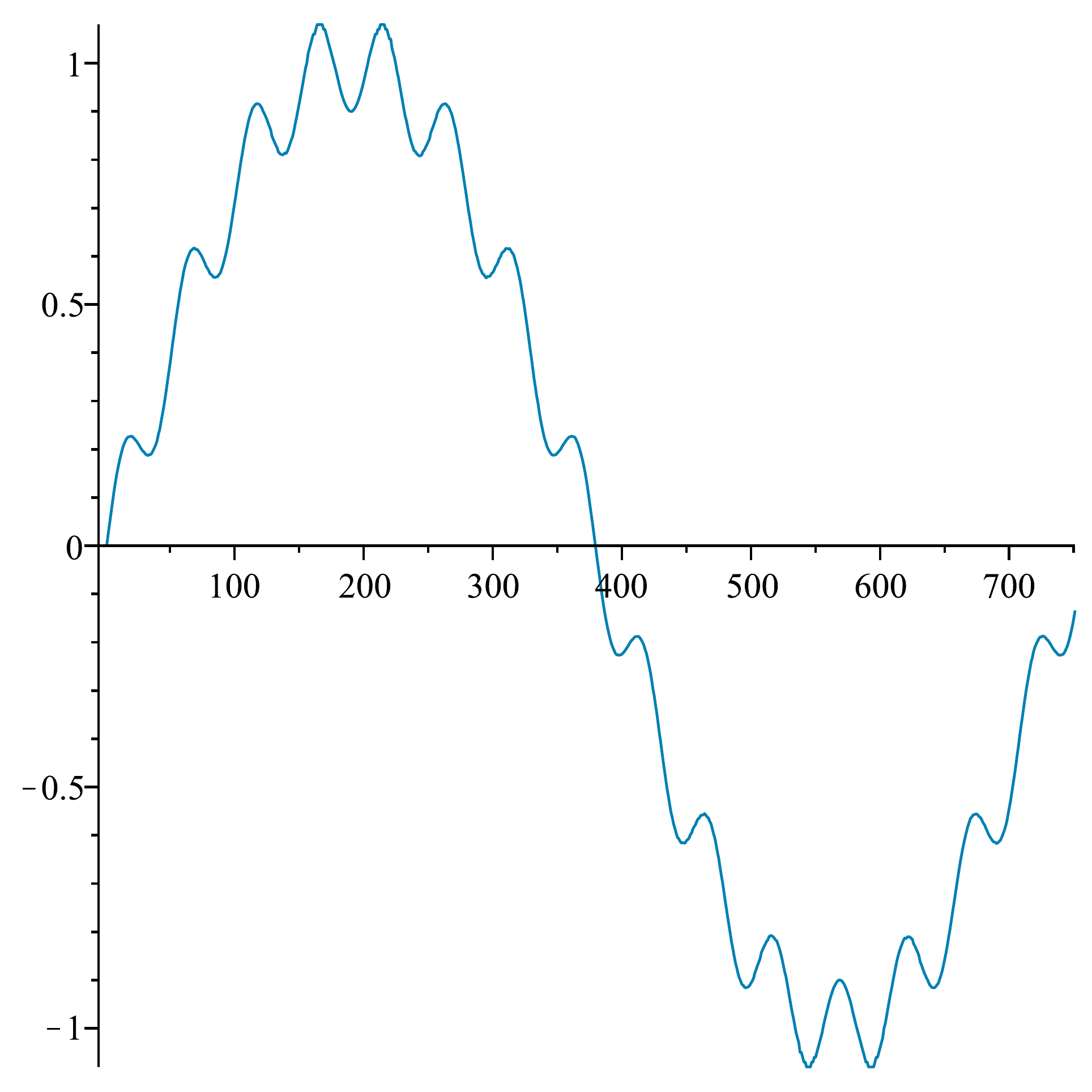}
\end{figure}
\begin{figure}[H]\caption{Gl\"attung mittels MA ($n+1=40$ und $n+1=100$)}\label{fig32}\centering
 \includegraphics[width=0.48\textwidth]{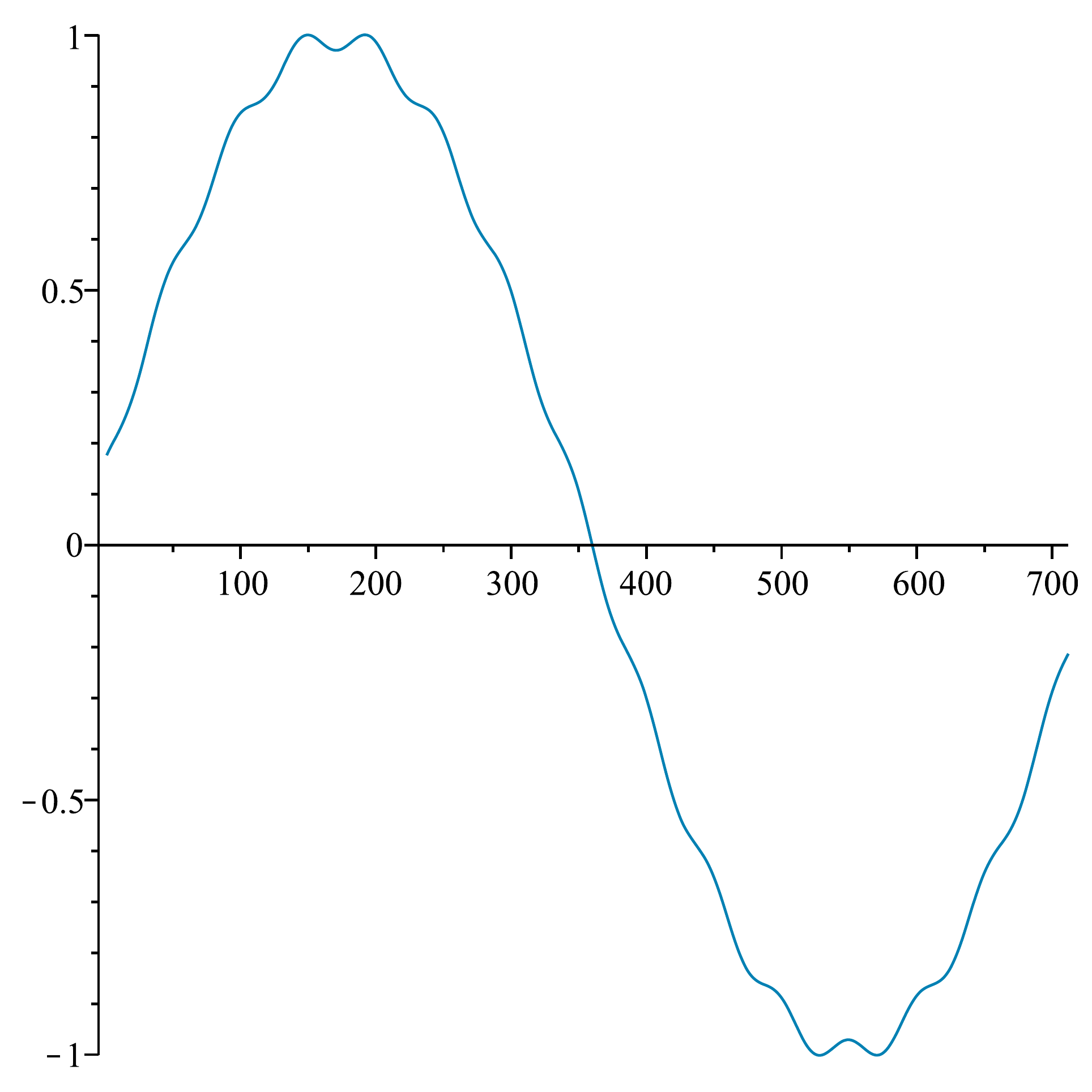}
 \includegraphics[width=0.48\textwidth]{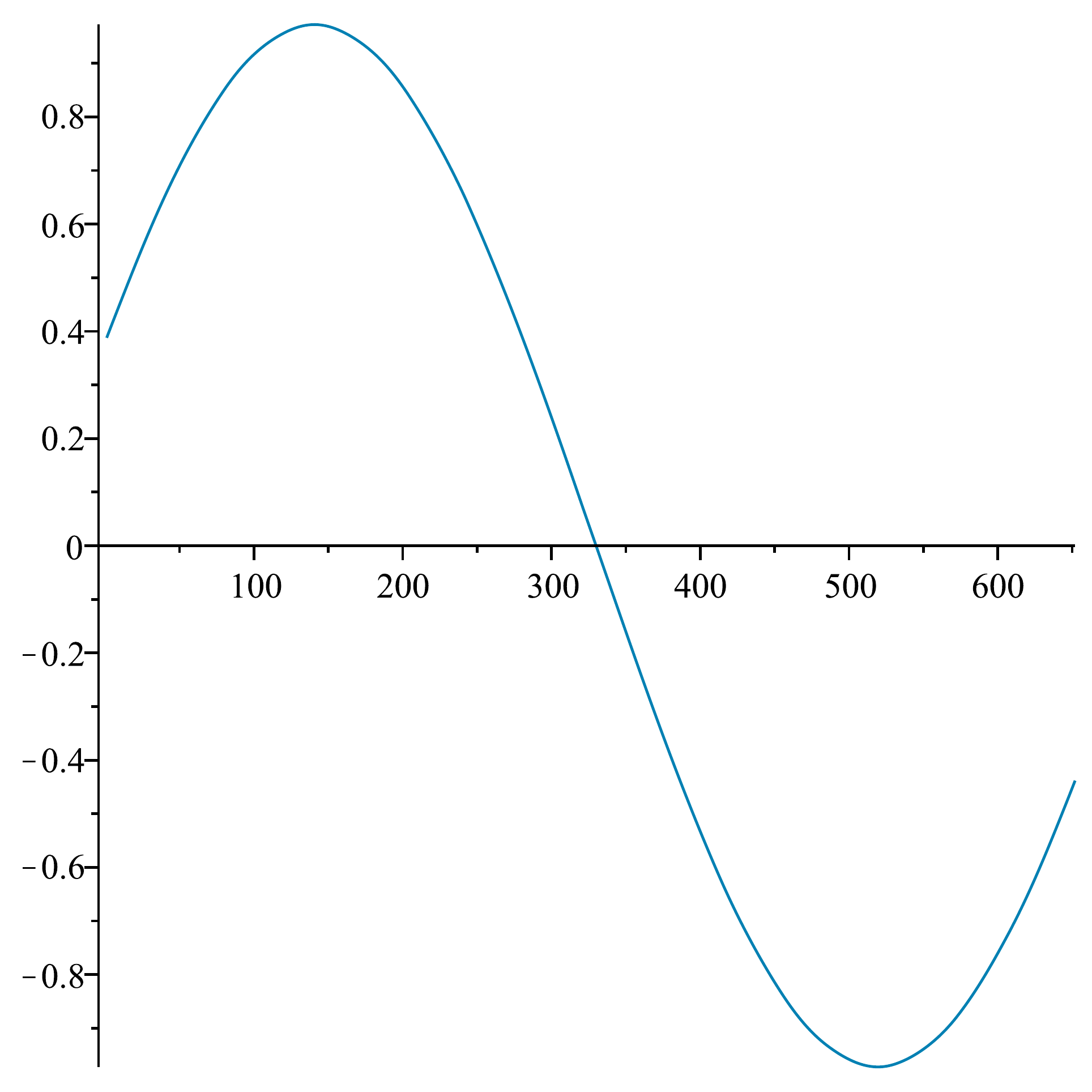}
\end{figure}
\begin{figure}[H]\caption{Gl\"attung mittels EA und rEA ($\alpha^{-1}=40$)}\label{fig33}\centering
 \includegraphics[width=0.48\textwidth]{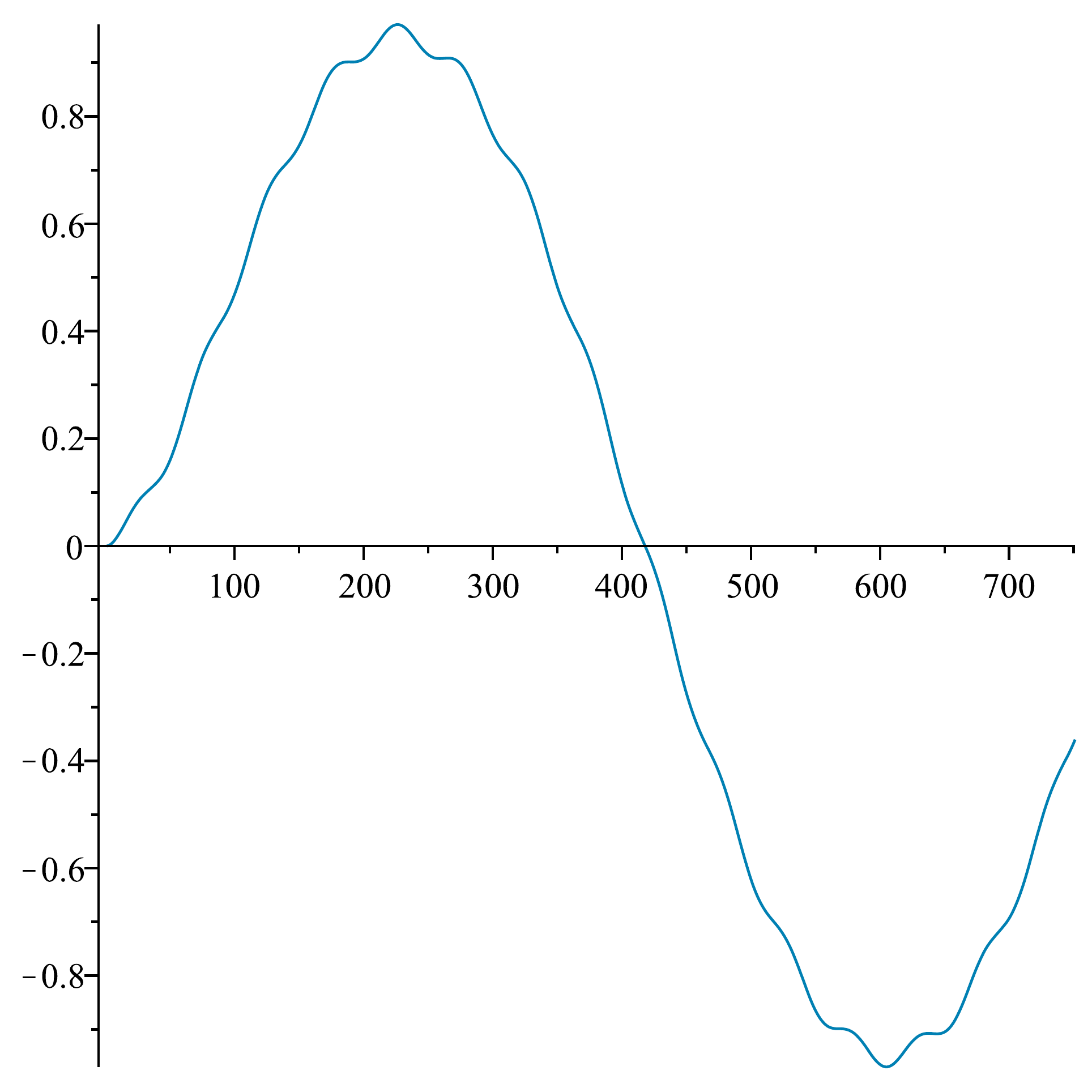}
 \includegraphics[width=0.48\textwidth]{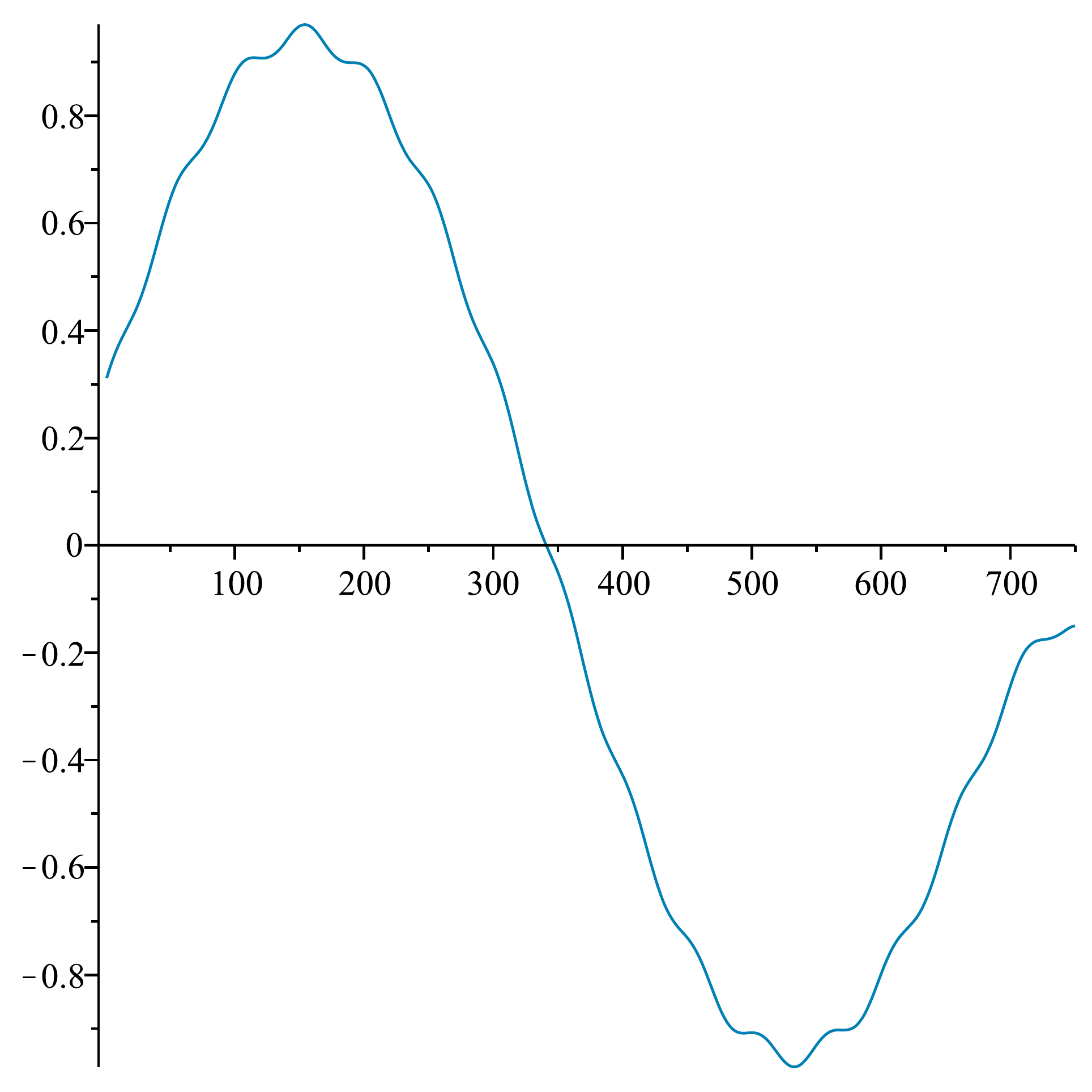}
\end{figure}
\begin{figure}[H]\caption{Gl\"attung mittels SEA ($\alpha^{-1}=40$)} \label{fig34}\centering
 \includegraphics[width=0.48\textwidth]{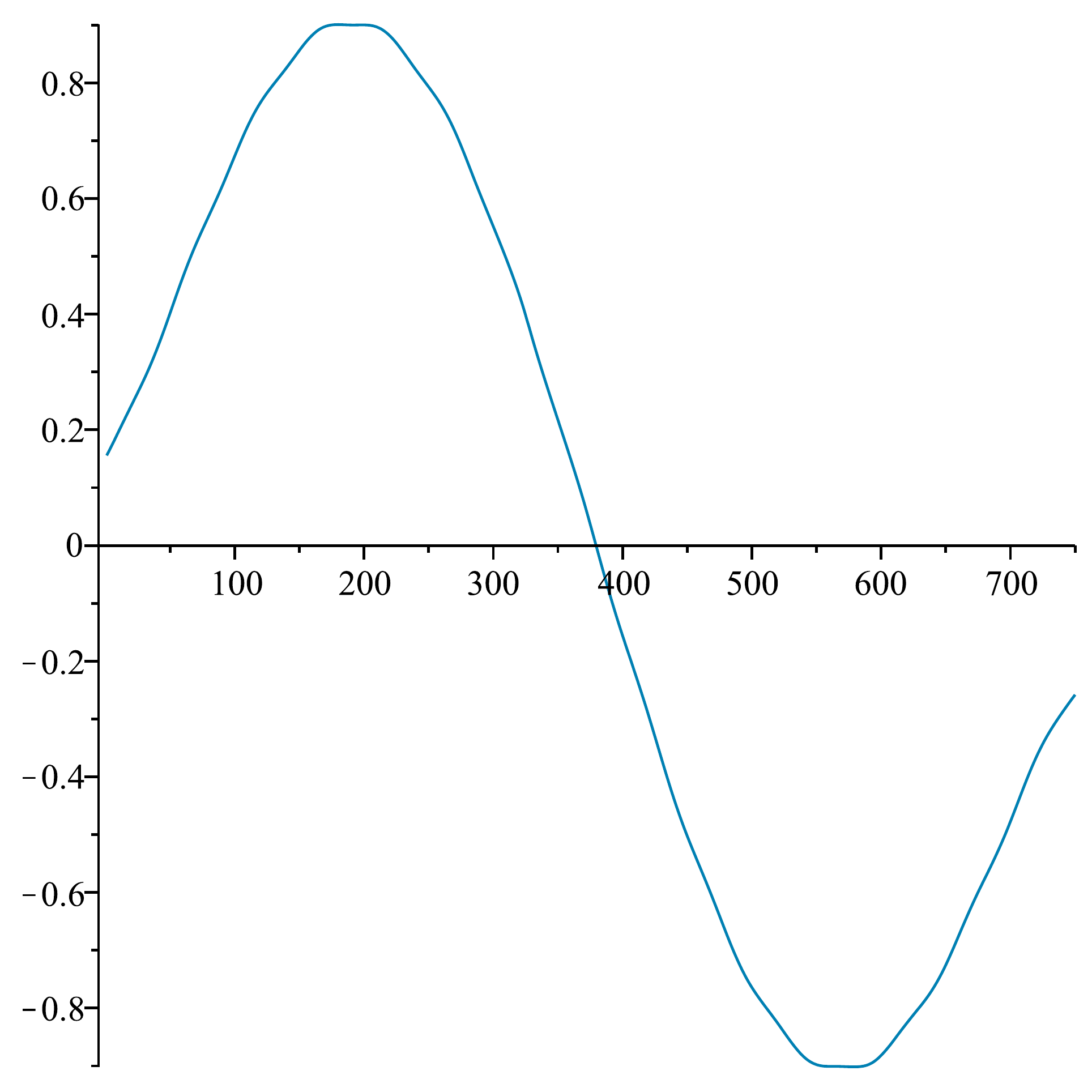}
\end{figure}

Die Gl\"attungsverfahren MA, EA und rEA haben in dem Beispiel der \"uberlagerten Schwingung Probleme. MA schafft es  mit steigender Gleitl\"ange jedoch eher -- wenn auch nicht viel eher -- als rEA und EA die Oberschwingungen auszugl\"attet. In diesem Beispiel ist die gleiche Gewichtung aller Beitr\"age zur Gl\"attung ein Vorteil des MA gegen\"uber dem EA.
Wieder liefert auch in diesem Beispiel der SEA das beste Ergebnis, denn schon bei vergleichbarer mittlerer Gleitl\"ange von 40 werden die Oberschwingungen beseitigt.

Wir m\"ochten an dieser Stelle betonen, dass  der Charakter der Kurve als \"Uberlagerung von Schwingungen  (Abb.~\ref{fig31}) beim \"Ubergang zur Gl\"attung (Abb.~\ref{fig34}) selbstverst\"andlich ver\"andert wird. Im Hinblick auf unsere Fragestellung wird der Charakter der Kurve jedoch allein durch die Grundschwingung bestimmt, und die Oberschwingung wird als St\"orung angesehen. Insofern erh\"alt auch hier die Gl\"attung den Charakter der Ausgangskurve.

Diese Ver\"anderung des Schwingungscharakters, also das Ausblenden der Oberschwingungen, bedeutet, dass die Gl\"attungsverfahren im Fall hoher Gleitl\"angen als Tiefpassfilter wirken. Zur Diskussion des MA und des EA als Tiefpassfilter und zu Eigenschaften der zugeh\"origen Filterfunktionen siehe zum Beispiel ~\cite{KN}. Ein Tiefpassfilter ruft in der Regel eine Phasenverschiebung der Grundschwingug hervor. Dies kann man hier gut durch das Verschieben der Extrema nach rechts bzw.\ links bei der Gl\"attung mittels EA  bzw.\ rEA (Abb.~\ref{fig33}) erkennen. Vergleicht man die Abbildungen \ref{fig31} und \ref{fig34} mit \ref{fig33} so sieht man dass die Symmetrisierung bei der Anwendung des SEA die zwei Phasenverschiebungen kompensiert und somit keine nennenswerte Verschiebung der Extrema auftritt.

 \subsection{Rechtecksignal mit endlicher Fourierentwicklung}\label{sub4}

\begin{figure}[H]\caption{Originaldaten der Fourierentwicklung}\label{fig41}\centering
 \includegraphics[width=0.48\textwidth]{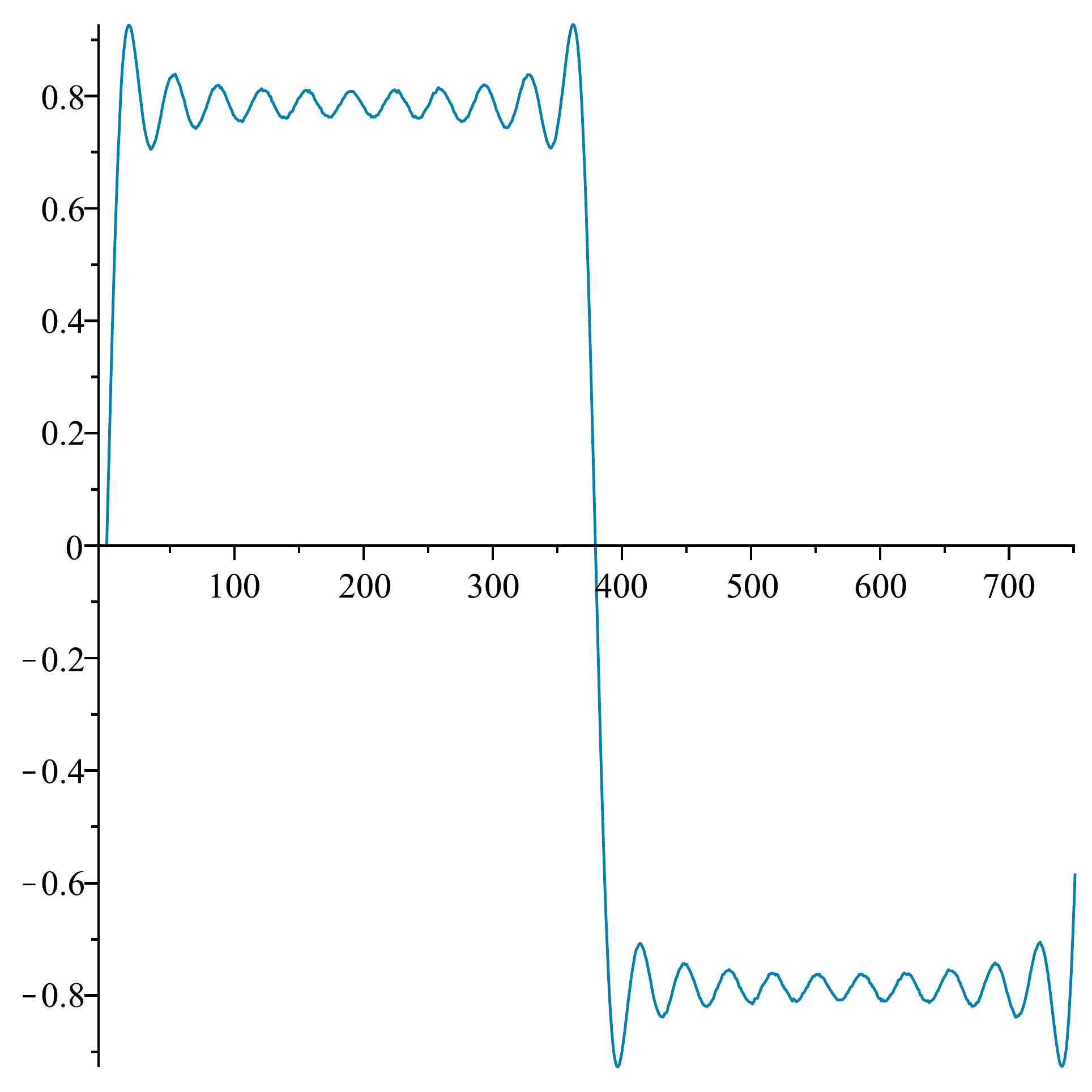}
\end{figure}

Wir betrachten als weiteres Beispiel die Fourierentwicklung eines Rechtecksignals -- genauer die Teilsumme zehnter Ordnung der Fourierentwicklung. Deutlich erkennt man in Abbildung \ref{fig41} das so genannte Gibbssche Ph\"anomen, das durch die hochfrequente \"Uberschwingungen an den Spr\"ungen des Signals charakterisiert ist. Dieses Ph\"anomen der \"Uberschwingung tritt unabh\"angig davon auf, wie sp\"at die Fourierreihe abgebrochen wird, und ihre Amplitude h\"angt von der Amplitude des Rechtecksignals ab. Das Ph\"anomen tritt nicht nur bei Rechecksignalen auf, sondern ist grunds\"atzlich zu finden, wenn die Ausgangsfunktion einen Sprung aufweist. Dabei nimmt die Amplitude der \"Uberschwingung  zwar mit steigender Ordnung der Fourierentwicklung ab, unterschreitet jedoch niemals ca.\ $9\%$ der Sprungh\"ohe. Vergleiche dazu den sch\"onen \"Ubersichtsartikel \cite{HH} oder den klassischen Artikel \cite{Gr} oder die Ausf\"uhrungen in \cite{How}.

\begin{figure}[H]\caption{Gl\"attung mittels MA ($n+1=40$ und $n+1=100$)}\label{fig42}\centering
 \includegraphics[width=0.48\textwidth]{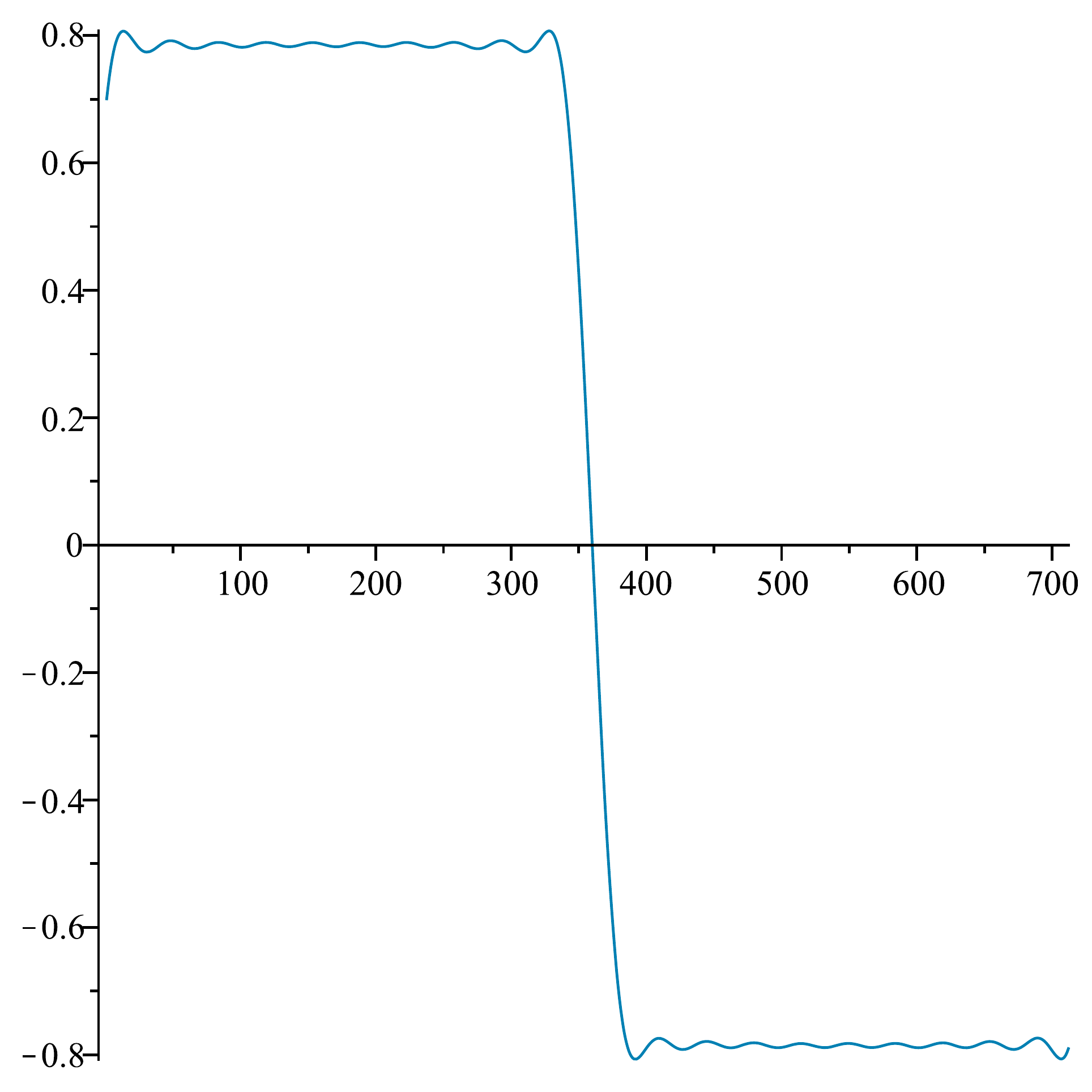}
 \includegraphics[width=0.48\textwidth]{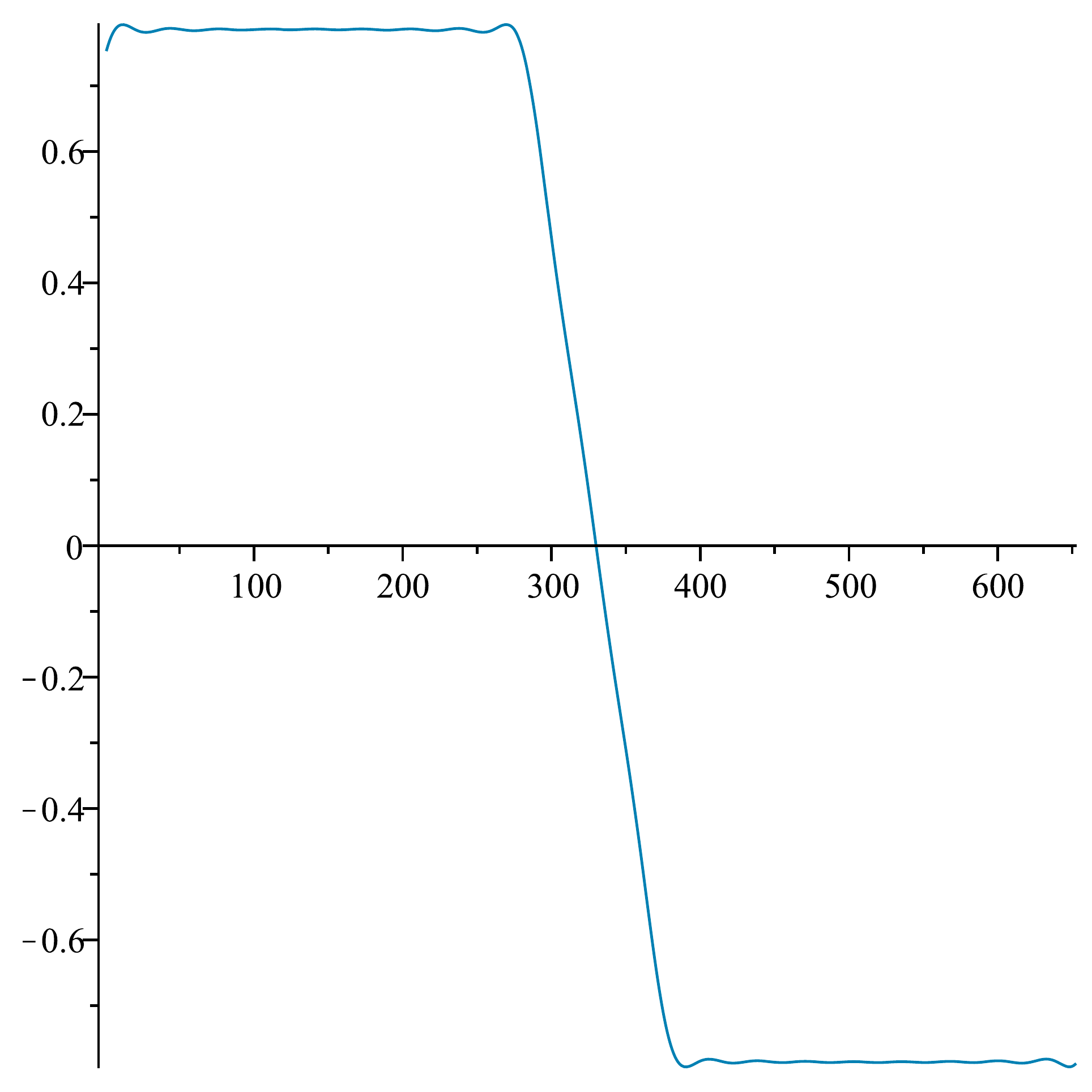}
\end{figure}
\begin{figure}[H]\caption{Gl\"attung mittels EA und rEA ($\alpha^{-1}=40$)}\label{fig43}\centering
 \includegraphics[width=0.48\textwidth]{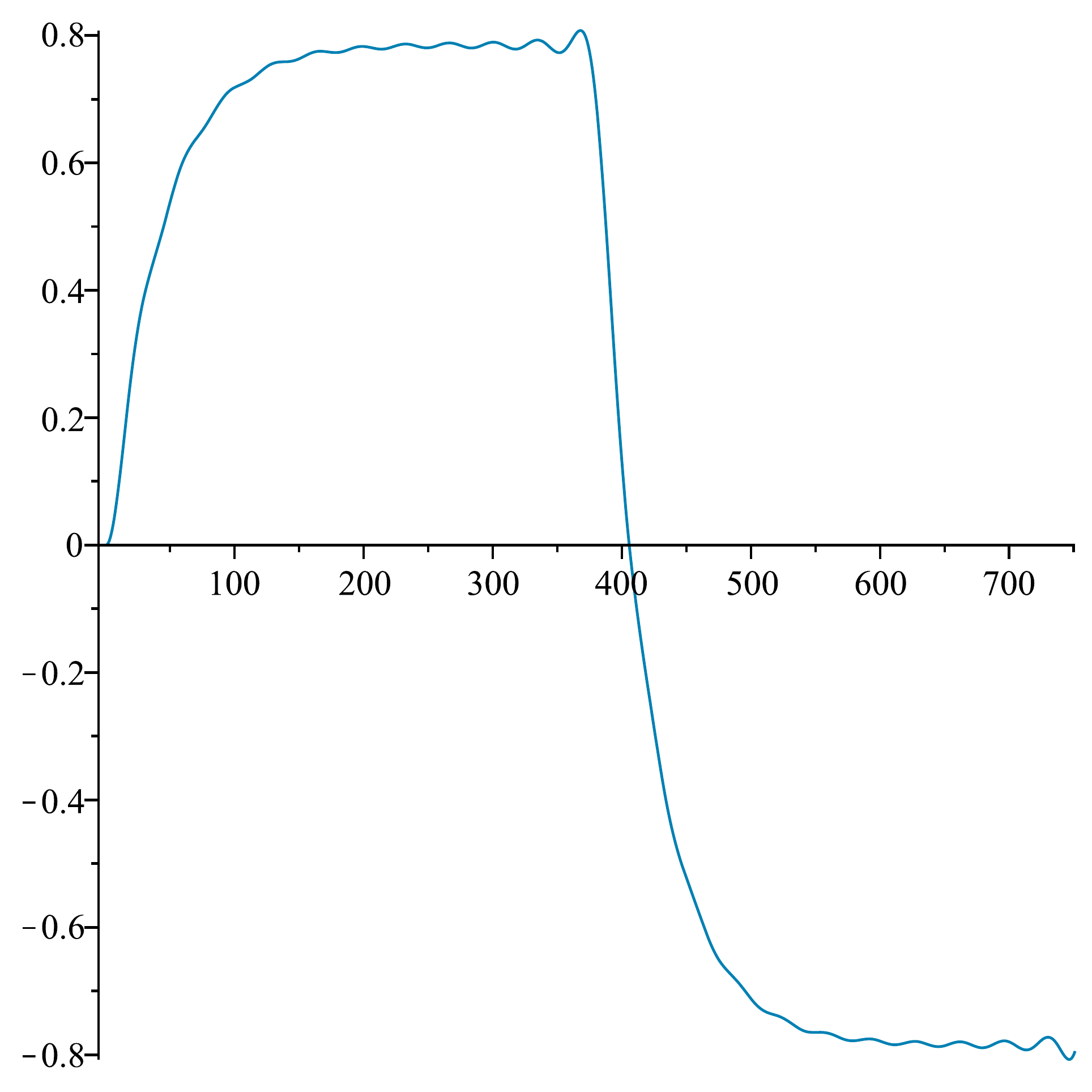}
 \includegraphics[width=0.48\textwidth]{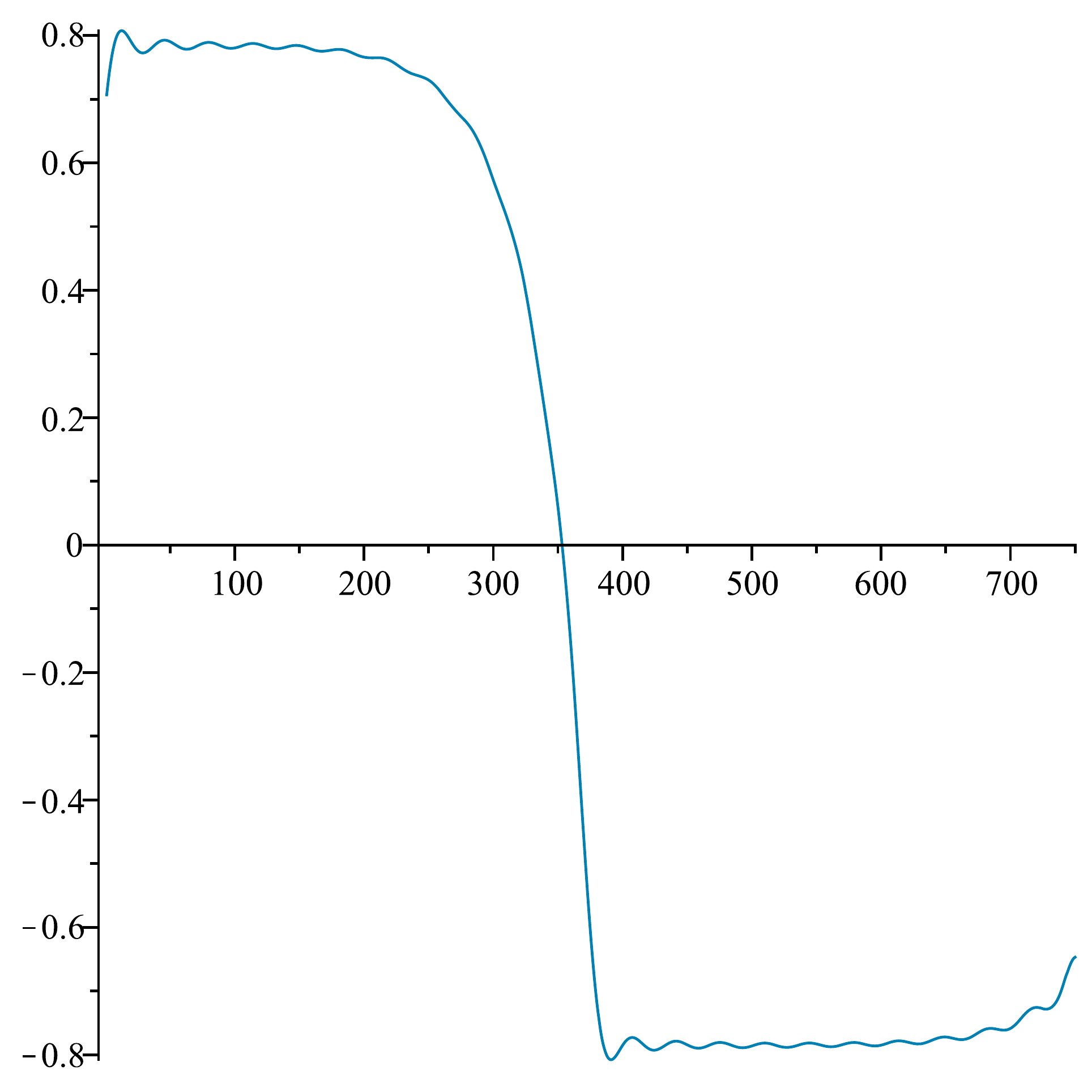}
\end{figure}
\begin{figure}[H]\caption{Gl\"attung mittels SEA ($\alpha^{-1}=40$)} \label{fig44}\centering
 \includegraphics[width=0.48\textwidth]{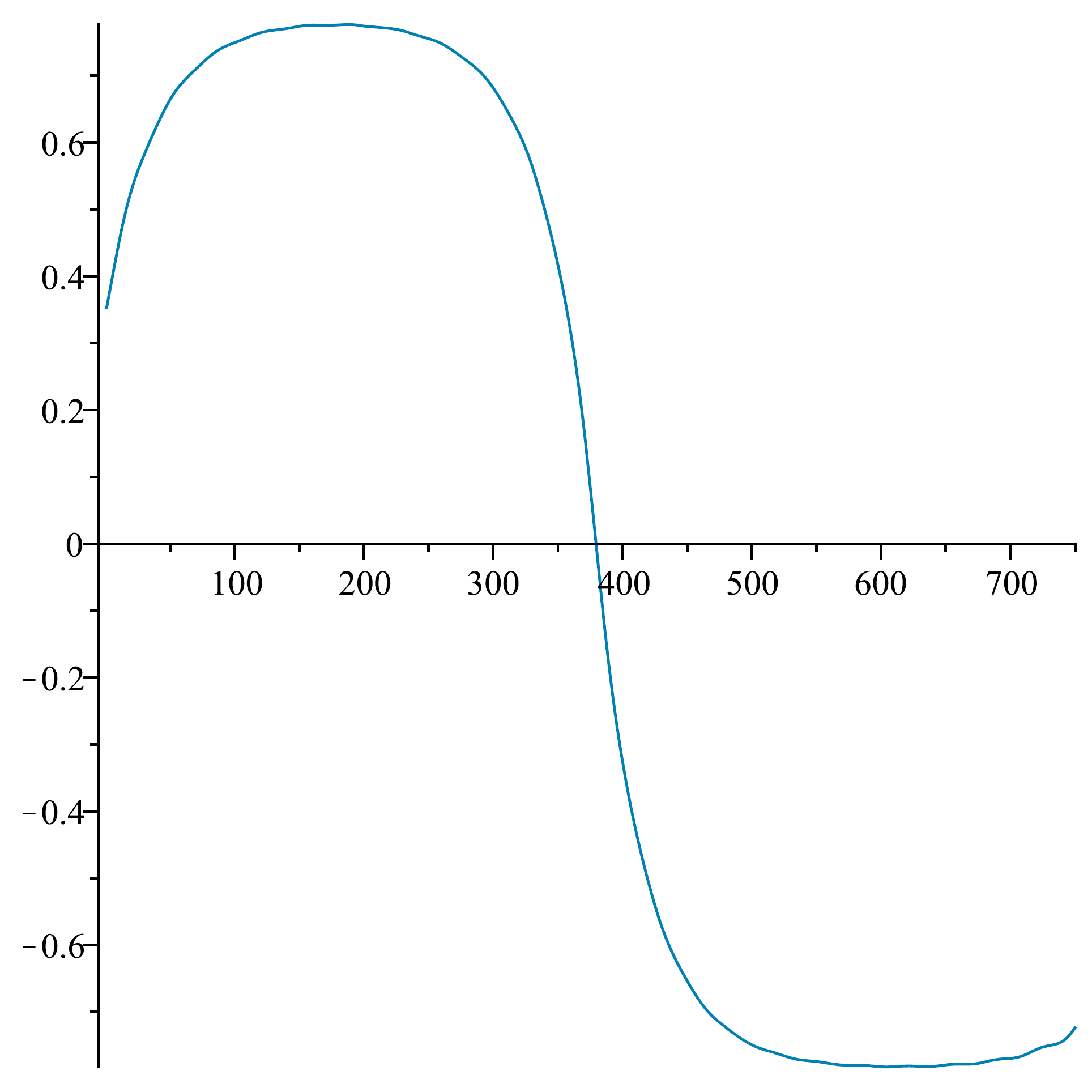}
\end{figure}

Die Gl\"attungsverfahren MA, EA und rEA liefern auch hier kein ausreichendes Ergebnis. Das wird auch mit steigender Gleitl\"ange nicht besser. Zwar werden die Amplituden der Schwingungen auf den Plateaus des Rechtecksignals verringert, verschwinden jedoch auch bei Erh\"ohung der Gleitl\"ange nicht in einem befriedigendem Ma\ss e. Auff\"allig ist jedoch insbesondere, dass die Gibbsschen \"Uberschwinger von keinem der drei Verfahren aufgefangen werden. An diesem Beispiel sieht man auch sehr gut, dass der EA die Symmetrie  der Ausgangsdaten nicht wiederspiegelt.

Ein optimales Ergebnis liefert auch in diesem Beispiel der SEA sowohl hinsichtlich der Gl\"attung schon bei mittlerer Gleitl\"ange (hier 40) als auch hinsichtlich der Symmetrieerhaltung.

\section{Die Behandlung von Ausrei\ss ern}

\subsection{Der Median und der gleitende Median}\label{ausrei}

Gegeben sei die Menge  von $2\ell+1$ Werten, $\{z_1,\ldots,z_{2\ell+1}\}$. Zur Bildung des {\em Medians} dieser Menge geht man wie folgt vor. Man sortiere die Werte der Gr\"o\ss e nach und erh\"alt 
$\{z_{i_1},z_{i_2},\ldots,z_{i_\ell},\ldots,z_{i_{2\ell+1}}\}$ mit
$z_{i_1}\leq z_{i_2}\leq\ldots\leq z_{i_{2\ell+1}}$. Der Median ist dann der Wert in der Mitte: $m:=z_{i_{\ell+1}}$.

Es sei nun wie in Abschnitt \ref{glaettung} der Datensatz $(y_1,y_2,\ldots, y_N,\ldots)$ als eine Reihe von Messpunkten gegeben. Der {\em gleitende Median mit Gleitl\"ange $2\ell+1$} dieser Messreihe ist eine Reihe von Daten, $m_i$, die durch die folgende Vorschrift definiert sind:
\begin{equation}
m_i := \text{\ Median der Werte }\{y_{i-\ell},\ldots,y_i,\ldots y_{i+\ell}\}
\end{equation}
f\"ur $i=\ell+1,\ldots,N-\ell$. Auch hier entf\"allt nat\"urlich die obere Grenze, wenn der Ausgangsdatensatz endlich ist. Der Algorithmus zur Berechnung des Medians ben\"otigt ein Sortierverfahren, f\"ur welches man etwa in \cite{inf} Beispiele mit zugeh\"origen Algorithmen finden kann.

\subsection{Die Gl\"attung von Messreihen mit Ausrei\ss ern}\label{aus1}

\begin{figure}[H]\caption{Originaldaten einer Messung mit Ausrei\ss ern}\label{test-ausgang}\centering
\includegraphics[width=5cm]{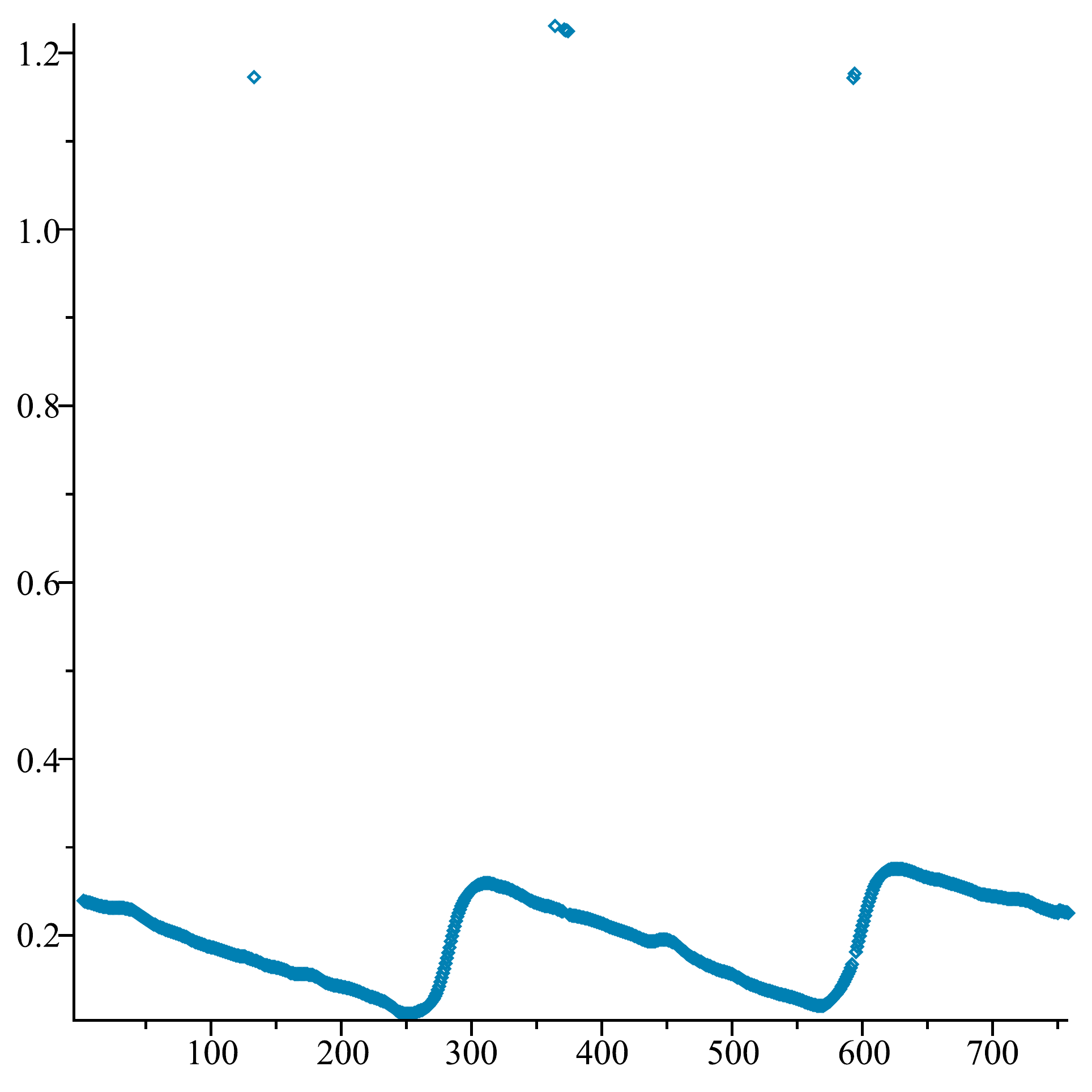}
\end{figure}

\begin{figure}[H]\caption{Gl\"attung mittels EA und rEA ($\alpha^{-1}=40$)}\label{test-EA-40}\centering
\includegraphics[width=0.48\textwidth]{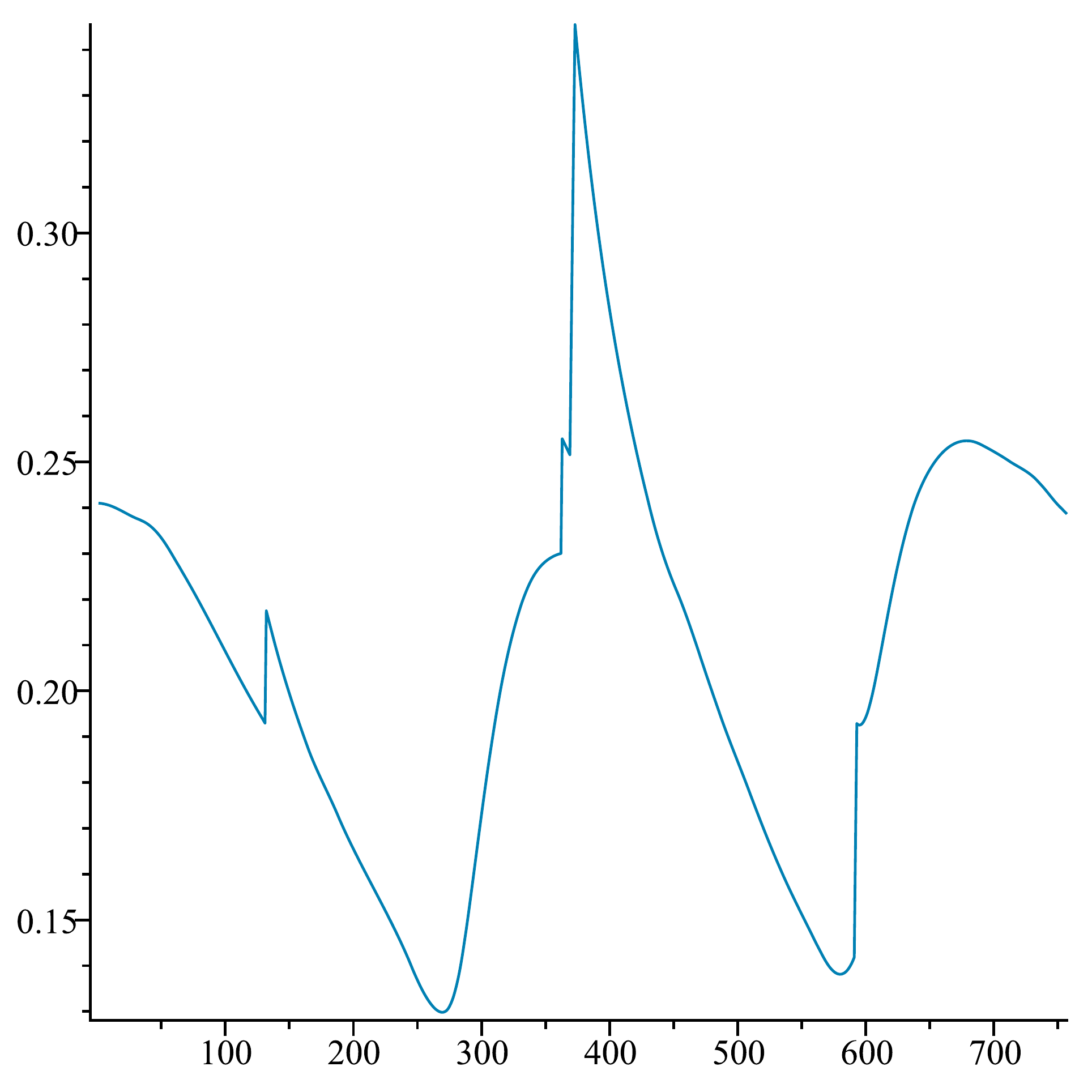}
\includegraphics[width=0.48\textwidth]{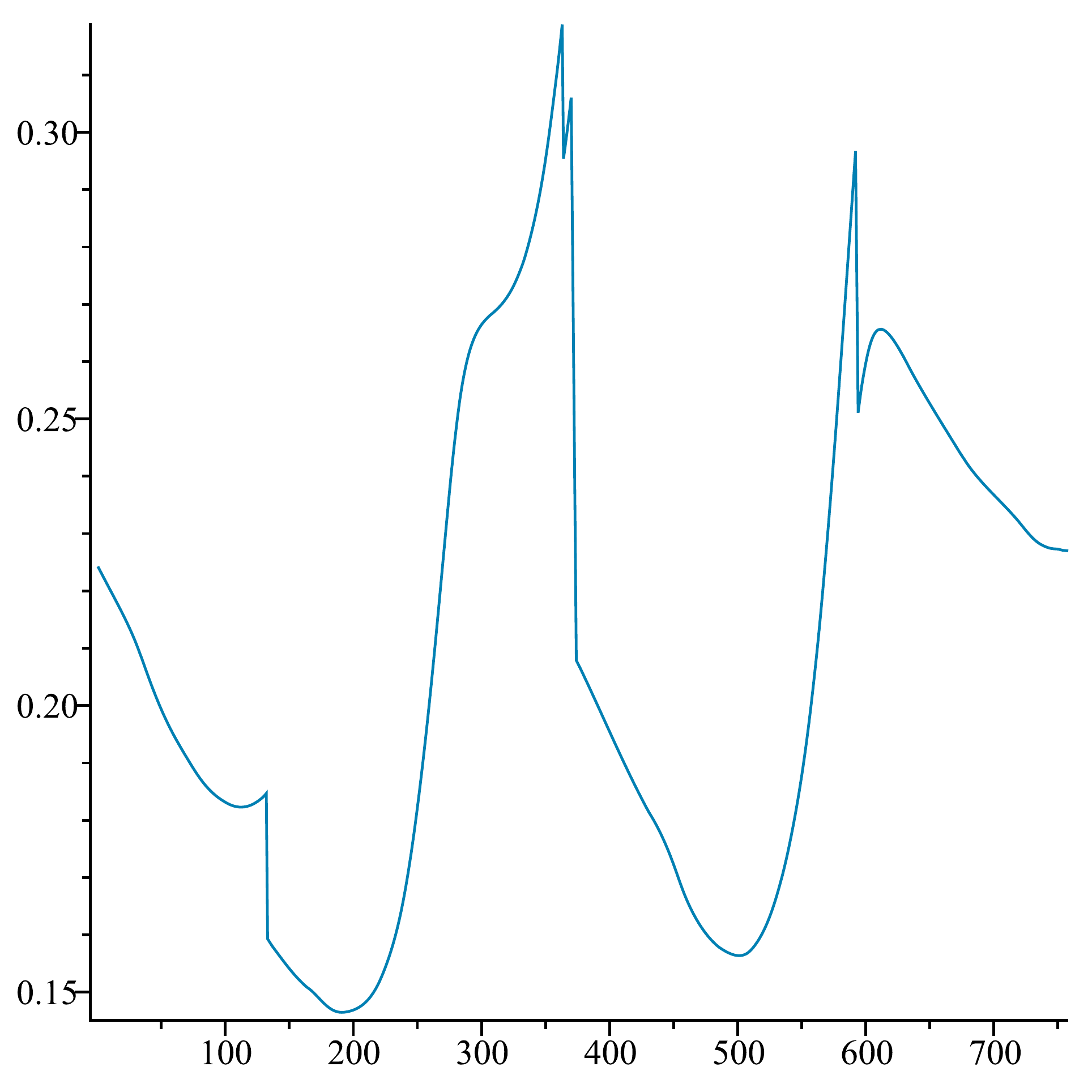}
\end{figure}

\begin{figure}[H]\caption{Gl\"attung mittels MA ($n+1=40$) und  mittels SEA ($\alpha^{-1}=40$)} \label{test-SEA-40}\centering
\includegraphics[width=0.48\textwidth]{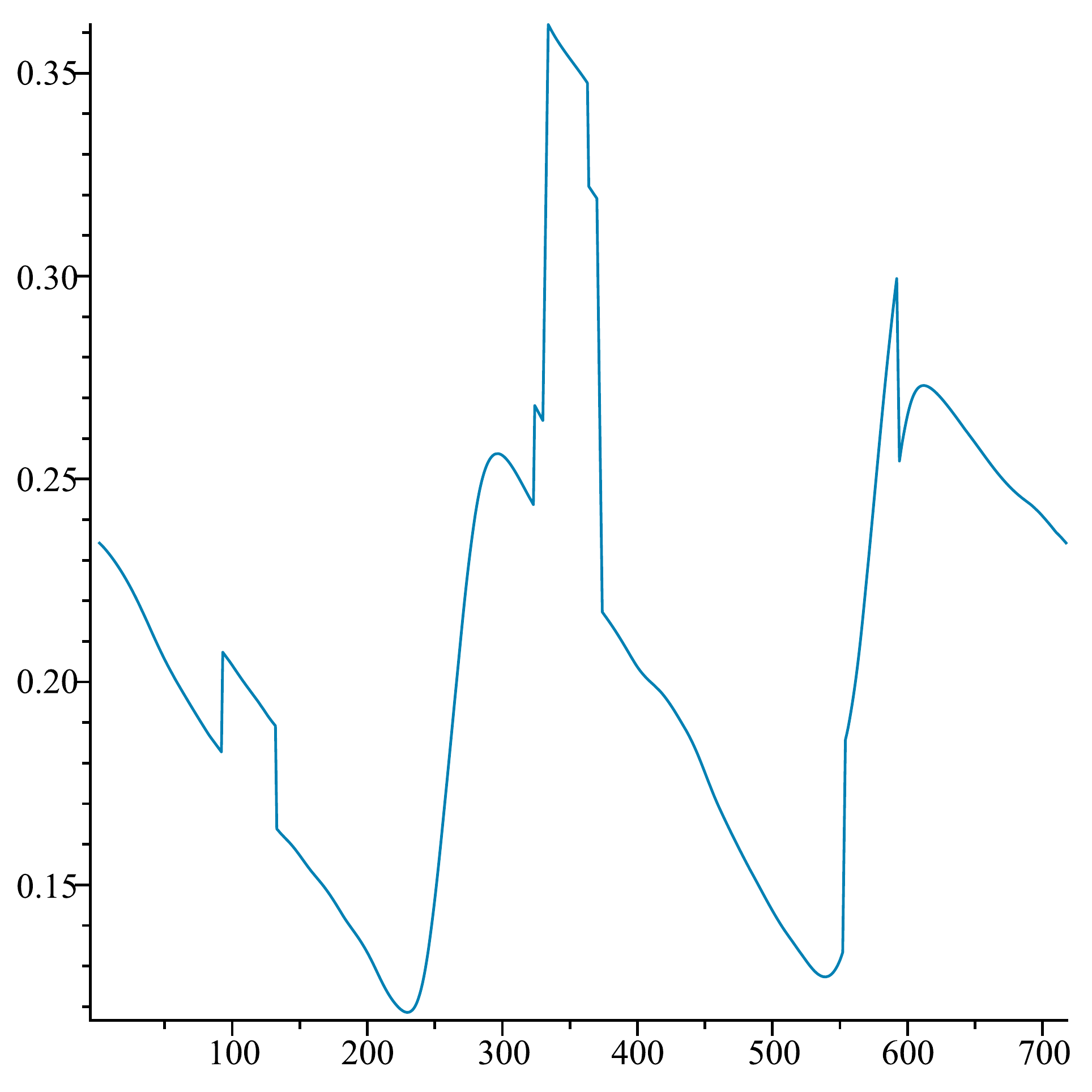}
\includegraphics[width=0.48\textwidth]{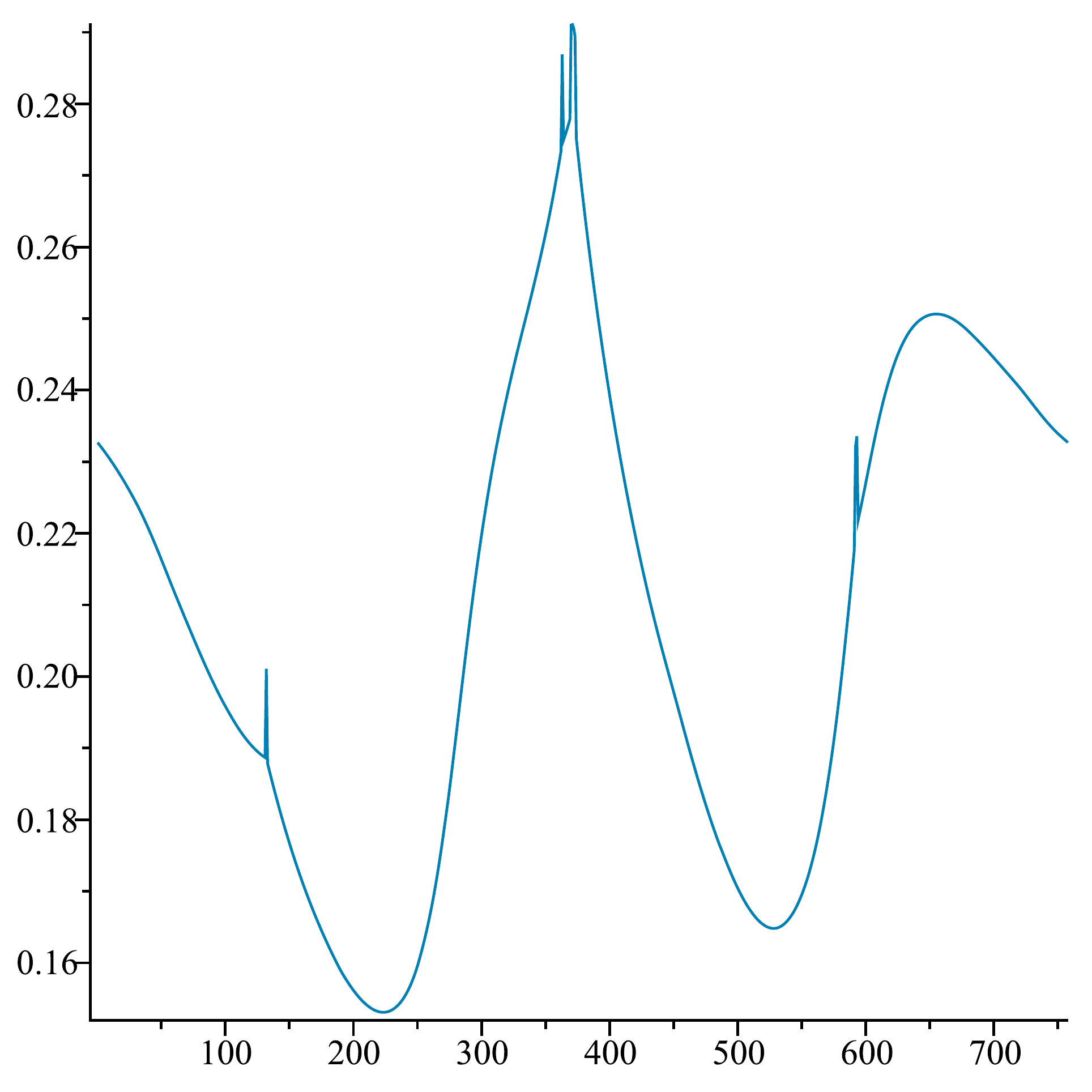}
\end{figure}

\begin{figure}[H]\caption{Medianbildung mit Gleitl\"ange $9$ und $40$} \label{test-Median}\centering
\includegraphics[width=0.48\textwidth]{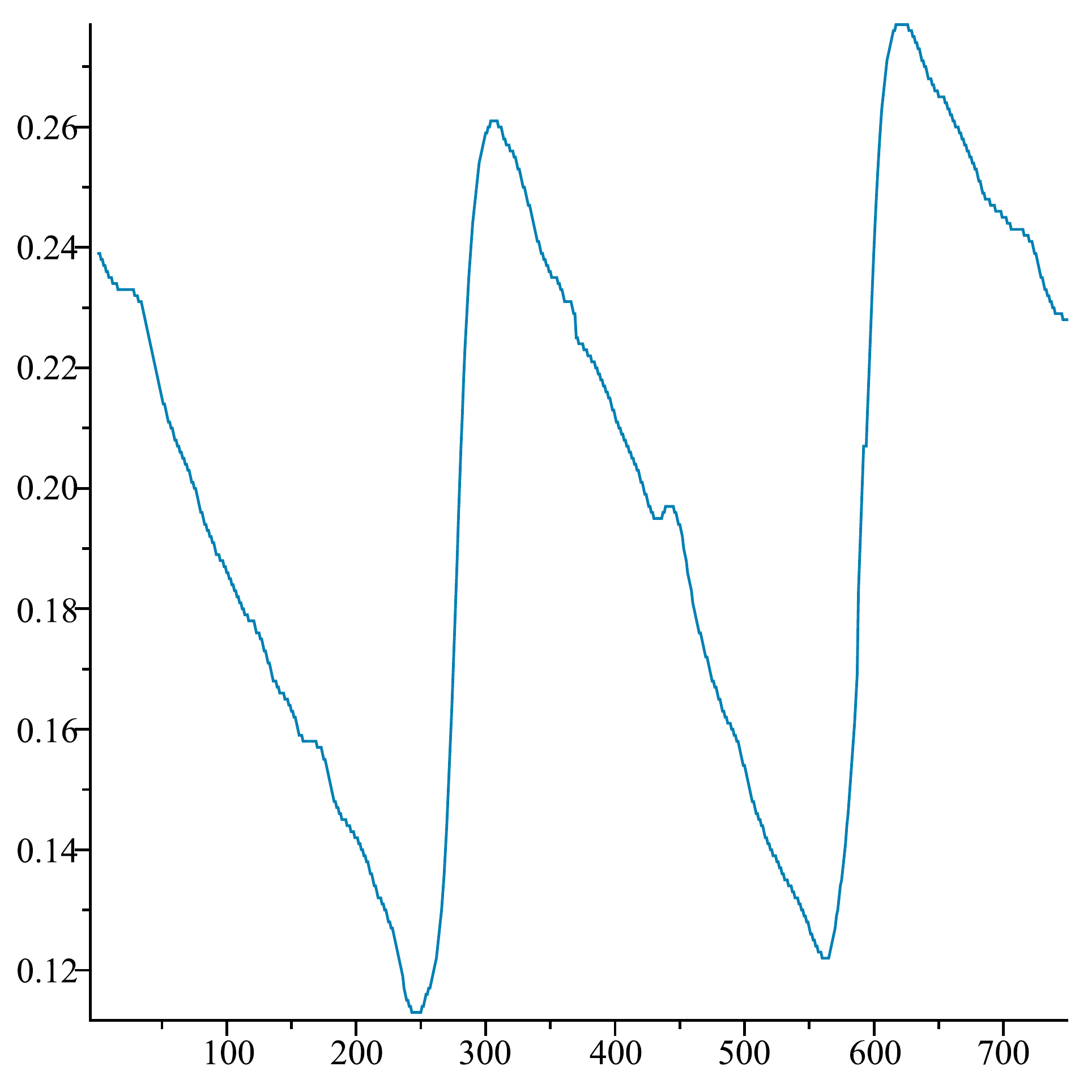}
\includegraphics[width=0.48\textwidth]{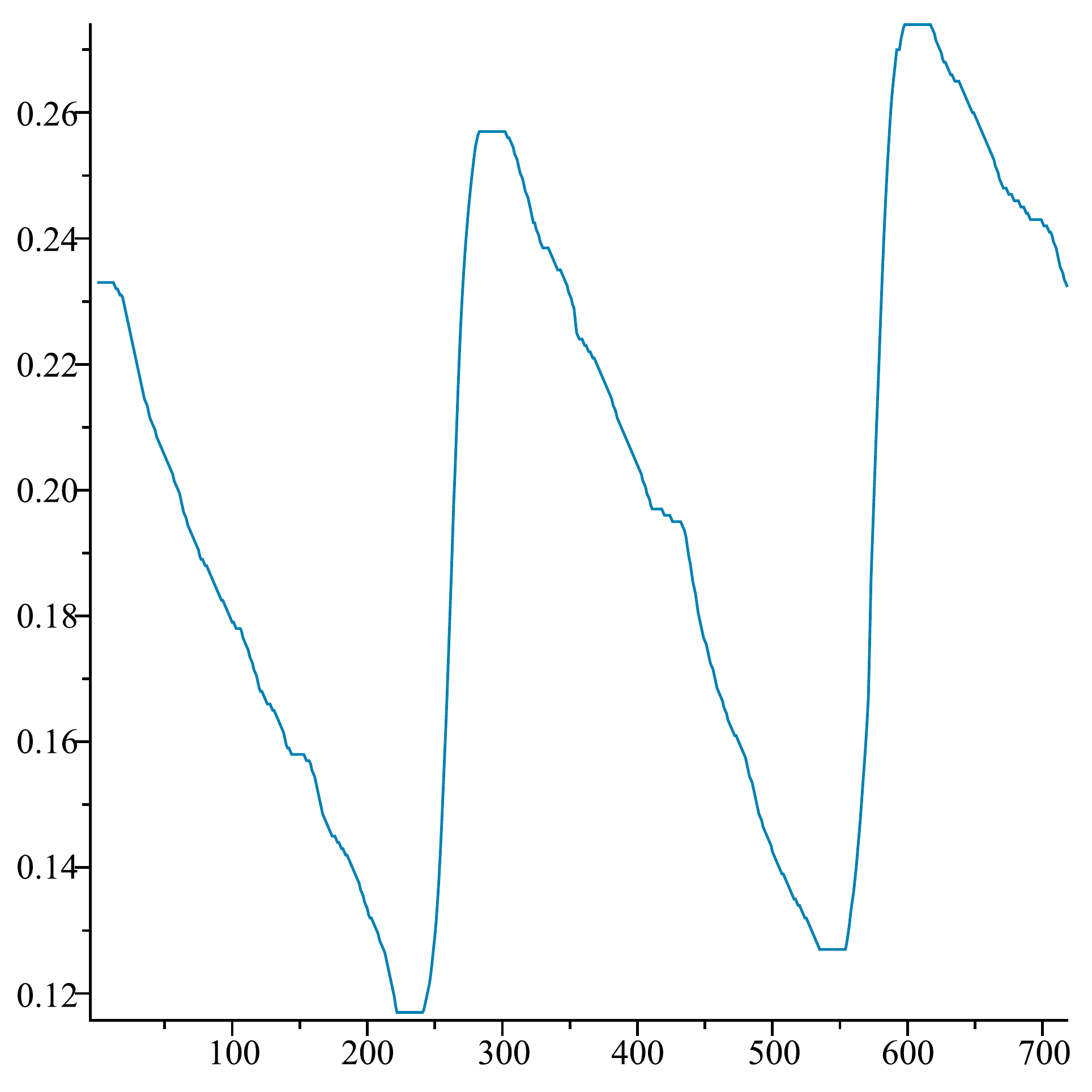}
\end{figure}

Ein Problem, das in der Praxis auftreten kann, sind Ausrei\ss er im Datensatz. Im Gegensatz zu den Artefakten aus Abschnitt \ref{sub2} sind Ausrei\ss er dadurch charakterisiert, dass sie eine wesentliche gr\"o\ss ere H\"ohe relativ zu den Extrema des Datensatzes haben, siehe\footnote{Hier haben wir in die Originalmessreihe aus Abbildung \ref{fig1} und Abschnitt \ref{sub1} k\"unstlich Ausrei\ss er eingef\"ugt.} Abbildung \ref{test-ausgang}.
Solche extremen Werte w\"urden bei den hier vorgestellten Verfahren immer Einfluss auf die Gl\"attung haben und weiterhin als zus\"atzliche Extrema zu sehen sein, siehe Abbildungen \ref{test-EA-40} und \ref{test-SEA-40}. Um diese Ausrei\ss er auszusortieren kann man zun\"achst -- bevor der eigentliche Gl\"attungsprozess mit den hier vorgestellten Methoden erfolgt -- den Datensatz mit Hilfe eines gleitenden Medians bearbeiten. Das Ergebnis dieser Medianbildung findet man in Abbildung \ref{test-Median} mit minimaler Gleitl\"ange (hier $9$) sowie mit Gleitl\"ange vergleichbar mit denen der verwendeten Gl\"attungsverfahren. Man sieht neben der Eigenschaft die Ausrei\ss er abzufangen, dass der Median gl\"attende Eigenschaften hat. Hierzu verweisen wir auf den n\"achsten Abschnitt.


Damit die Ausrei\ss er aus dem Datensatz wirklich eliminiert werden, muss die Gleit\-l\"ange des gleitenden Medians geeignet gew\"ahlt werden. Die Wahl der {\em minimalen Gleitl\"ange} zur Eliminierung aller Ausrei\ss er geschieht wie folgt:

W\"ahle die Gleitl\"ange $2\ell+1\leq N$ des Medians derart, dass sich f\"ur jede $(2\ell+1)$-elementige Teilmenge $\{y_k,\ldots,y_{k+2\ell}\}$  des Ausgangsdatensatzes h\"ochstens $\ell$ Ausrei\ss er in dieser Teilmenge befinden. Insbesondere reicht bei einzelnen Ausrei\ss ern, die nicht geh\"auft auftreten -- das hei\ss t, der Abstand zwischen einzelnen Ausrei\ss er ist mindestens zwei -- eine Gleitl\"ange von $3$. In unserem Beispiel gem\"a\ss~Abbildung \ref{test-ausgang} gibt es vier aufeinanderfolgende Ausrei\ss er, so dass eine Gleitl\"ange von $9$ notwendig ist.

\subsection{Der gleitende Median als Gl\"attungsverfahren}\label{aus2}

\begin{figure}[H]\caption{MA und Medianbildung mit Gleitl\"ange $20$} \label{11neu}\centering
\includegraphics[width=0.48\textwidth]{ydata1-MA-20.pdf}
\includegraphics[width=0.48\textwidth]{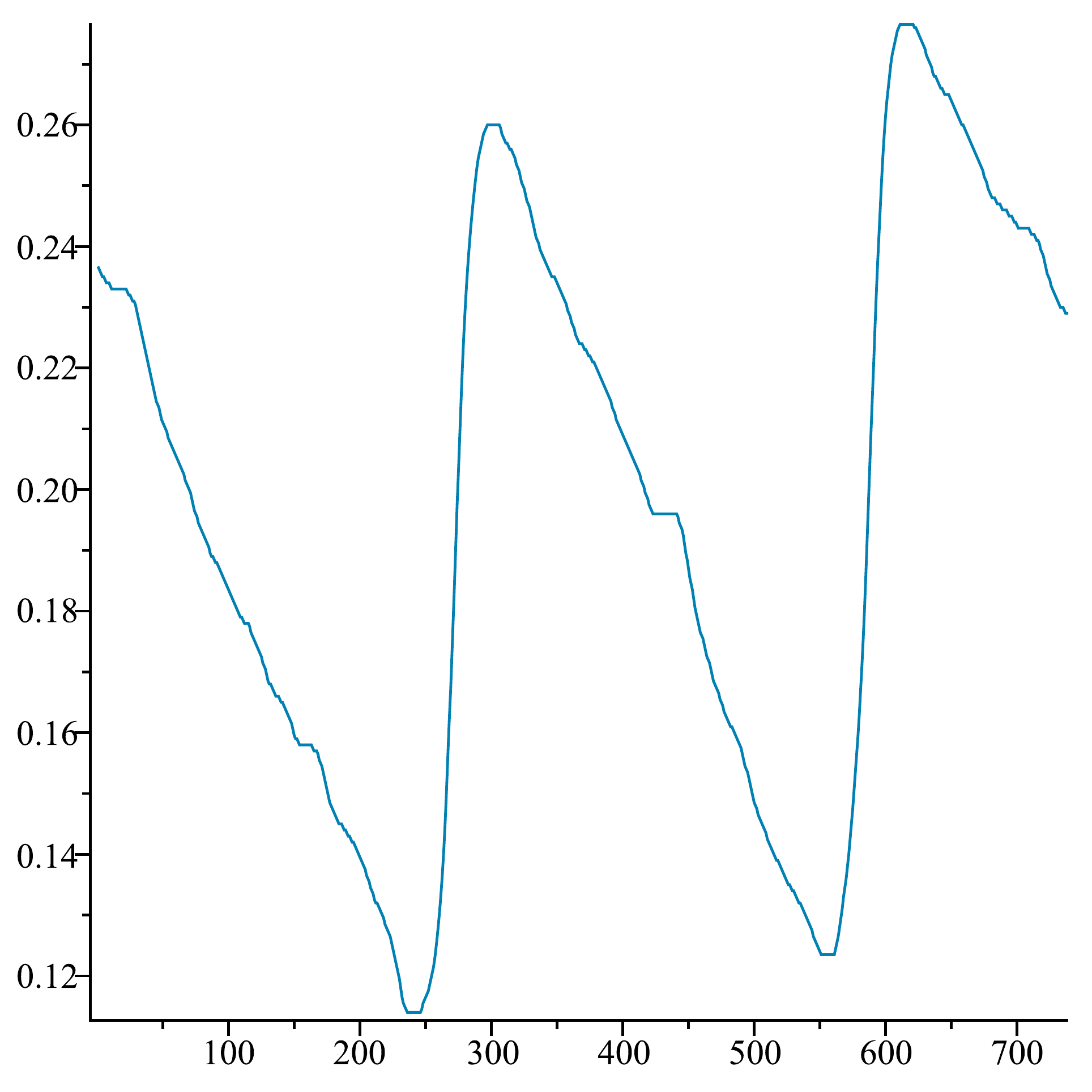}
\end{figure}

\begin{figure}[H]\caption{MA und Medianbildung mit Gleitl\"ange $40$} \label{median}\centering 
\includegraphics[width=0.48\textwidth]{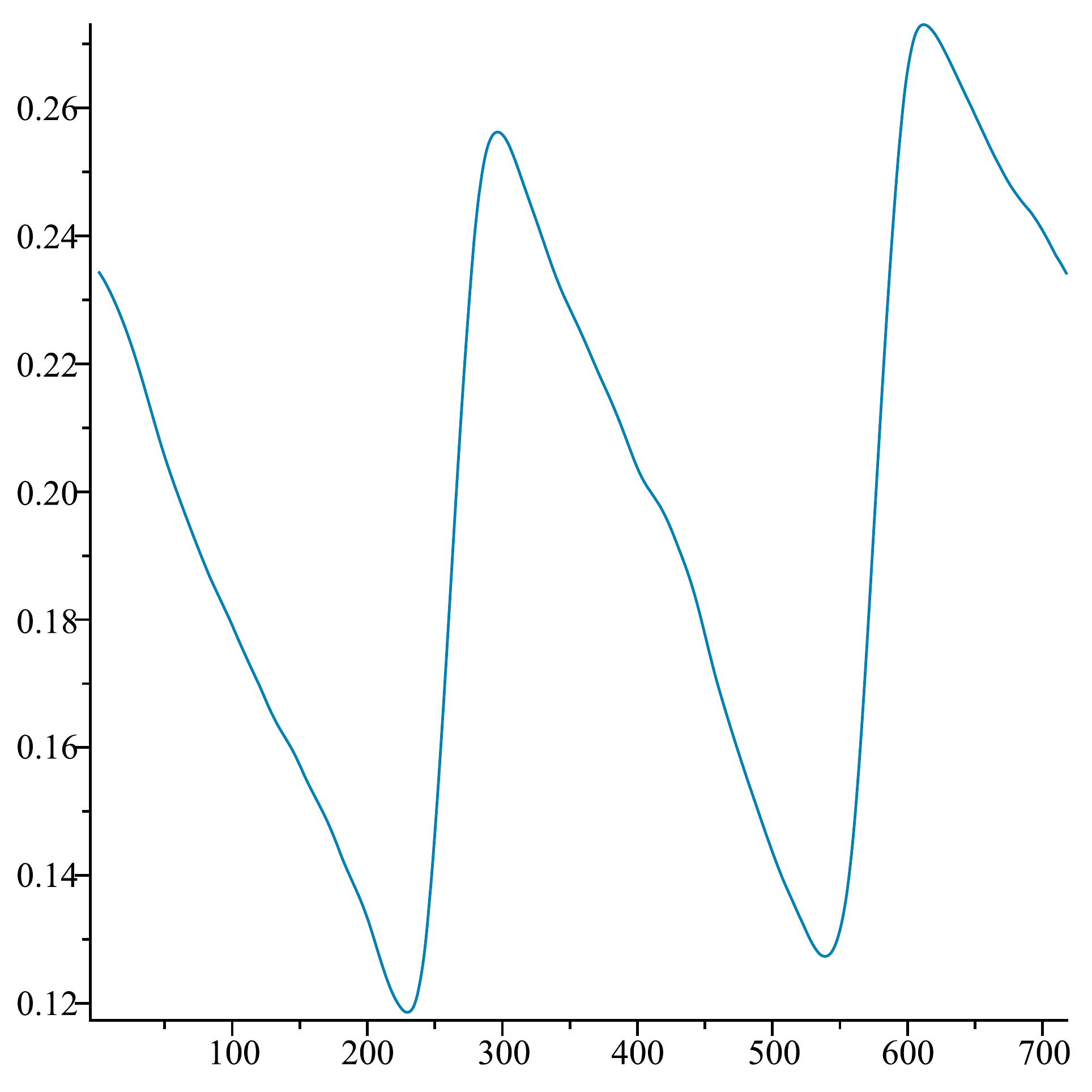}
\includegraphics[width=0.48\textwidth]{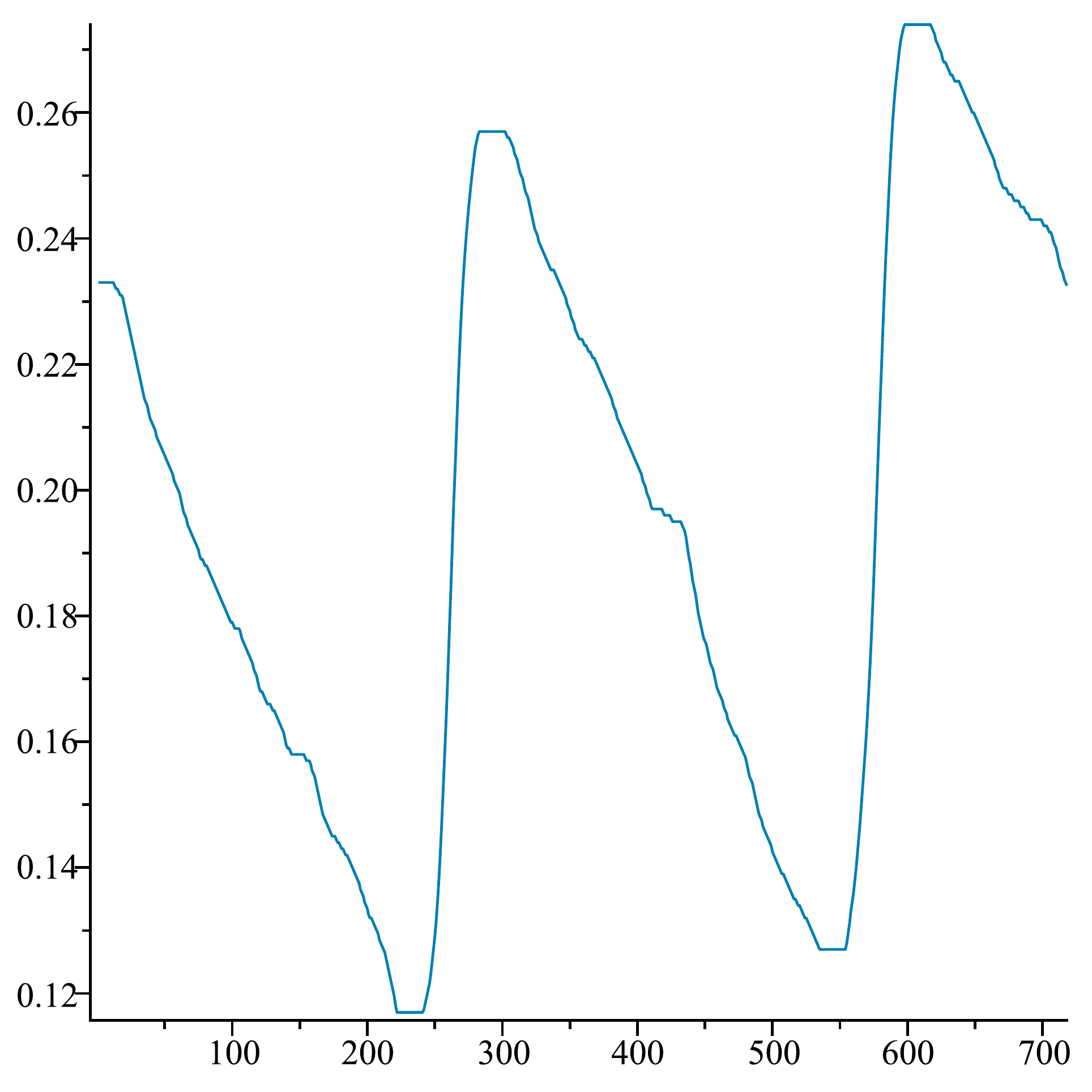}
\end{figure}

Mit Blick auf Abbildung \ref{test-Median} kann man geneigt sein, den gleitenden Median selbst als Gl\"attungsverfahren in Betracht zu ziehen. Aufgrund der bisherigen Diskussion und der folgenden Betrachtung ist davon allerdings abzuraten.  
\"Ahnlich wie beim MA muss man auch beim gleitenden Median eine Verk\"urzung des Datensatzes in Kauf nehmen, wenn man ihn zur Gl\"attung benutzt.
In den folgenden Abbildungen sind zum Vergleich der MA (links) und der gleitende Median (rechts) mit Gleitl\"ange $20$ bzw. $40$ zu sehen (Abbildungen \ref{11neu} und \ref{median}). Als Grundlage haben wir hier die Originaldaten der Laserausmessung aus Abschnitt \ref{sub1} gew\"ahlt.
Wie man schon an diesem Beispiel erkennt, besitzt der MA im Gegensatz zum gleitenden Median die besseren Gl\"attungseigenschaften. Dar\"uberhinaus neigt der Datensatz nach Anwendung des gleitenden Medians bei steigender Gleitl\"ange zur Plateaubildung in den Extrema. 


\section{Praxisdaten im Vergleich}\label{sub6}

\subsection{Vorbereitung der Praxisdaten}\label{sub7}

In diesem Abschnitt stellen wir anhand von vier Messreihen zu der in Abschnitt \ref{einl} beschriebenen technischen Fragestellung unsere diskutierten Verfahren vor und begr\"unden anhand der Ergebnisse unser Fazit. Desweiteren werden wir eine notwendige Vorbereitung der Datens\"atze aus der Praxis begr\"unden.
Wir betrachten hier jeweils zwei Messungen f\"ur die grobe Seite und f\"ur die feine Seite der Zahnscheibe.

Die Messreihen sind derart, dass vor und nach den eigentlichen Messwerten ein Nulllauf auftritt. Die L\"ange dieser Nulll\"aufe ist f\"ur die vorliegenden Messreihen wie folgt:
\[\text{\renewcommand{\arraystretch}{1.5}
\begin{tabular}{|l|c|c|
				c|c|
				}\hline
Messung 		& I 	& II 	
				& III 	& IV 	
				\\\hline
Nullvorlauf 	& 40 	& 150 	
				& 80 	& 30 	
				\\\hline
Nullnachlauf 	& 110 	& 110 	
				& 40 	& 80 	
				\\\hline
\end{tabular}
}\]
Da der Nulllauf verh\"altnism\"a\ss ig lang ist, kann man ihn nicht mehr als Menge von Ausrei\ss ern behandeln, vergleiche Abschnitt \ref{ausrei}. Eine Medianbildung ist also zur Beseitigung nicht angezeigt. Wie unsere Diskussion zeigen wird, ist eine Abtrennung des Nulllaufs als Vorbereitung des Datensatzes sinnvoll und sollte der Behandlung des Nulllaufs als Teil des Datensatzes vorgezogen werden.

In den Abbildungen \ref{101}-\ref{105} und \ref{201}-\ref{505} ist f\"ur die Messreihe  jeweils links die Messreihe inklusive Nulllauf sowie ihre Gl\"attungen zu sehen. Dabei haben wir bei den gegl\"atteten Reihen zur besseren Vergleichbarkeit auf den relevanten Teil der Grafik gezoomt. Rechts findet man dann jeweils die um den Nulllauf bereinigte Messreihe sowie ihre Gl\"attungen. 

Die erg\"anzende Abbildung \ref{106} zur Messreihe {\bf I} wird in der folgenden Bemerkung n\"aher erl\"autert.

Wir w\"ahlen eine kleine Gleit\"ange von 10 zum Vergleich der Gl\"attungsverfahren, und sehen, dass das in den meisten F\"allen gen\"ugt.

\subsection*{Bemerkungen zu den Abbildungen in Abschnitt \ref{Abbsec}}
\begin{figure}[H]\caption{{\bf Messreihe I}, Originaldaten und bereinigte Originaldaten}\label{101}
 \includegraphics[width=0.45\textwidth]{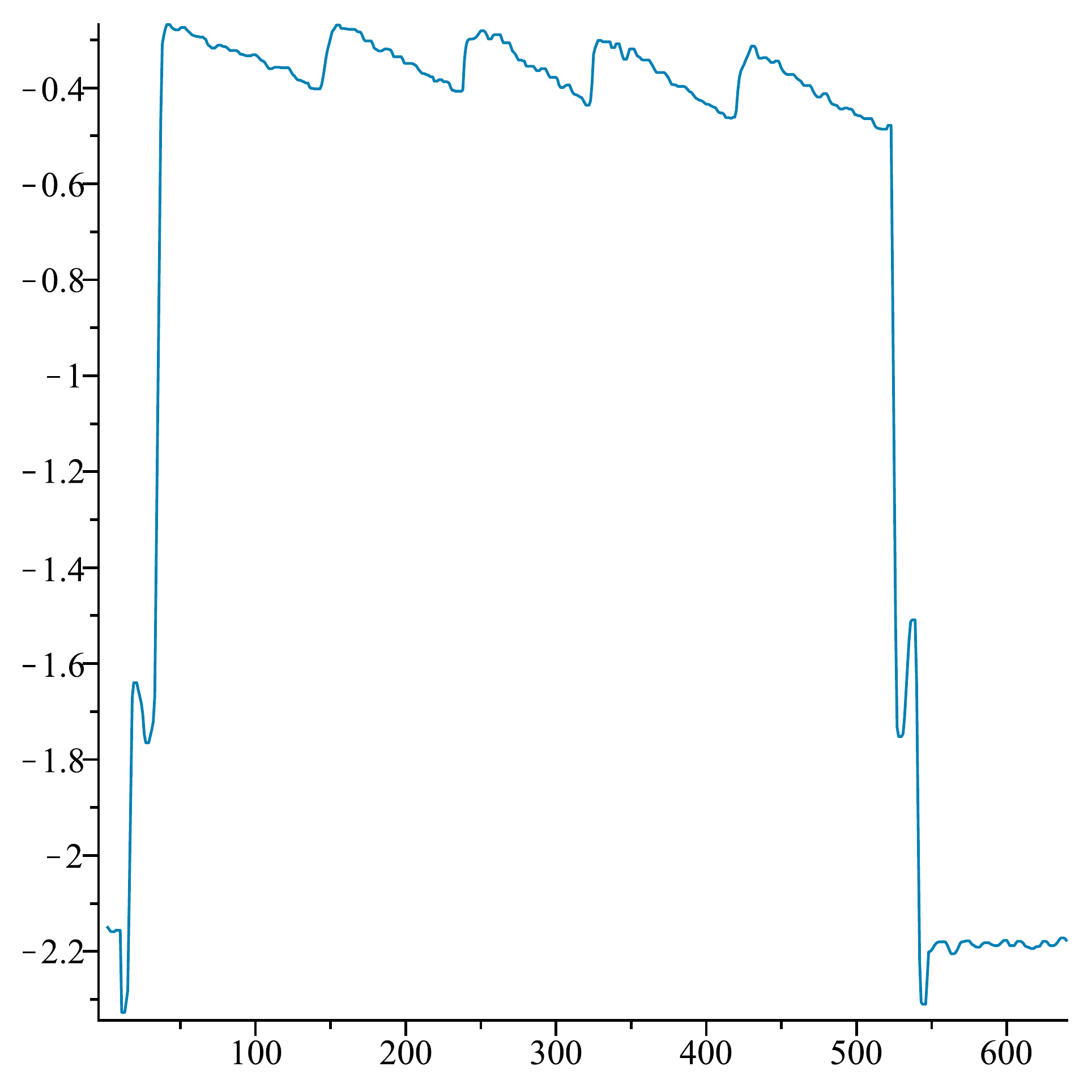}
 \includegraphics[width=0.45\textwidth]{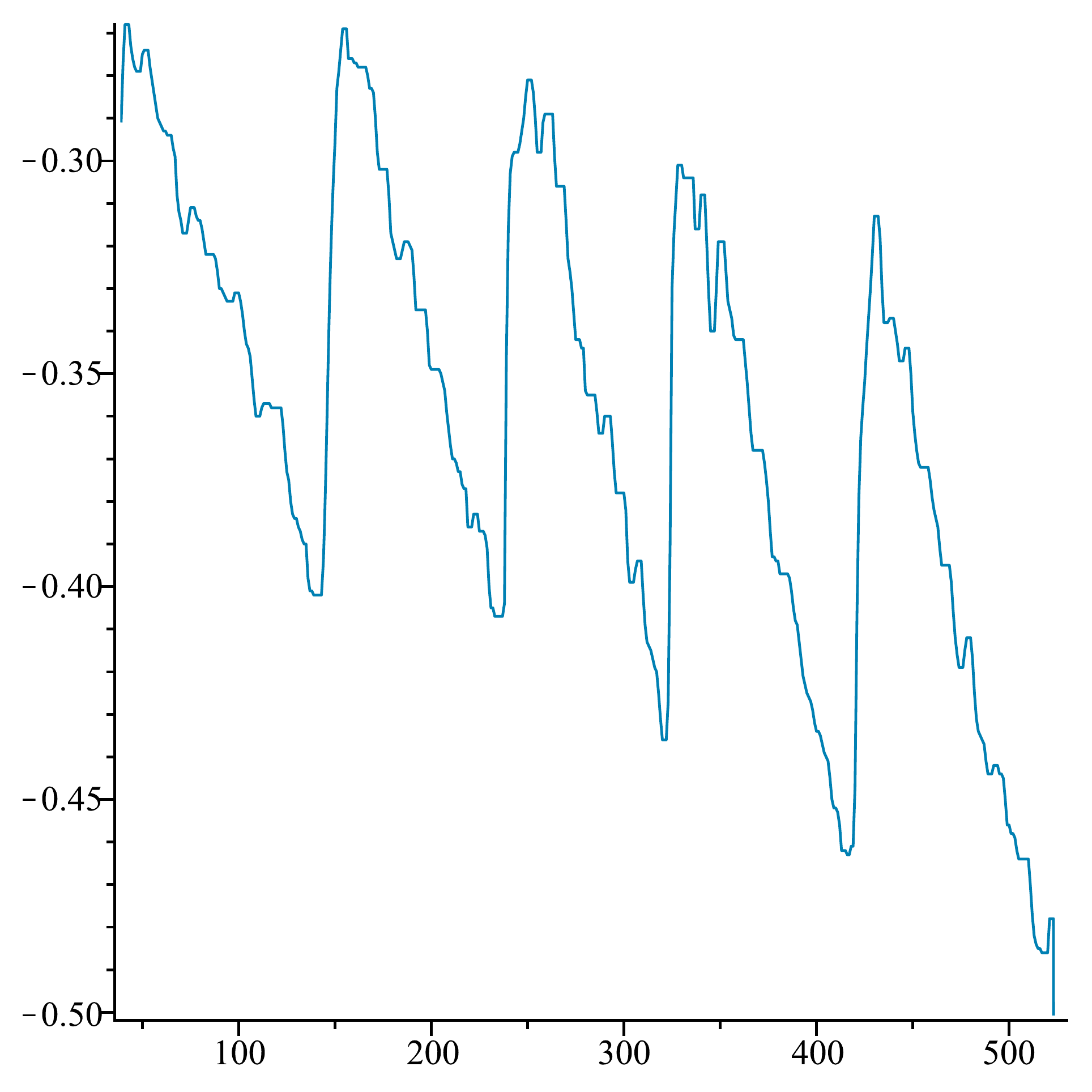}
\end{figure}

\begin{figure}[H]\caption{{\bf Messreihe I}, Gl\"attung mit MA, Gleitl\"ange $10$.}\label{102}
 \includegraphics[width=0.45\textwidth]{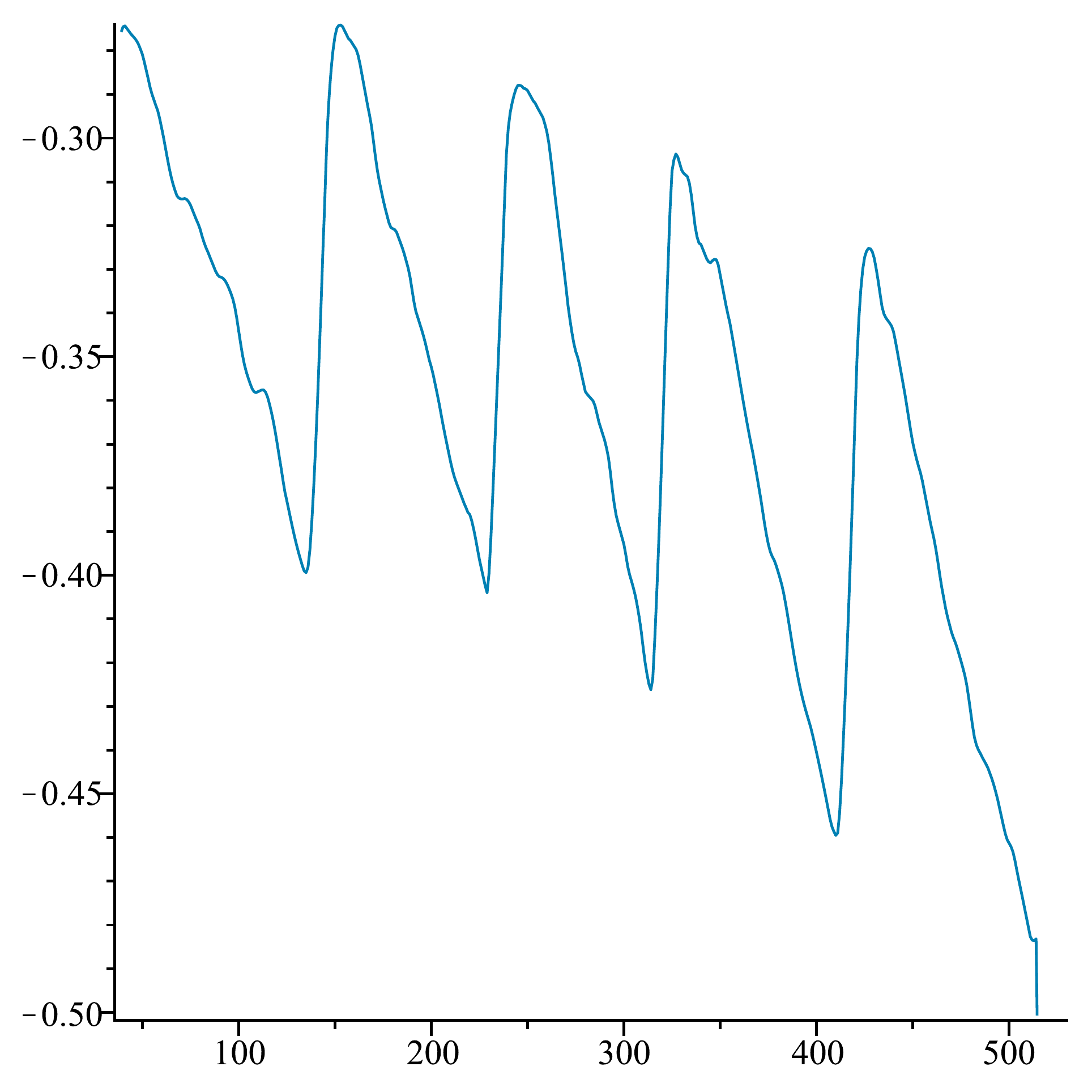}
 \includegraphics[width=0.45\textwidth]{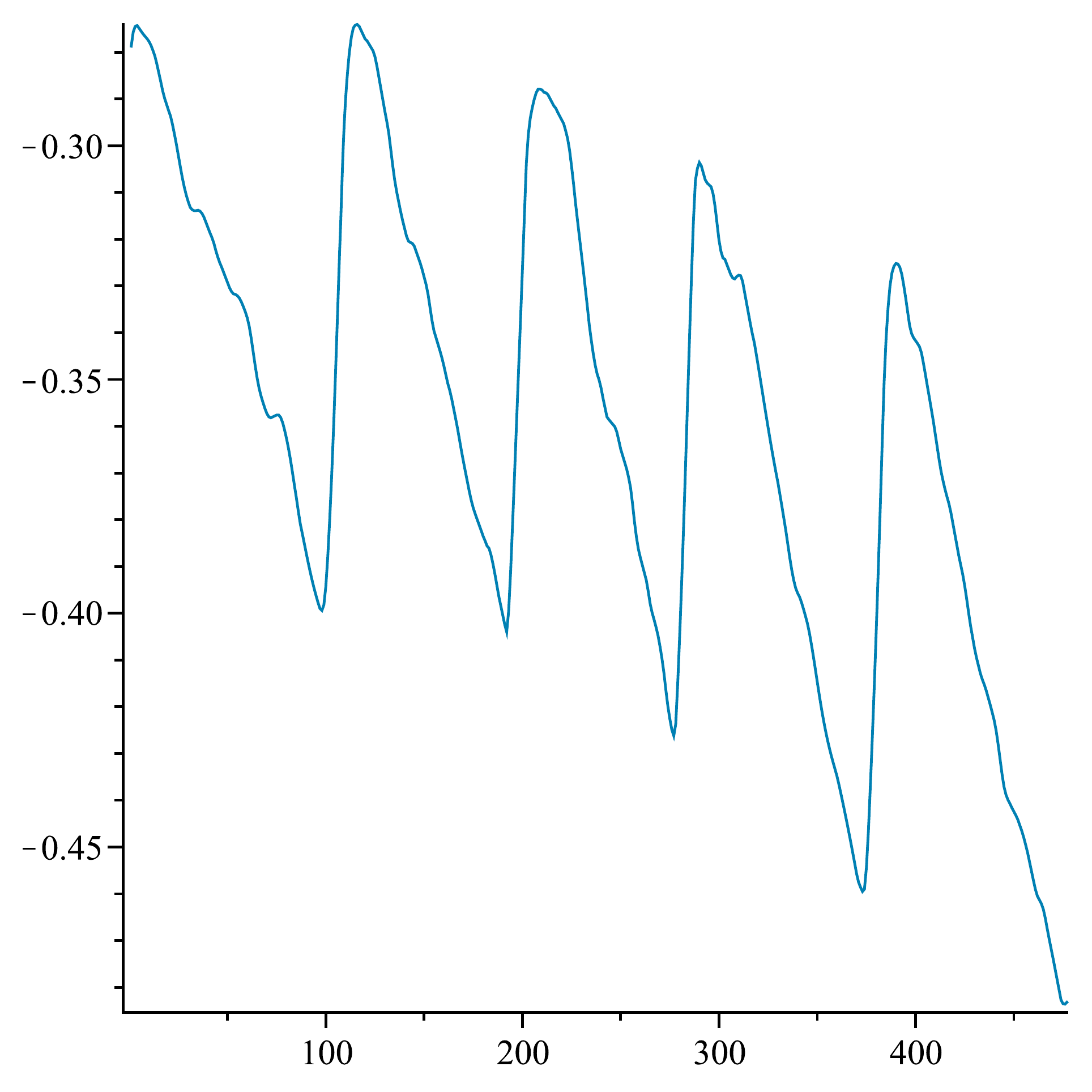}
\end{figure}

\begin{figure}[H]\caption{{\bf Messreihe I}, Gl\"attung mit EA, Gleitl\"ange $10$.}\label{103}
 \includegraphics[width=0.45\textwidth]{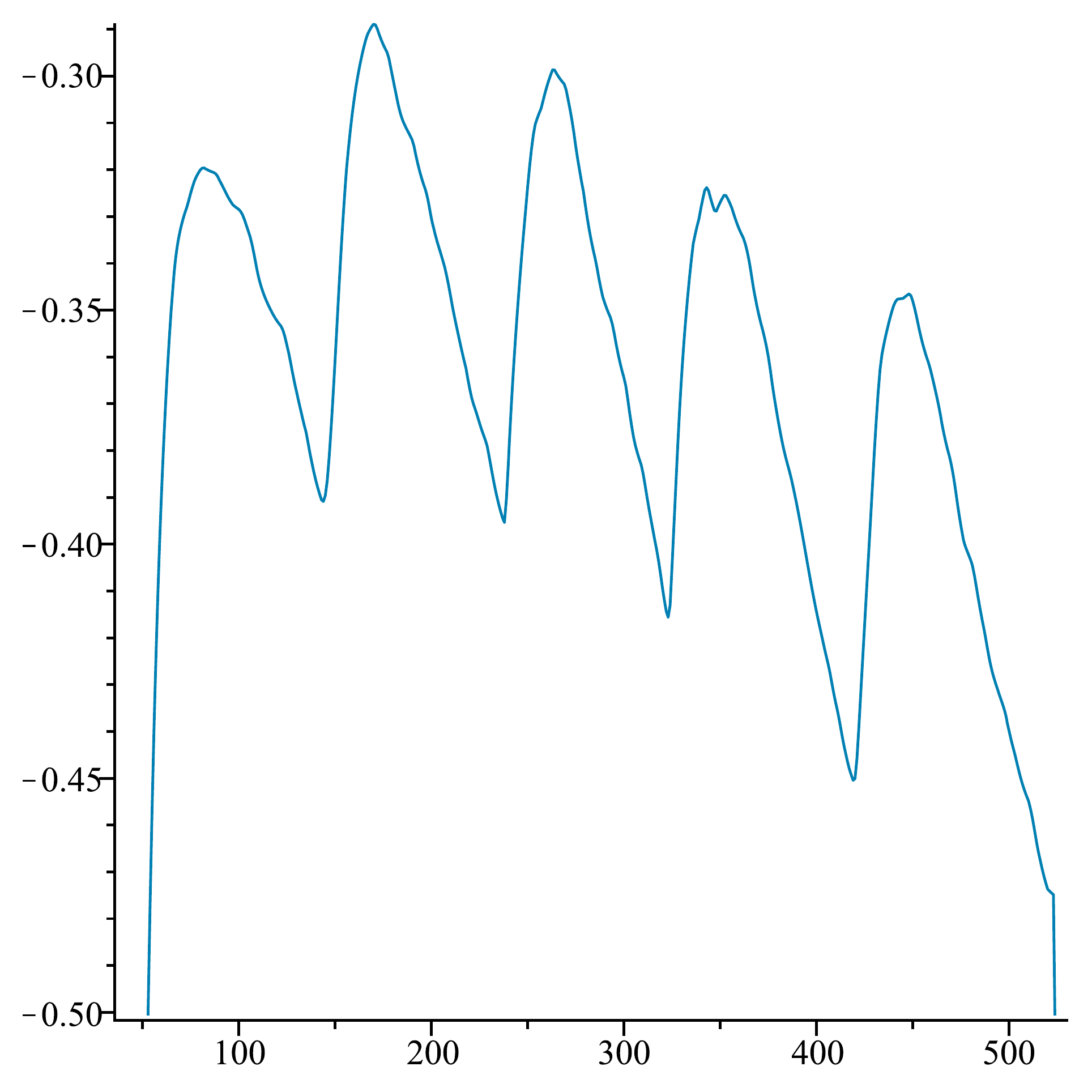}
 \includegraphics[width=0.45\textwidth]{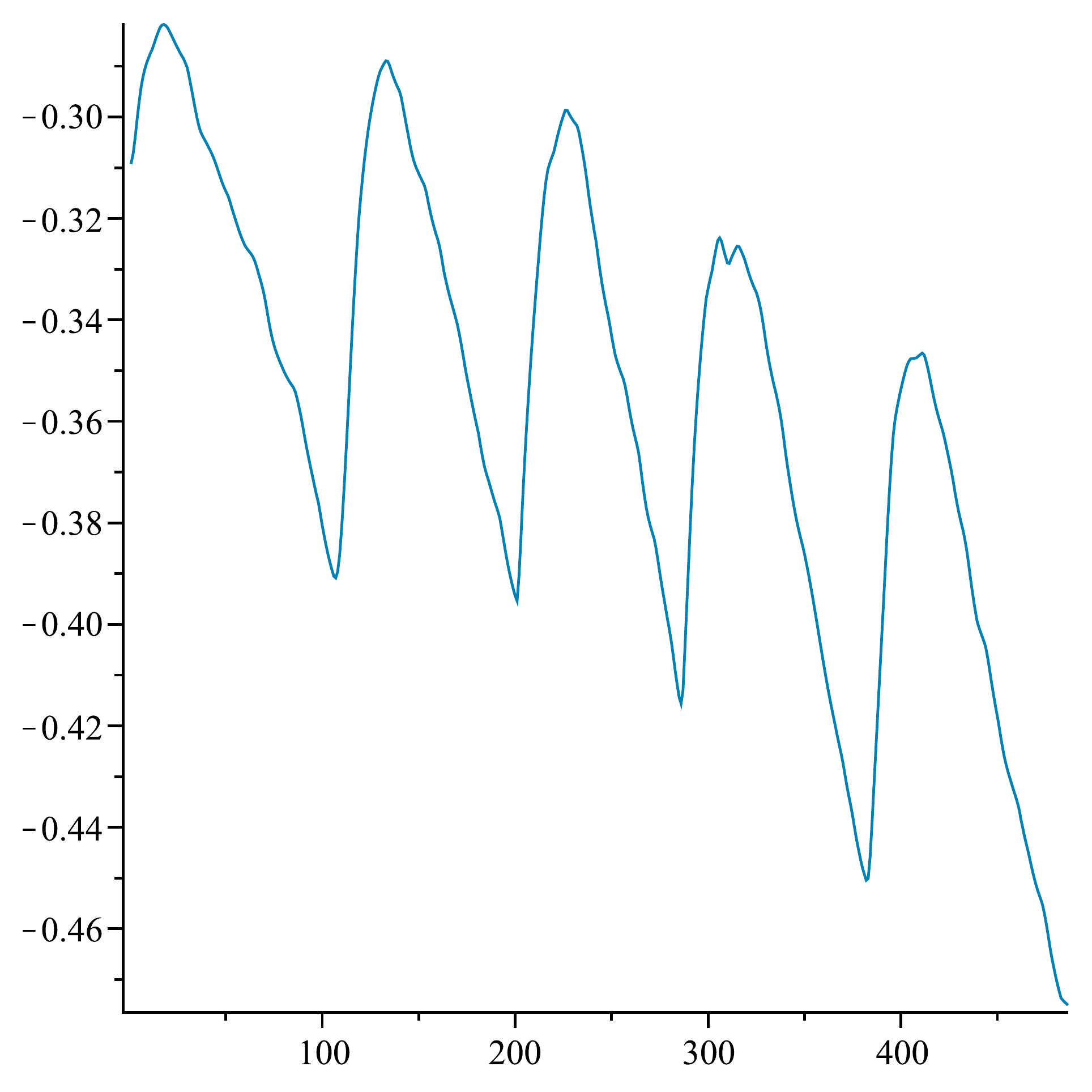}
\end{figure}

\begin{figure}[H]\caption{{\bf Messreihe I}, Gl\"attung mit rEA, Gleitl\"ange $10$.}\label{104}  
 \includegraphics[width=0.45\textwidth]{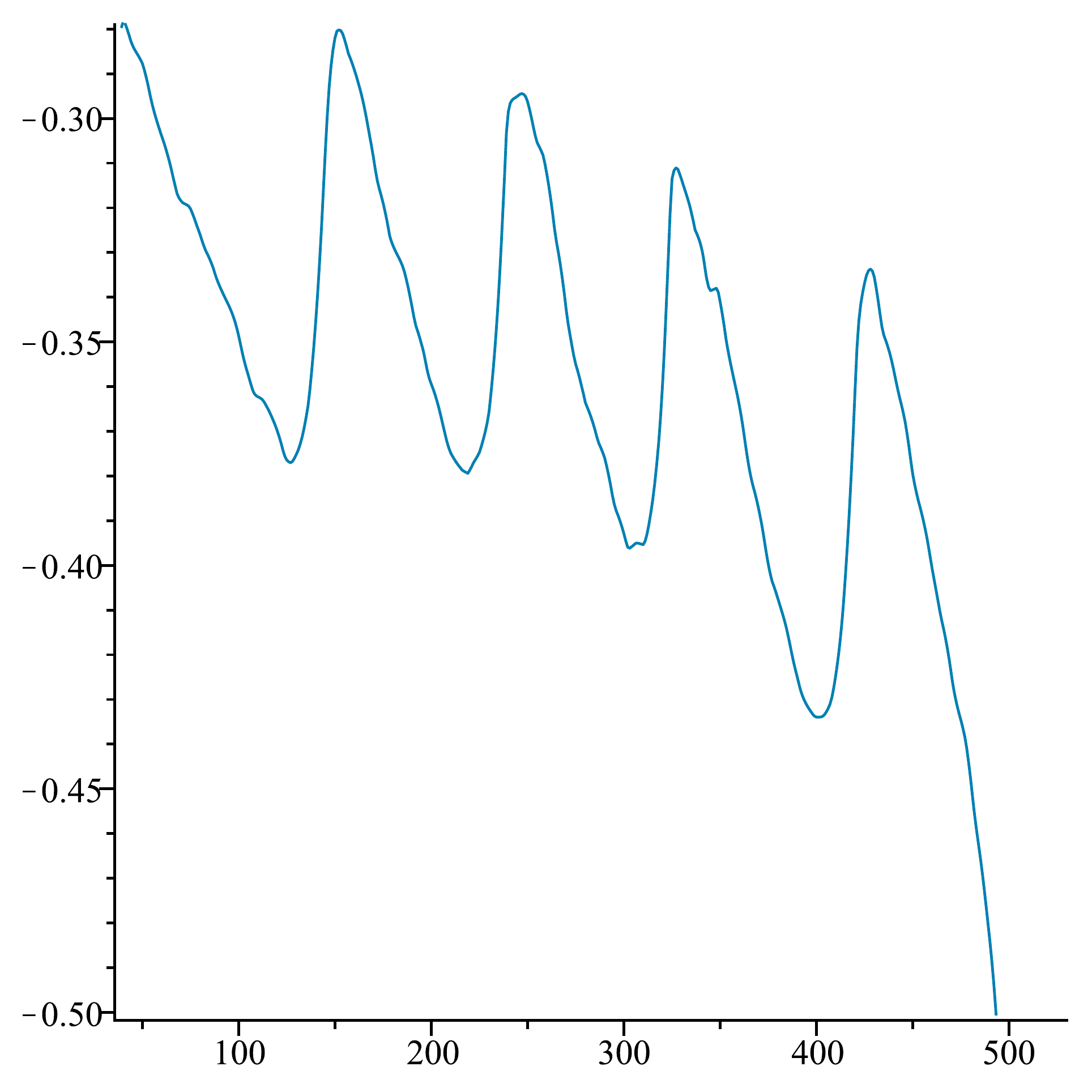}
 \includegraphics[width=0.45\textwidth]{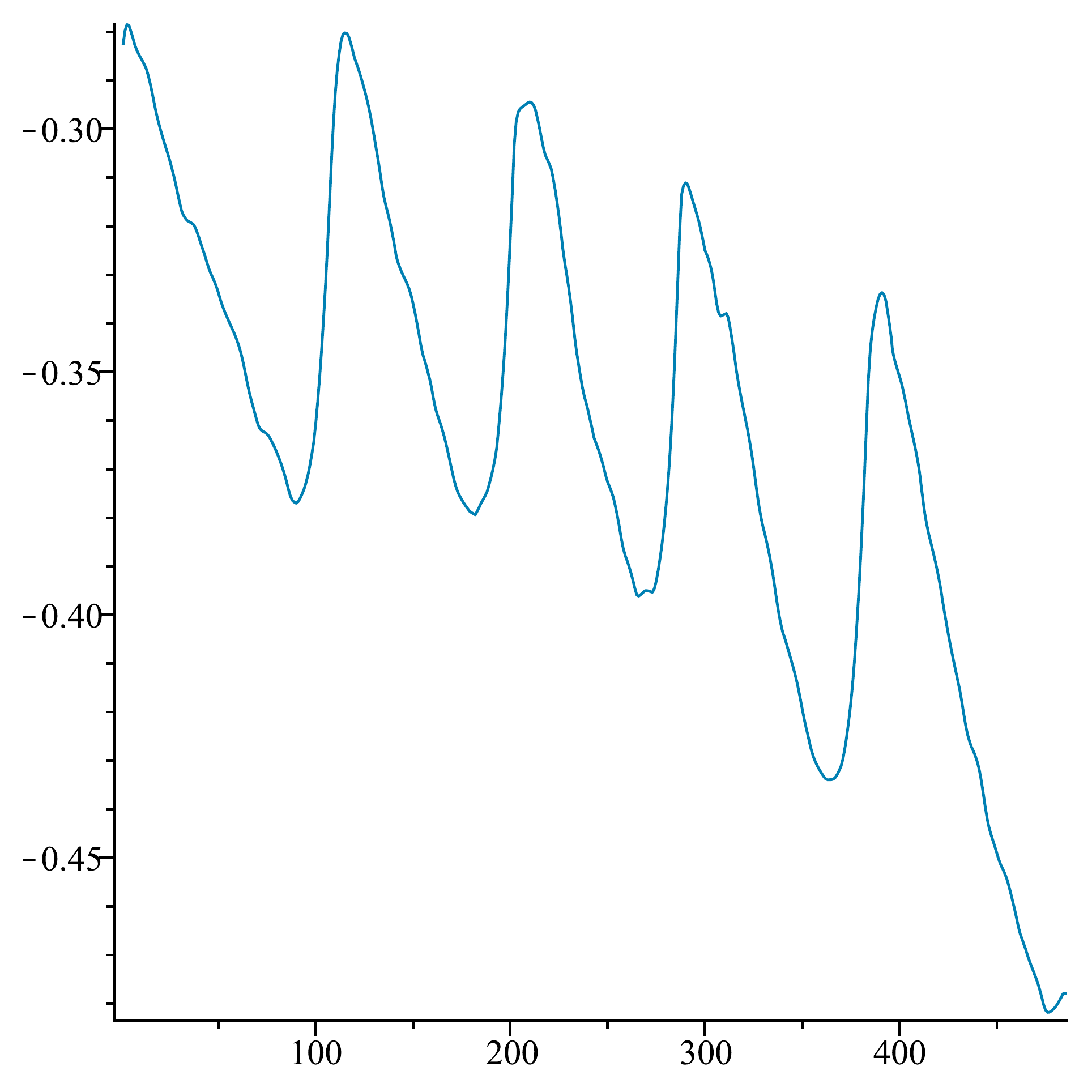}
\end{figure}

\begin{figure}[H]\caption{{\bf Messreihe I}, Gl\"attung mit SEA, Gleitl\"ange $10$.}\label{105}
 \includegraphics[width=0.45\textwidth]{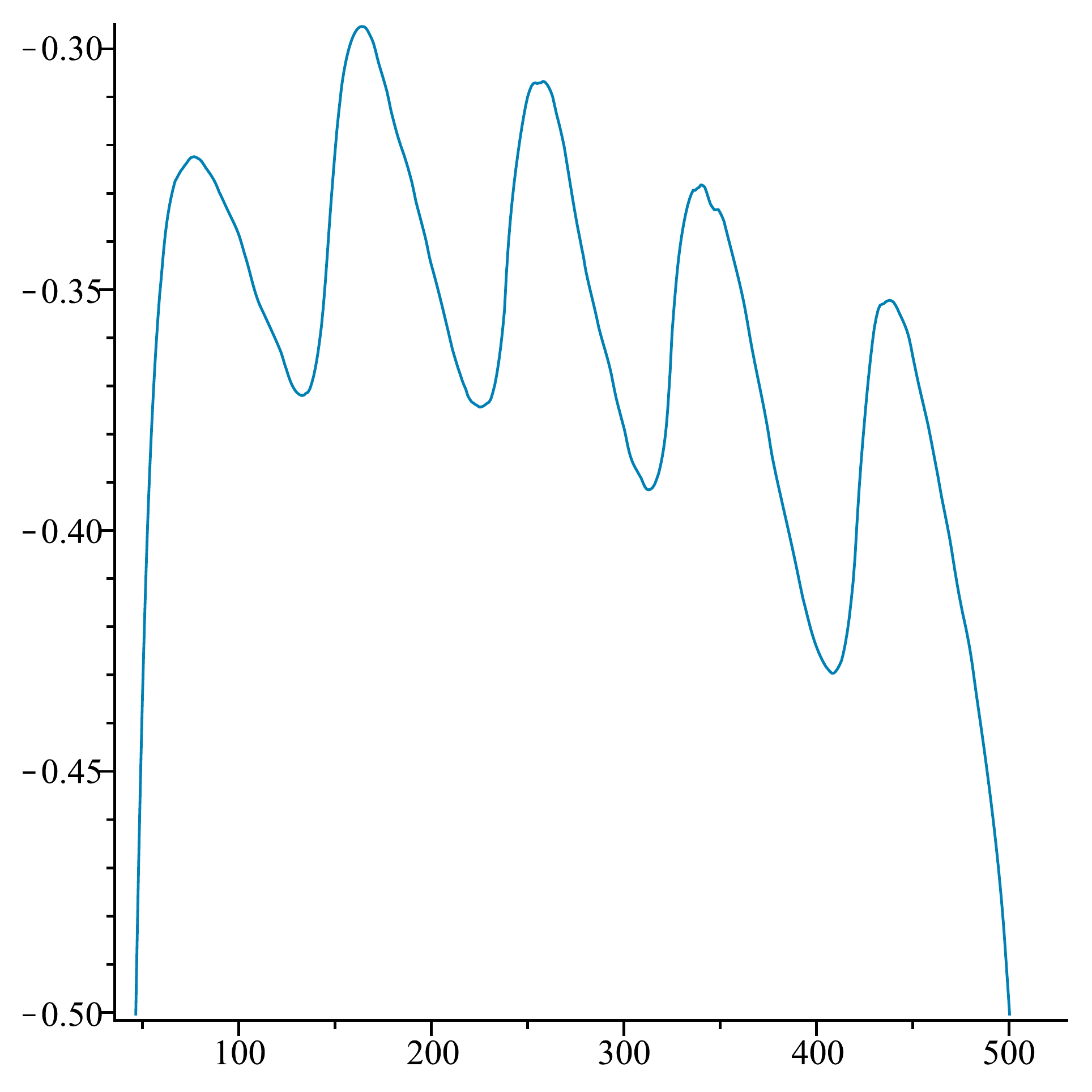}
 \includegraphics[width=0.45\textwidth]{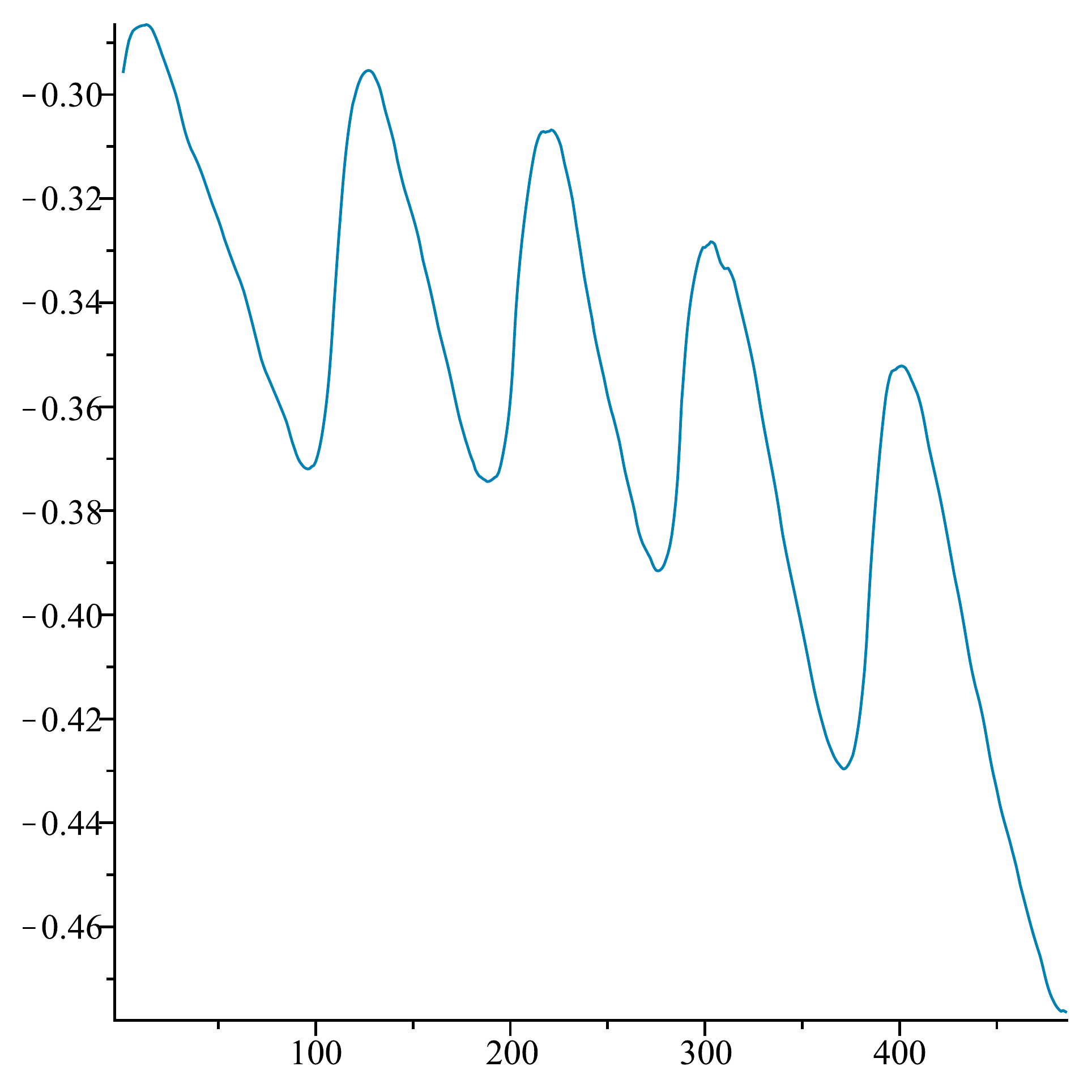}
\end{figure}

\begin{figure}[H]\caption{{\bf Messreihe I}, Originaldaten und Medianbildung mit Gleitl\"ange $9$.}\label{106}
 \includegraphics[width=0.45\textwidth]{real2-zoom-daten.pdf}
 \includegraphics[width=0.45\textwidth]{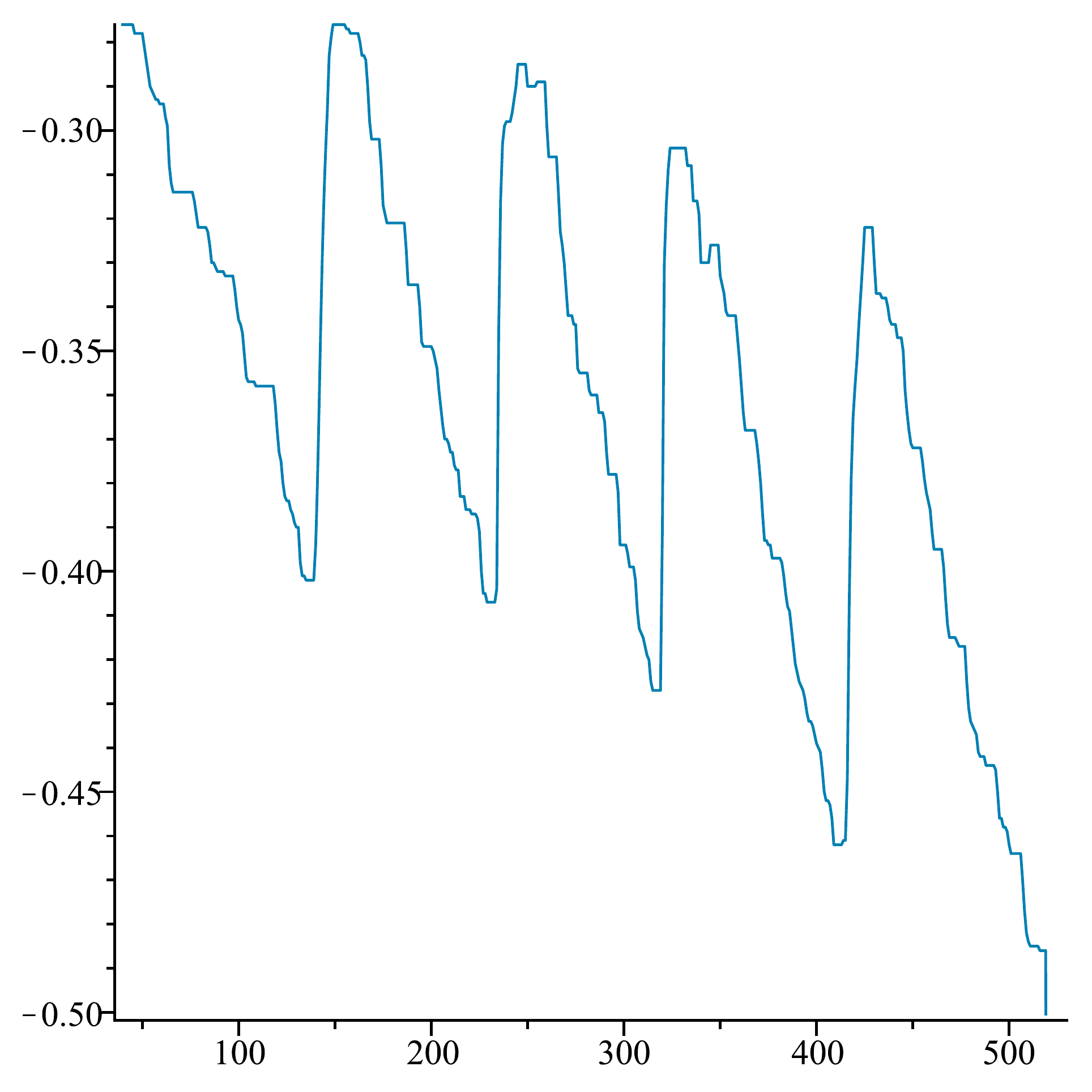}
\end{figure}

\begin{figure}[H]\caption{{\bf Messreihe II}, Originaldaten und bereinigte Originaldaten} \label{201}
 \includegraphics[width=0.45\textwidth]{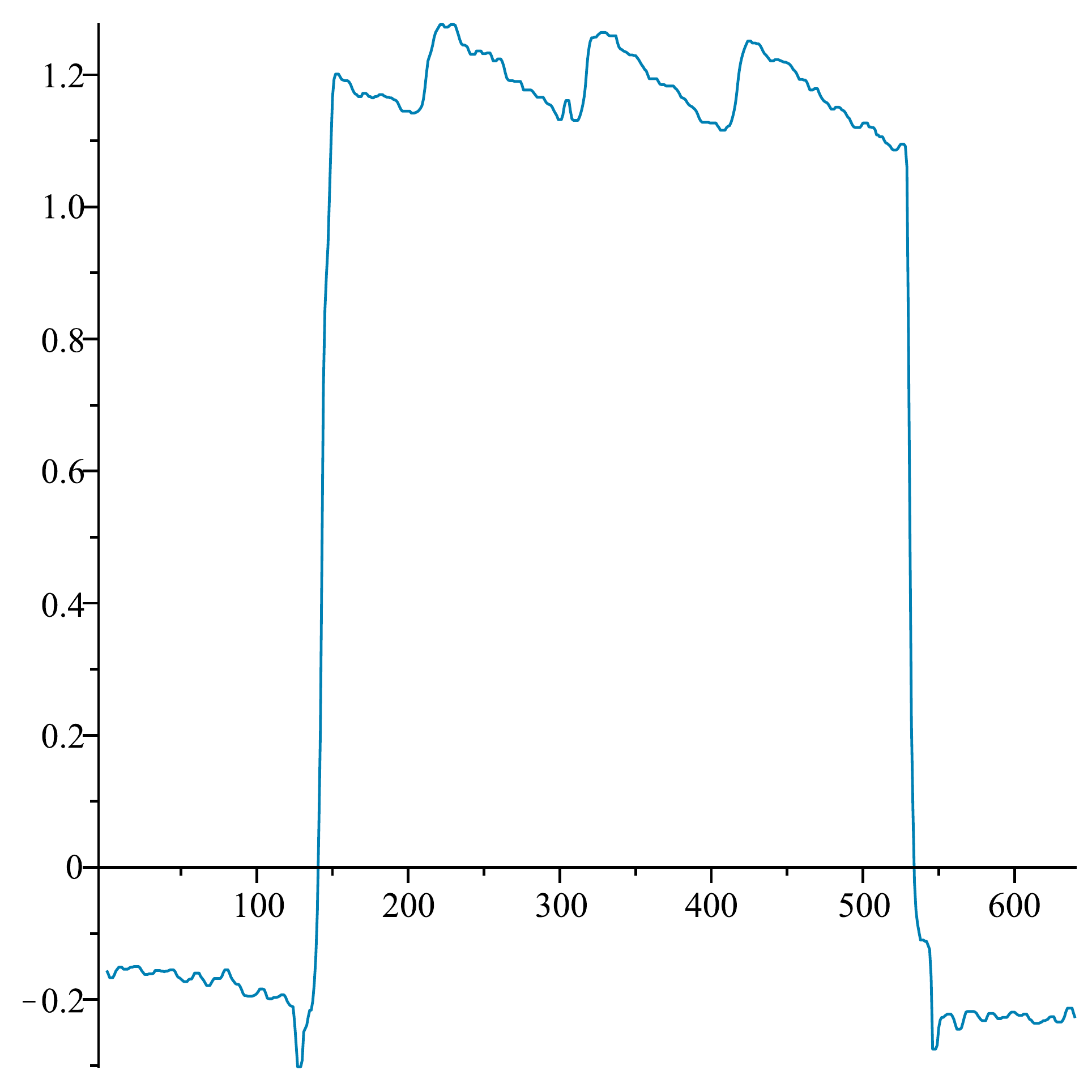}
 \includegraphics[width=0.45\textwidth]{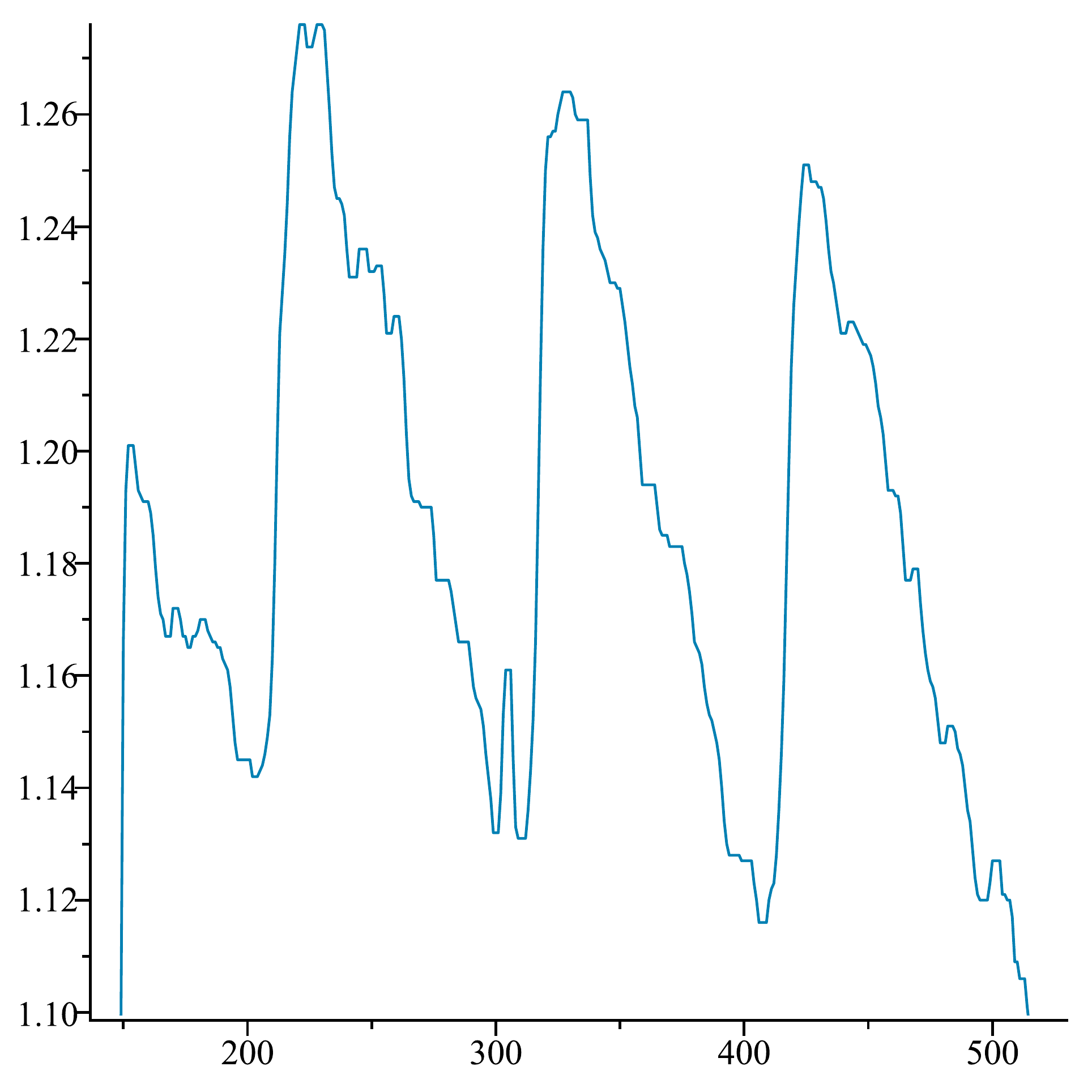}
\end{figure}

\begin{figure}[H]\caption{{\bf Messreihe II}, Gl\"attung mit MA, Gleitl\"ange $10$.} \label{202}
 \includegraphics[width=0.45\textwidth]{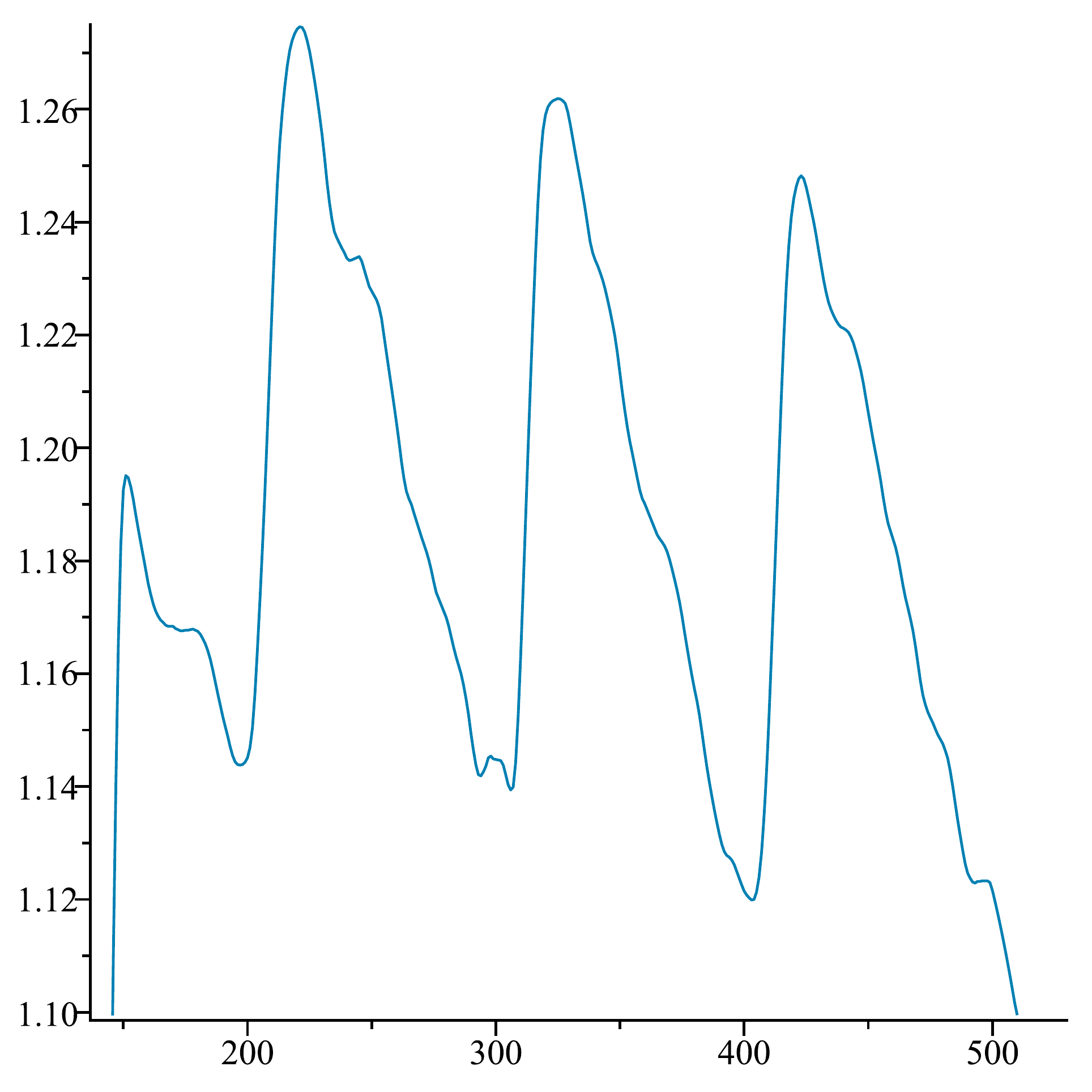}
 \includegraphics[width=0.45\textwidth]{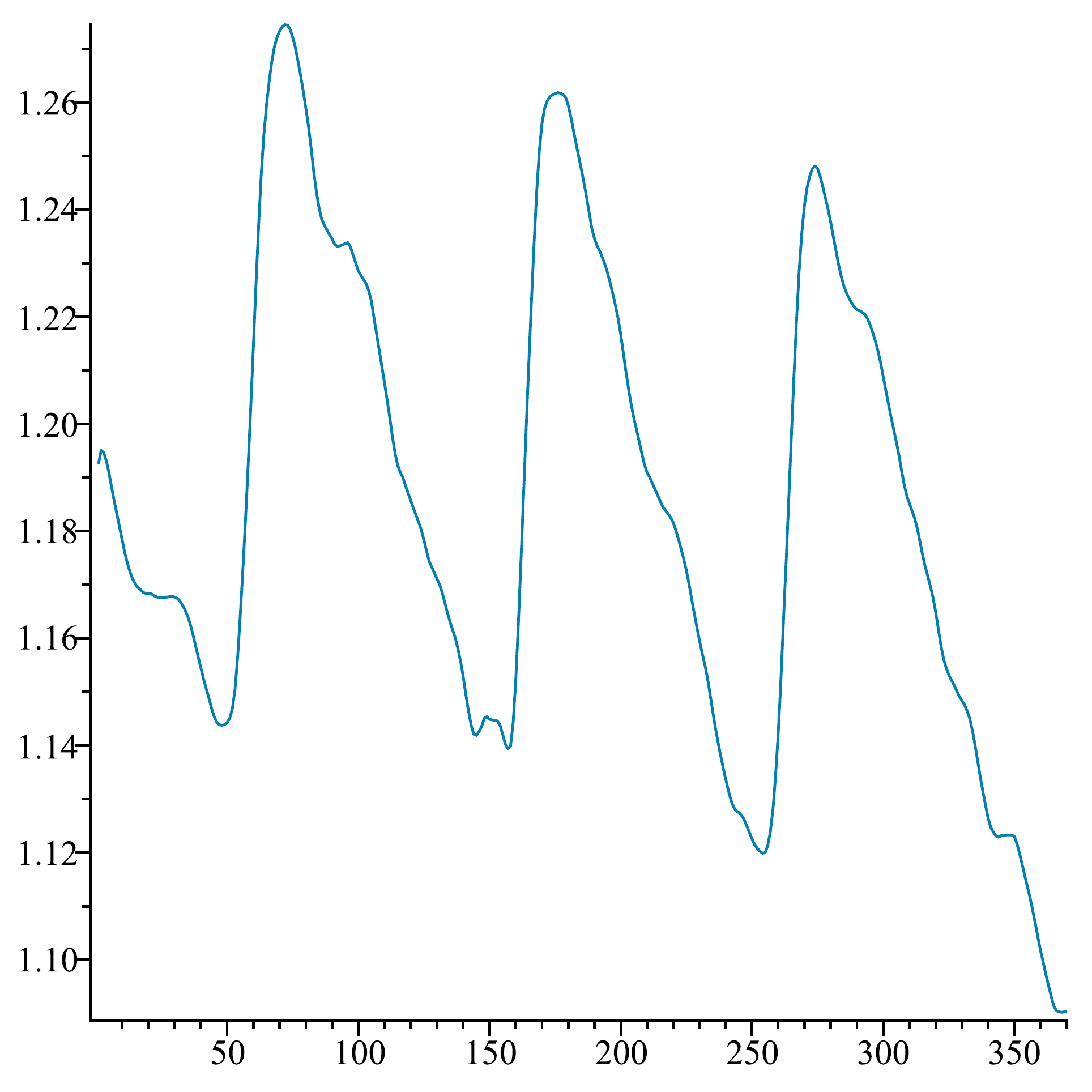}
\end{figure}

\begin{figure}[H]\caption{{\bf Messreihe II}, Gl\"attung mit EA, Gleitl\"ange $10$.}\label{203}
 \includegraphics[width=0.45\textwidth]{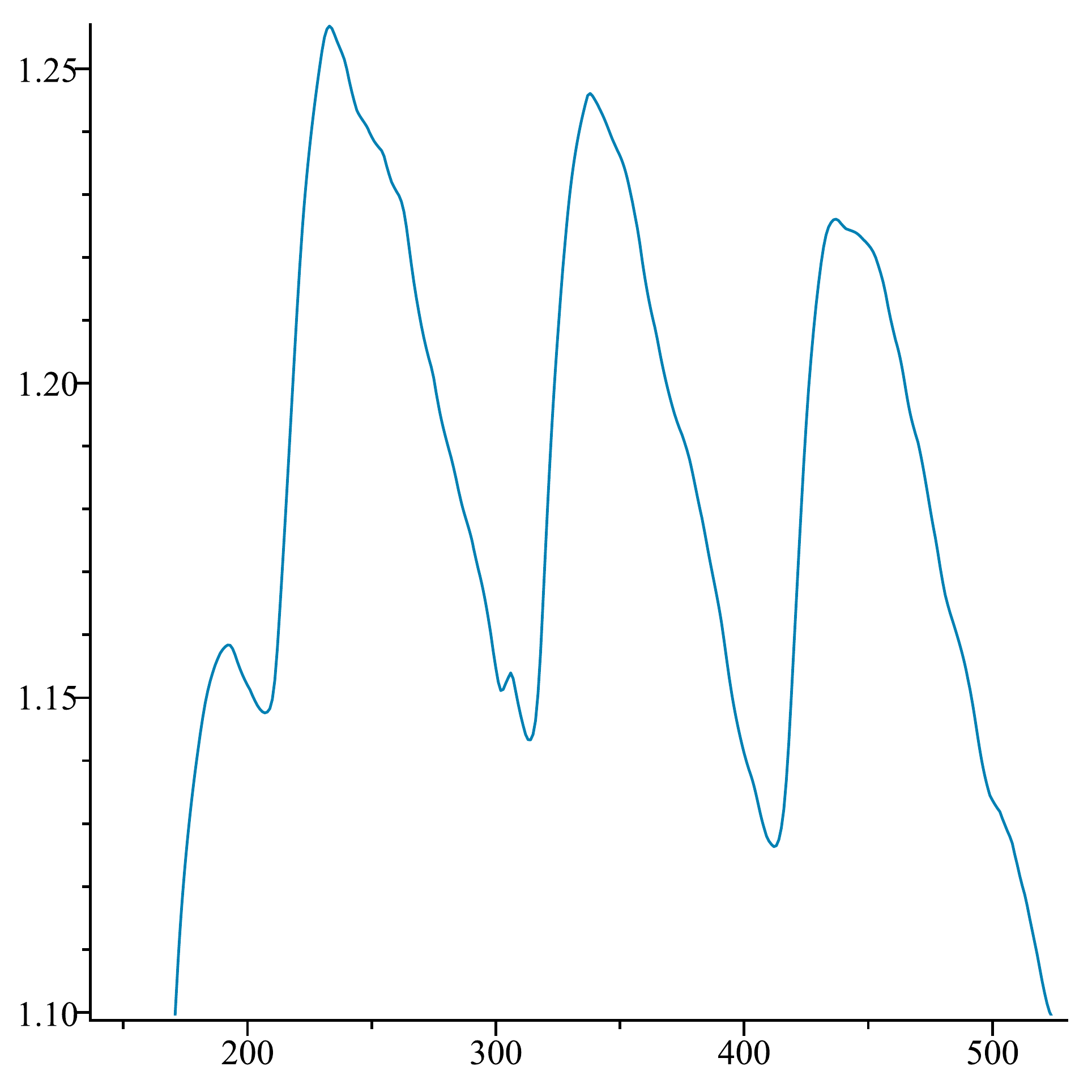}
 \includegraphics[width=0.45\textwidth]{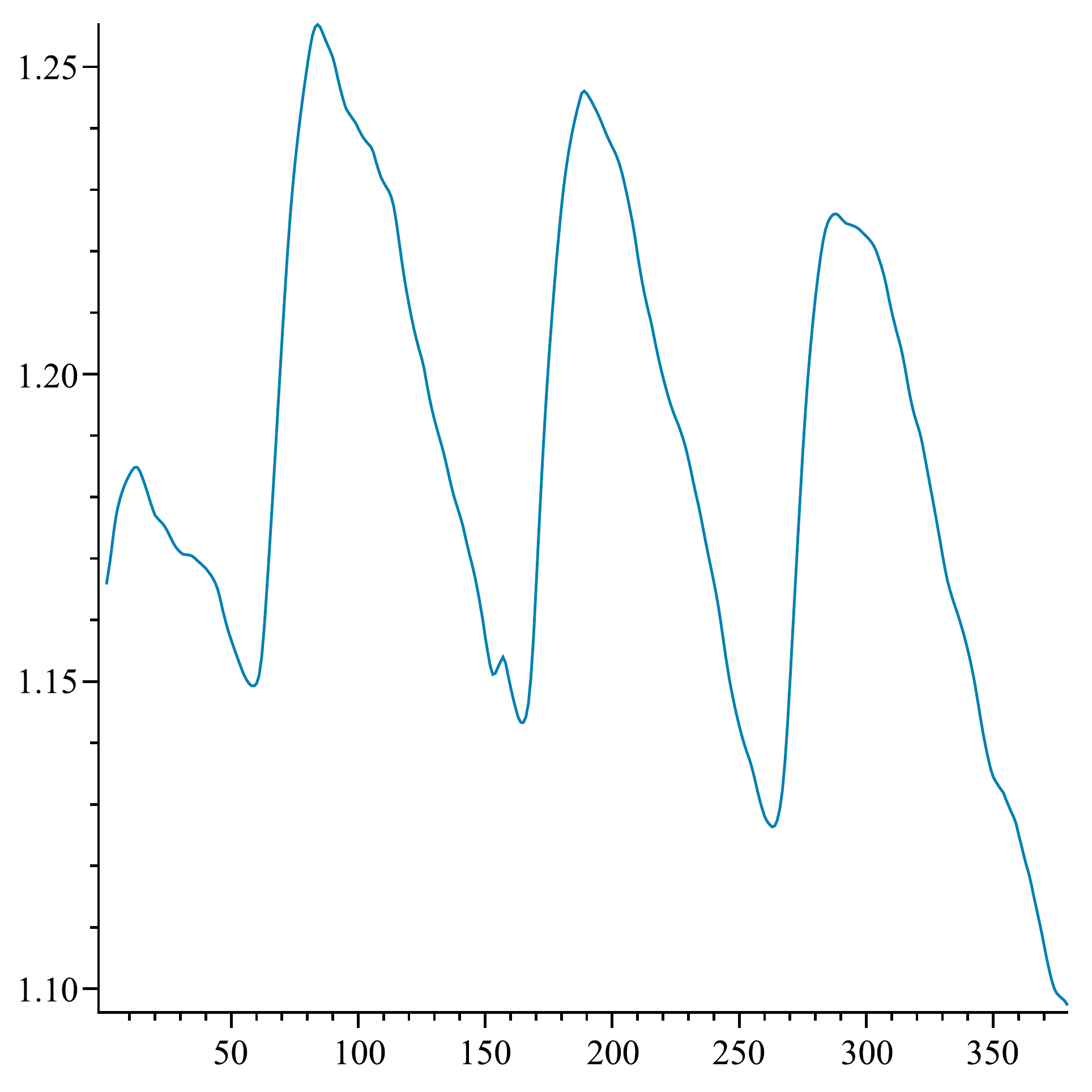}
\end{figure}

\begin{figure}[H]\caption{{\bf Messreihe II}, Gl\"attung mit rEA, Gleitl\"ange $10$.}\label{204}  
 \includegraphics[width=0.45\textwidth]{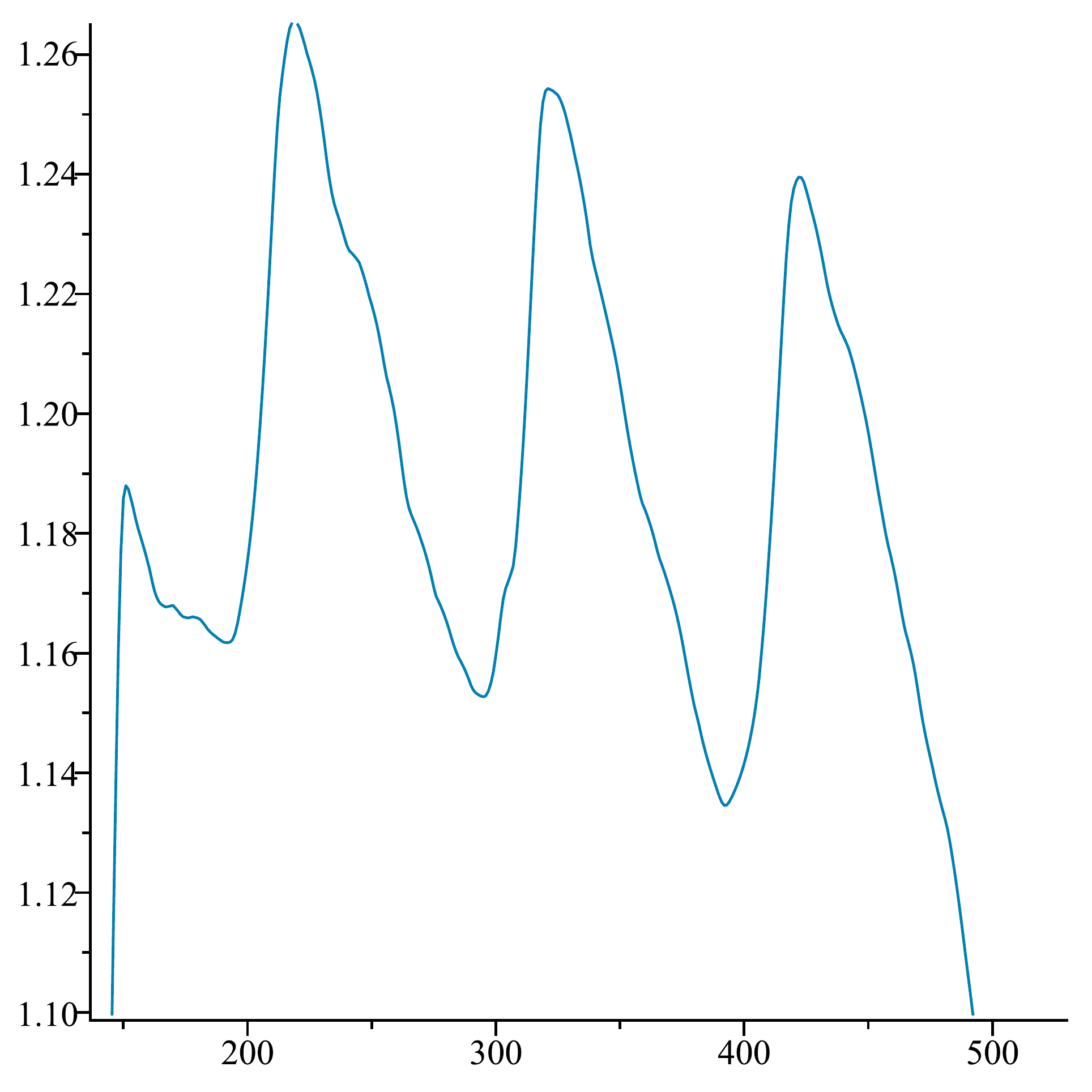}
 \includegraphics[width=0.45\textwidth]{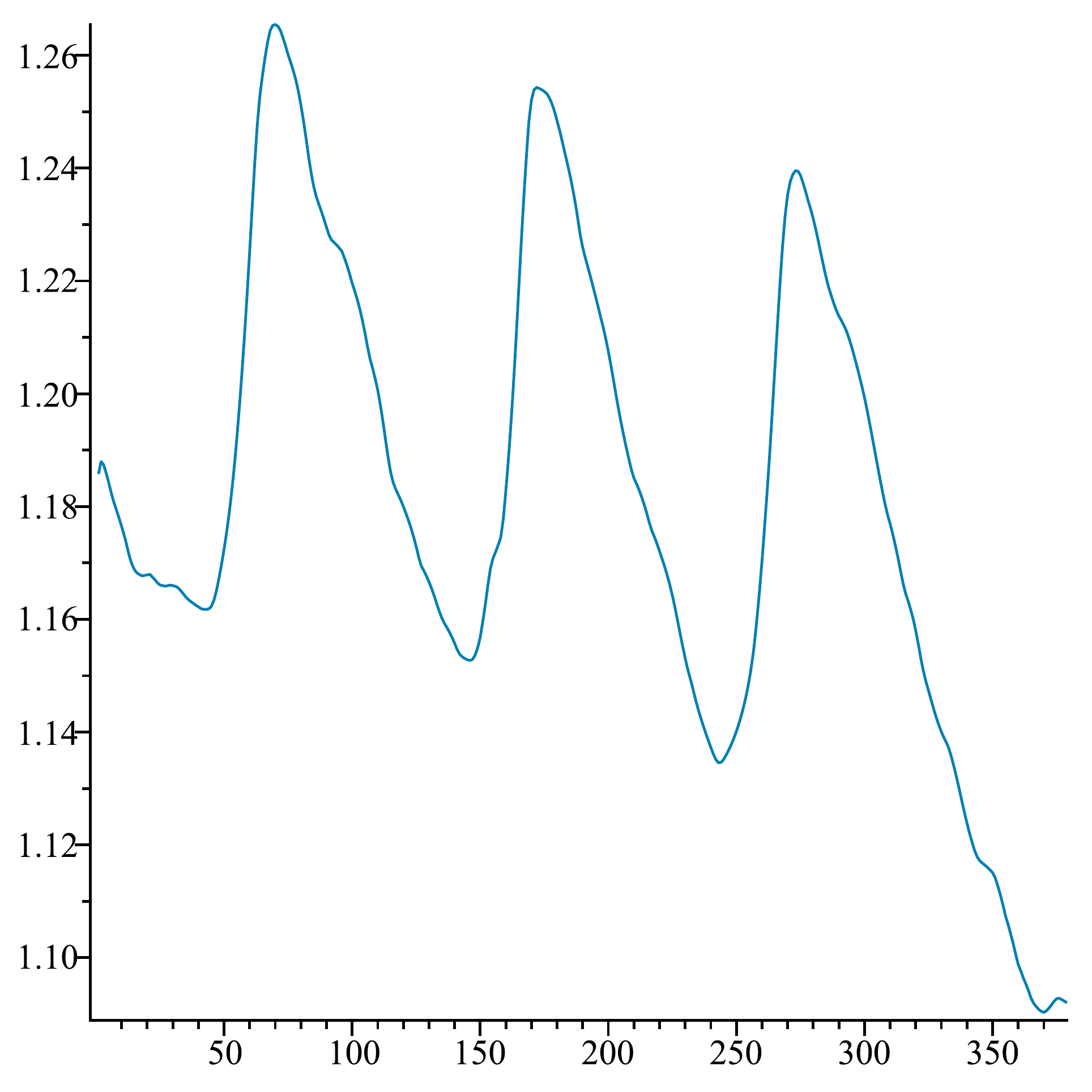}
\end{figure}

\begin{figure}[H]\caption{\small {\bf Messreihe II}, Gl\"attung mit SEA, Gleitl\"ange $10$.}\label{205}
 \includegraphics[width=0.45\textwidth]{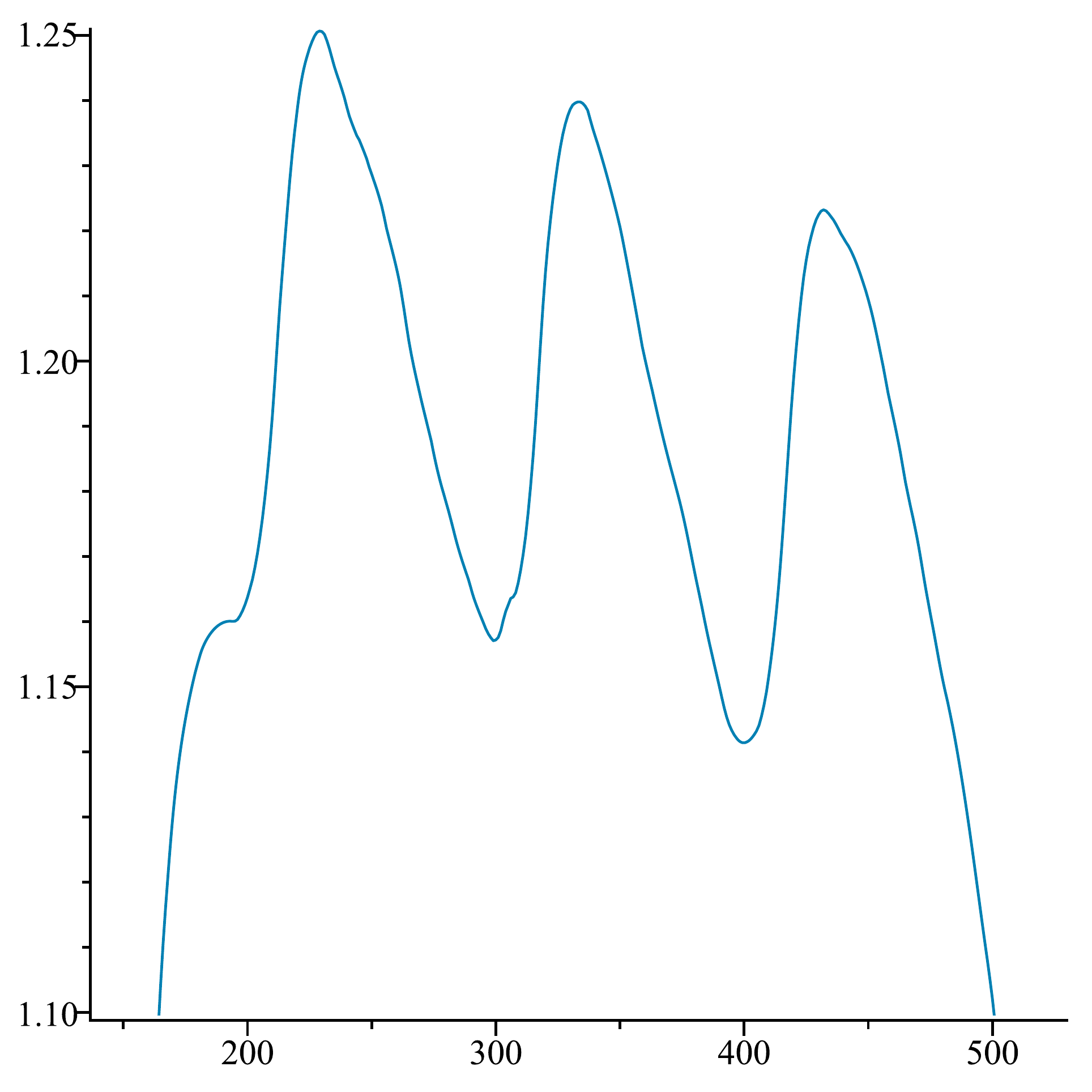}
 \includegraphics[width=0.45\textwidth]{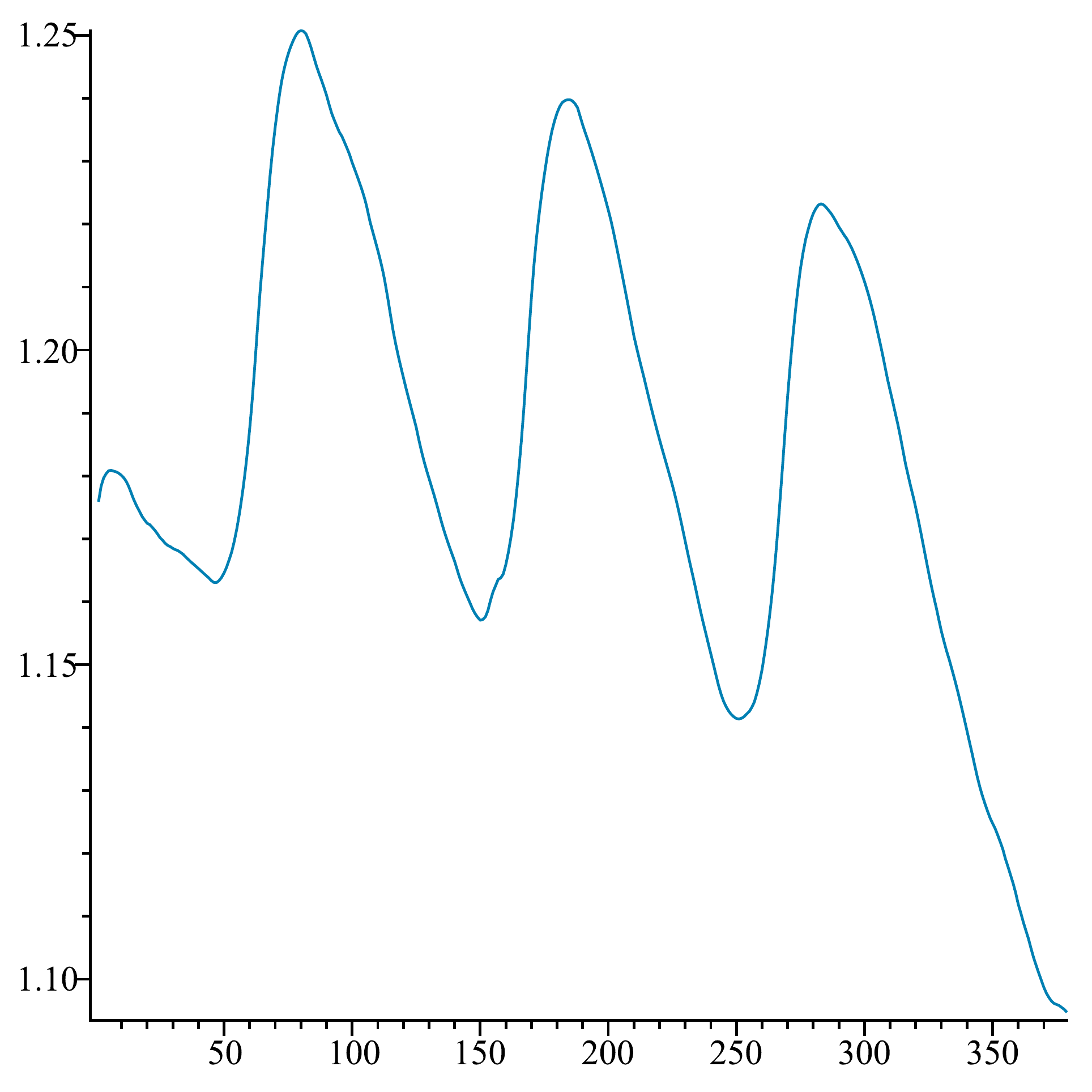}
\end{figure}

\begin{figure}[H]\caption{{\bf Messreihe III}, Originaldaten und bereinigte Originaldaten} \label{401}
 \includegraphics[width=0.45\textwidth]{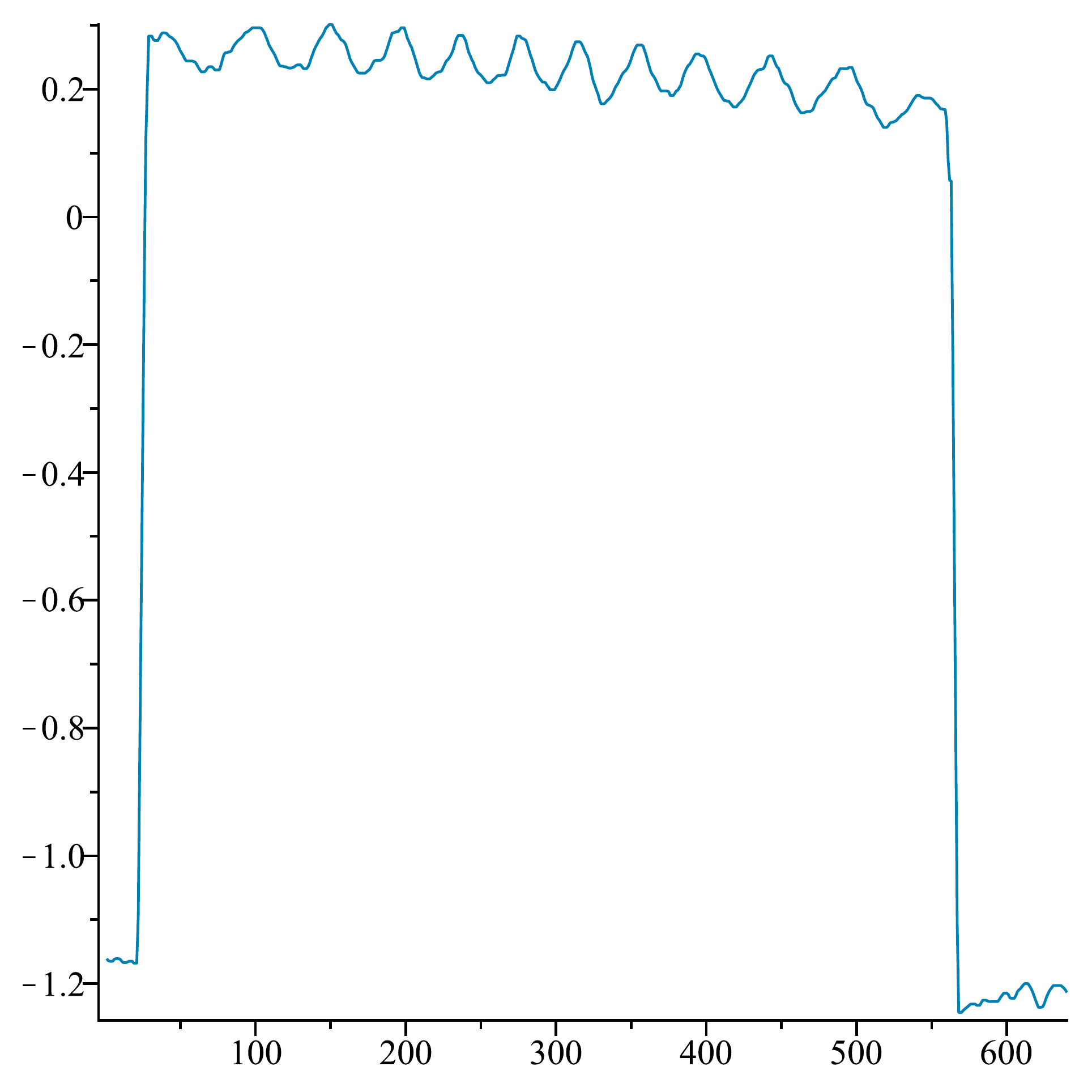}
 \includegraphics[width=0.49\textwidth]{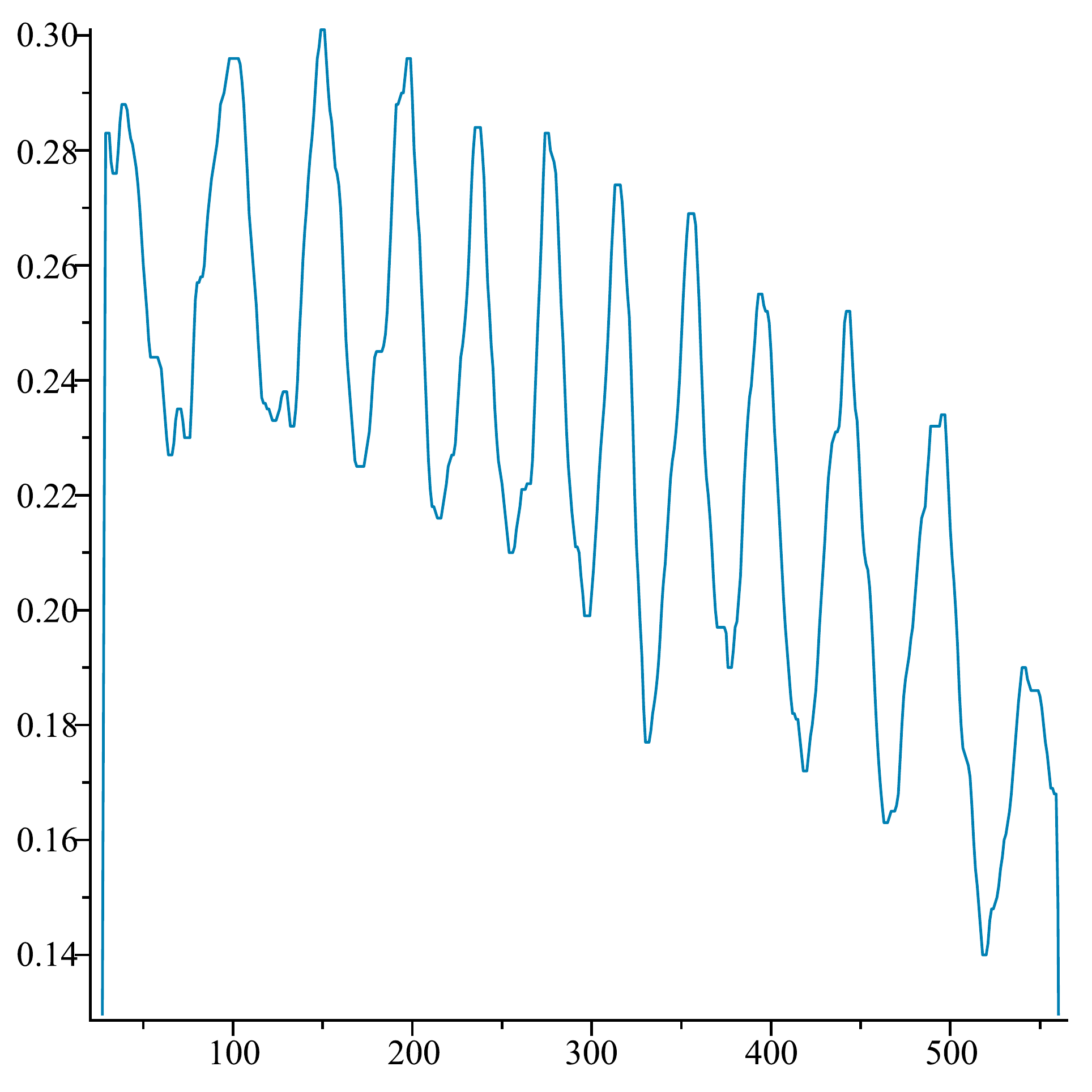}
\end{figure}

\begin{figure}[H]\caption{{\bf Messreihe III}, Gl\"attung mit MA, Gleitl\"ange $10$.} \label{402}
 \includegraphics[width=0.45\textwidth]{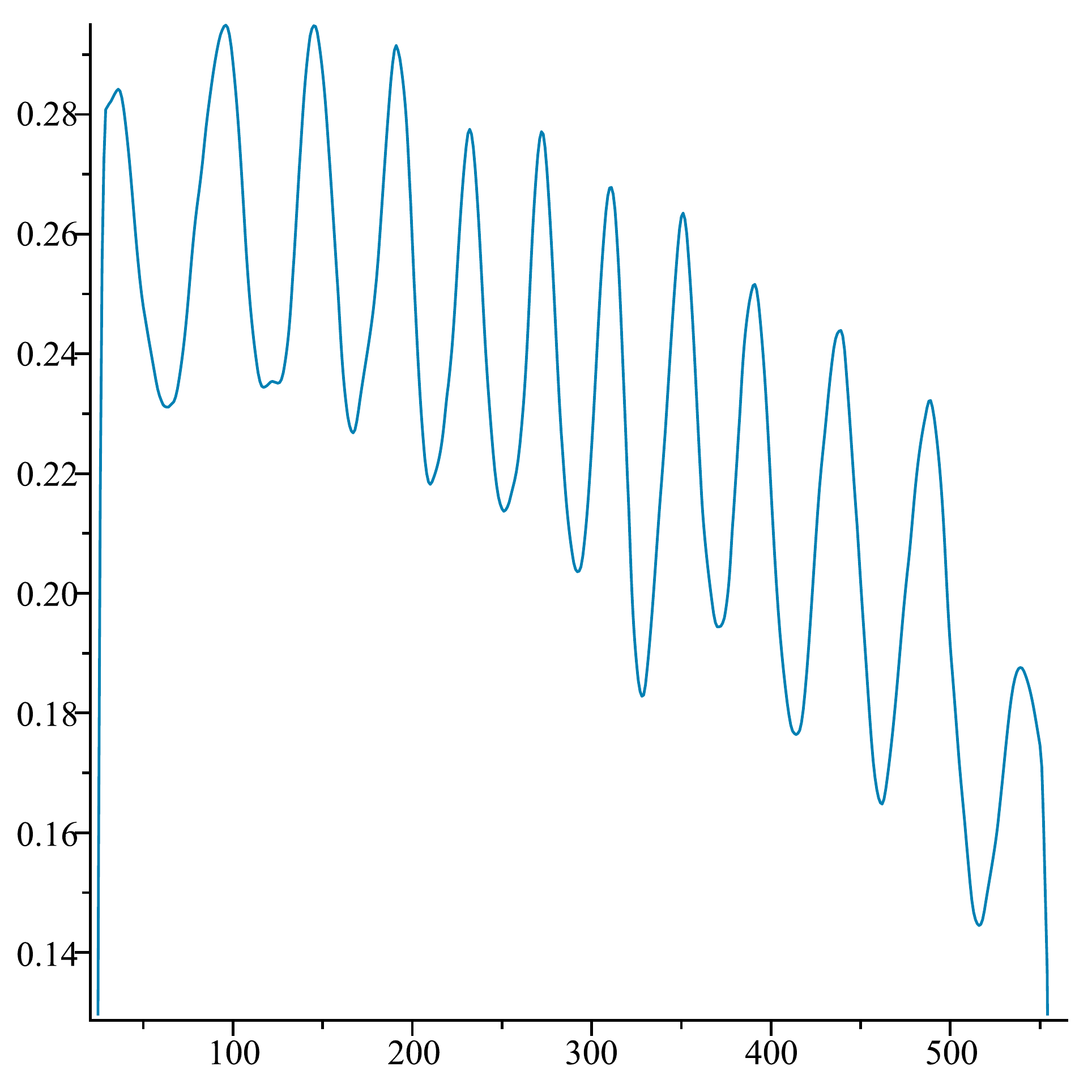}
 \includegraphics[width=0.45\textwidth]{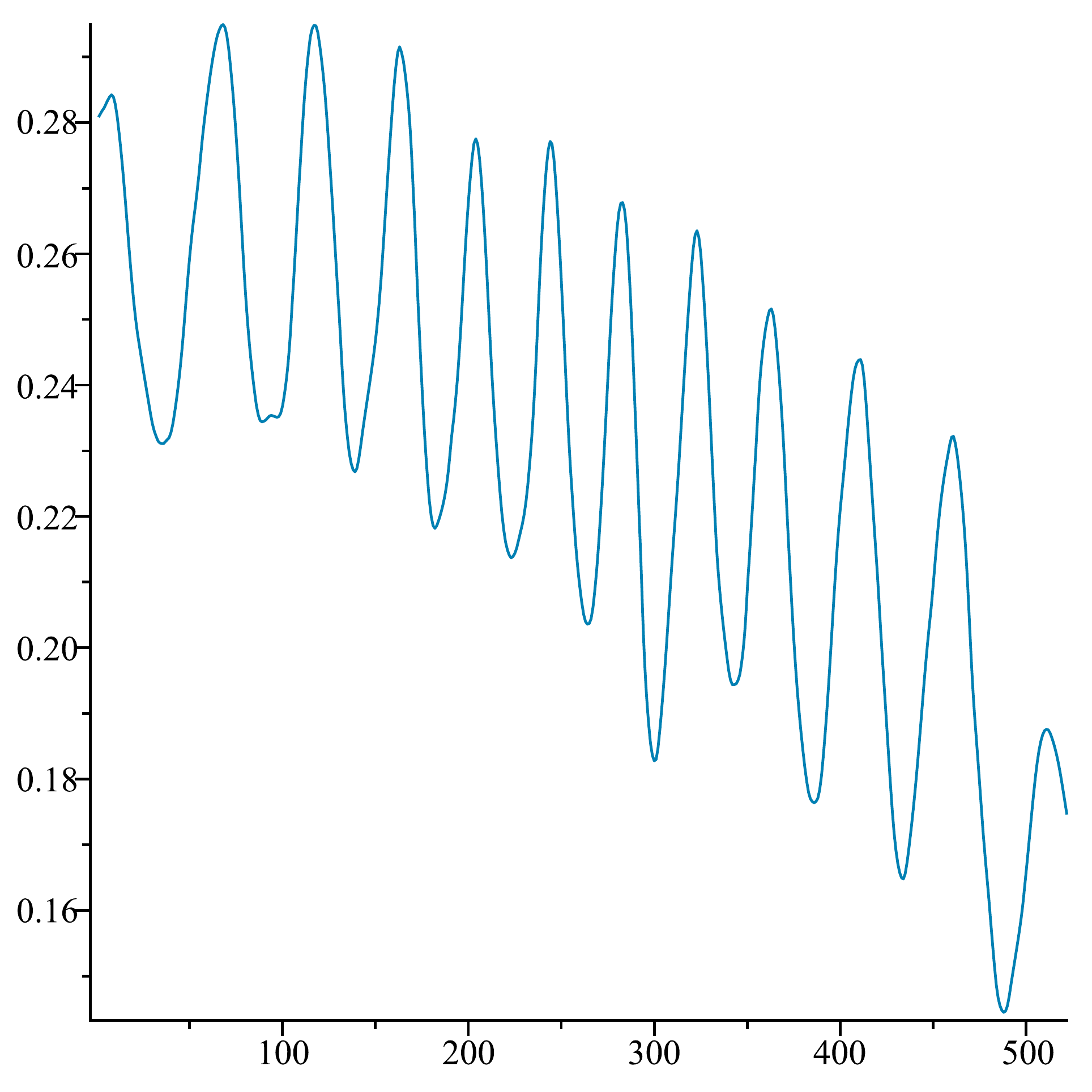}
\end{figure}

\begin{figure}[H]\caption{{\bf Messreihe III}, Gl\"attung mit EA, Gleitl\"ange $10$.}\label{403}
 \includegraphics[width=0.45\textwidth]{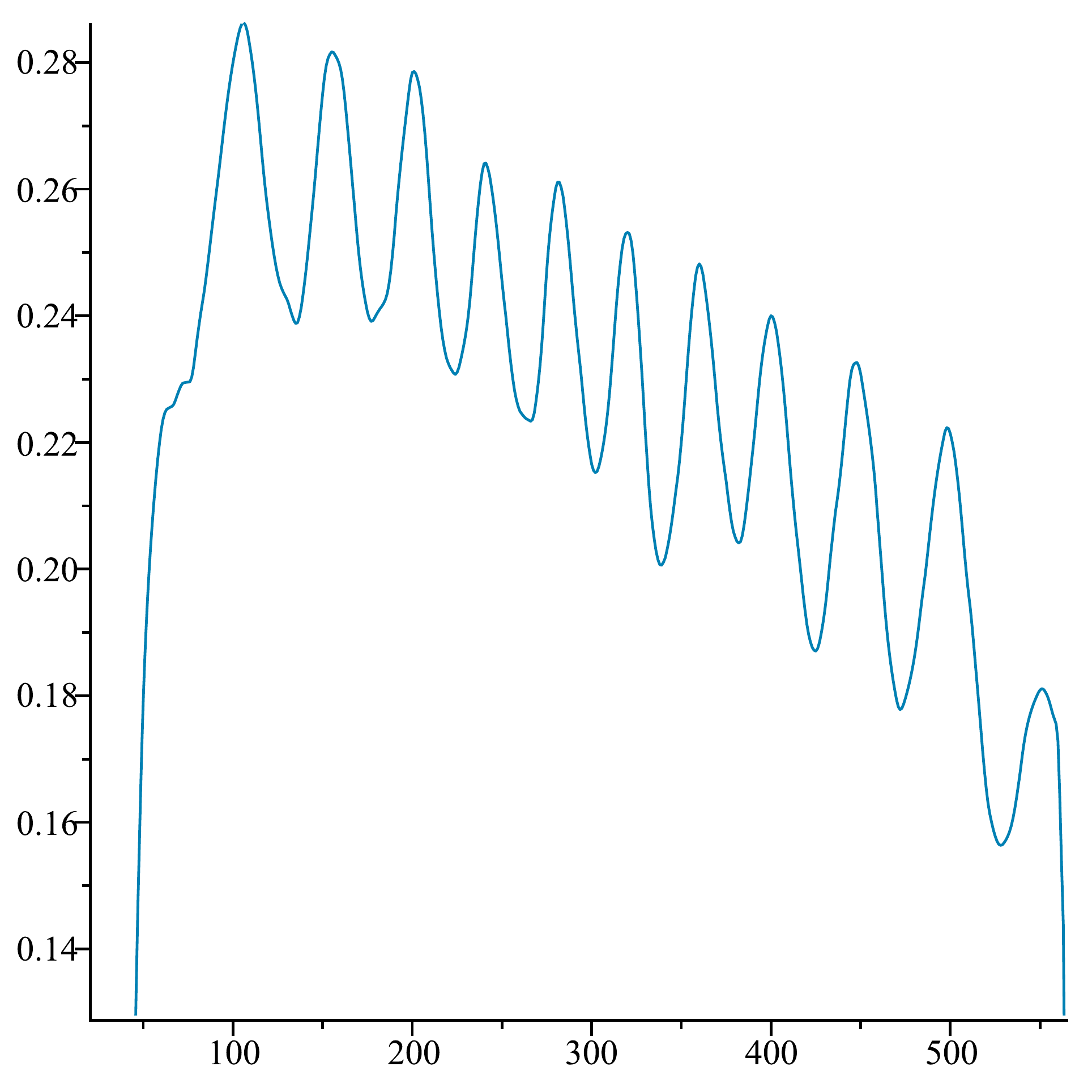}
 \includegraphics[width=0.45\textwidth]{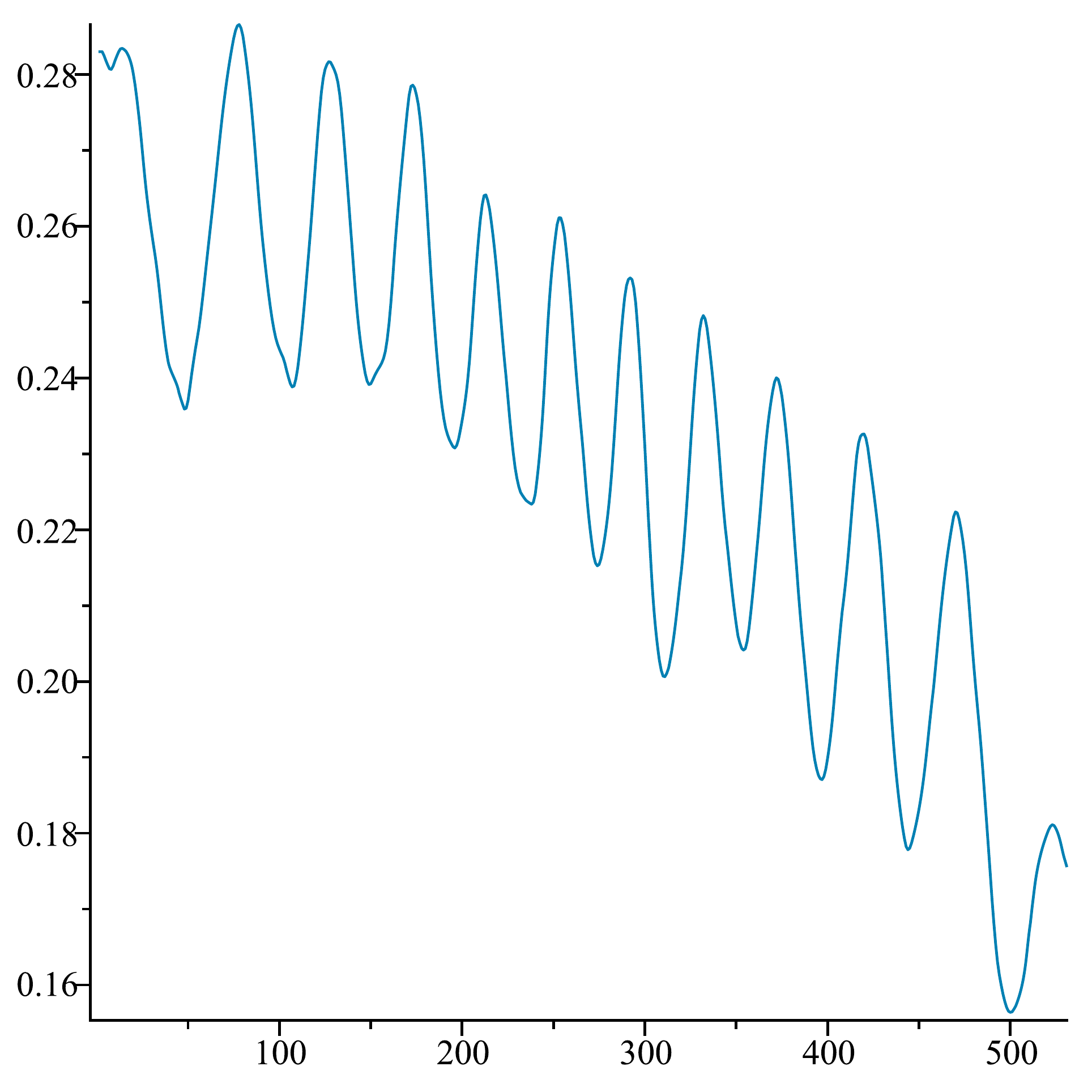}
\end{figure}

\begin{figure}[H]\caption{{\bf Messreihe III}, Gl\"attung mit rEA, Gleitl\"ange $10$.}\label{404}
 \includegraphics[width=0.45\textwidth]{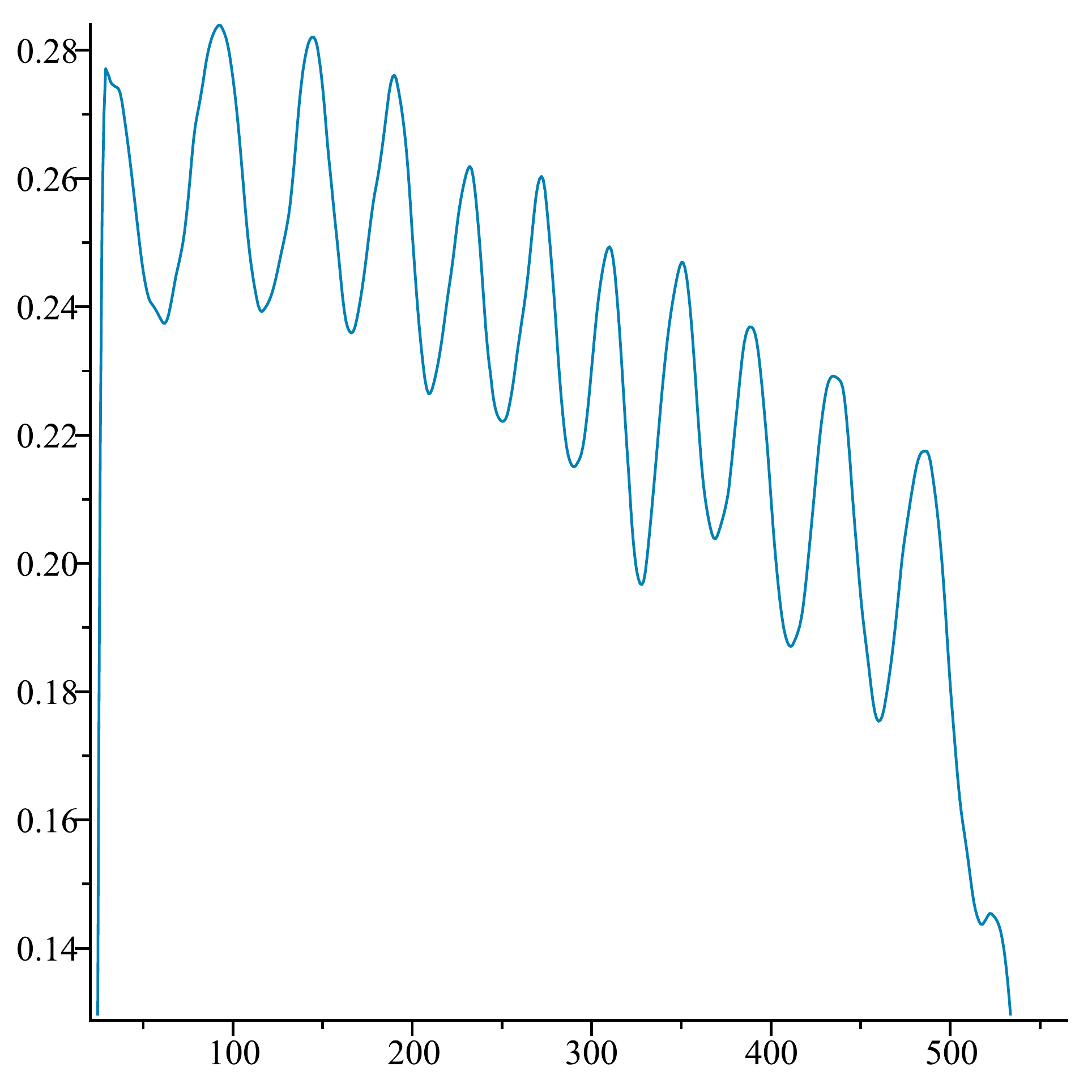}
 \includegraphics[width=0.45\textwidth]{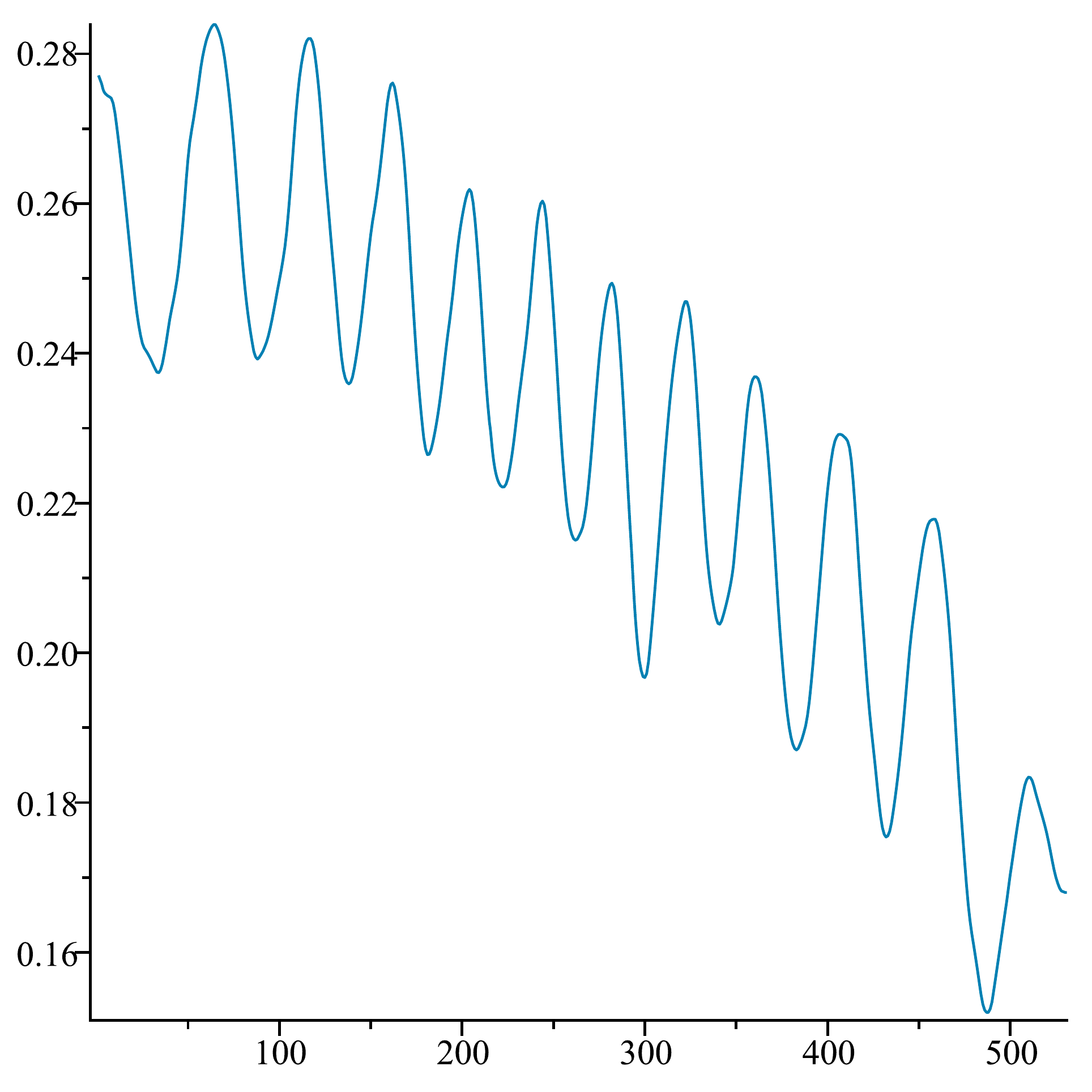}
\end{figure}

\begin{figure}[H]\caption{{\bf Messreihe III}, Gl\"attung mit SEA, Gleitl\"ange $10$.}\label{405}
 \includegraphics[width=0.45\textwidth]{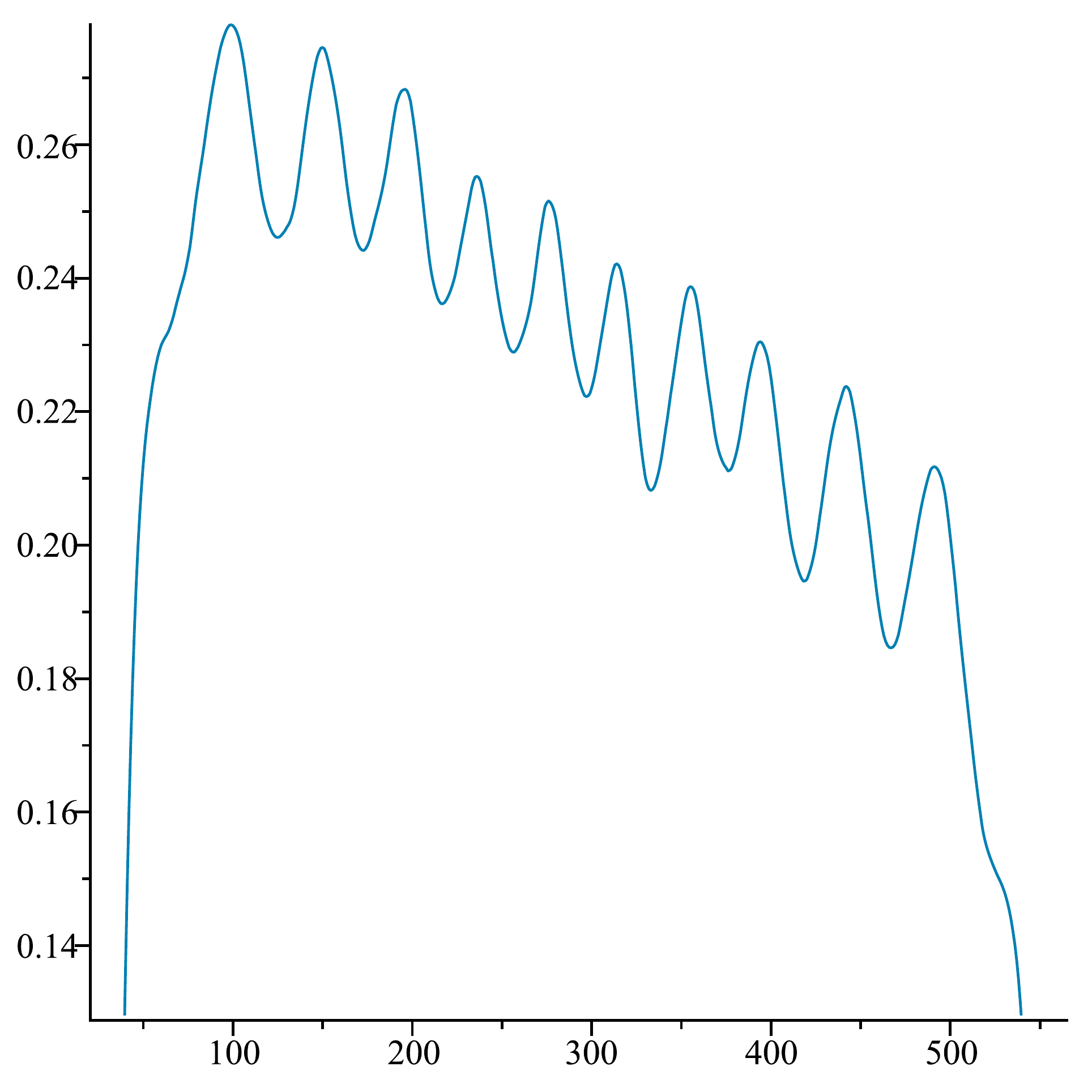}
 \includegraphics[width=0.45\textwidth]{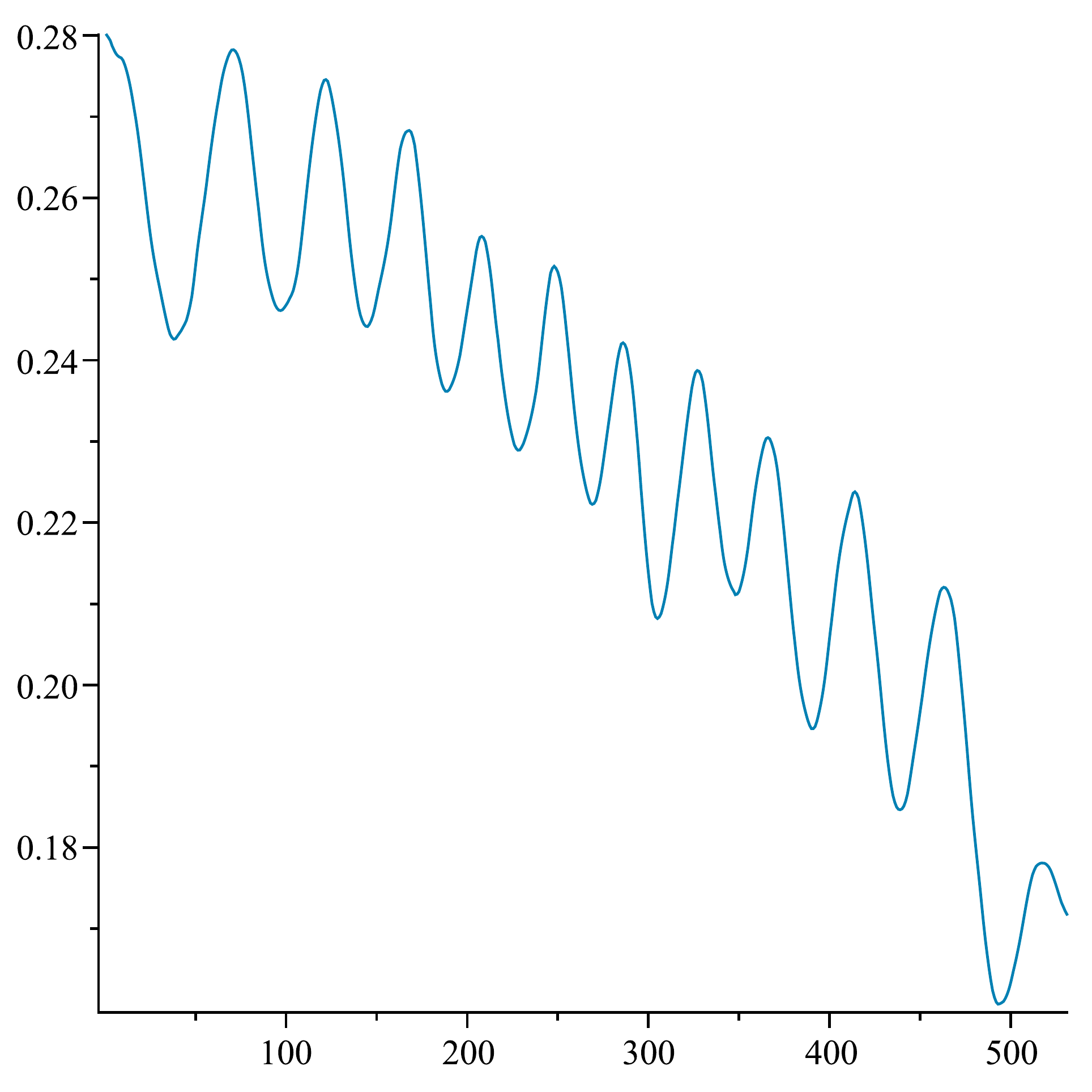}
\end{figure}

\begin{figure}[H]\caption{{\bf Messreihe IV}, Originaldaten und bereinigte Originaldaten} \label{501}
 \includegraphics[width=0.45\textwidth]{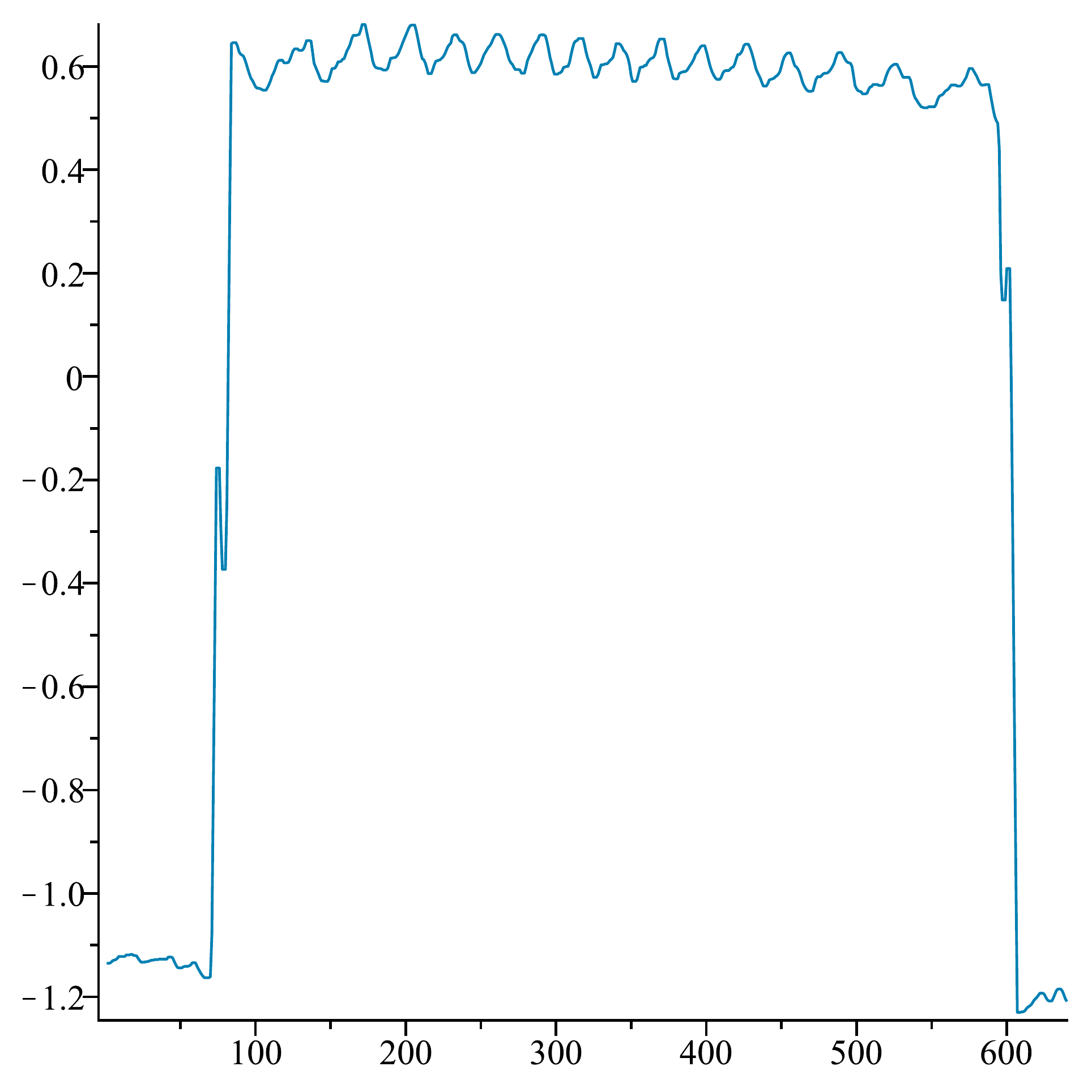}
 \includegraphics[width=0.45\textwidth]{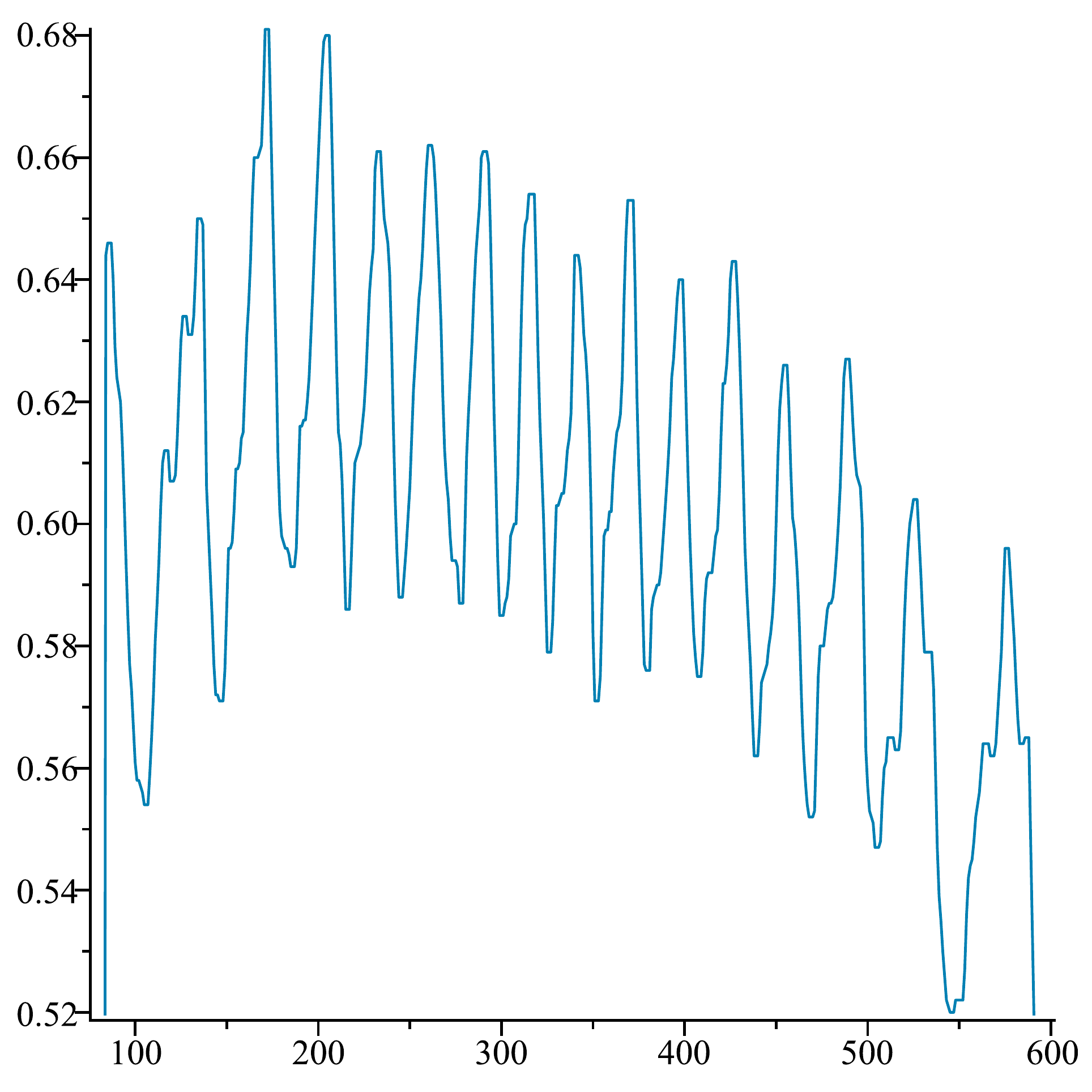}
\end{figure}

\begin{figure}[H]\caption{{\bf Messreihe IV}, Gl\"attung mit MA, Gleitl\"ange $10$.} \label{502}
 \includegraphics[width=0.45\textwidth]{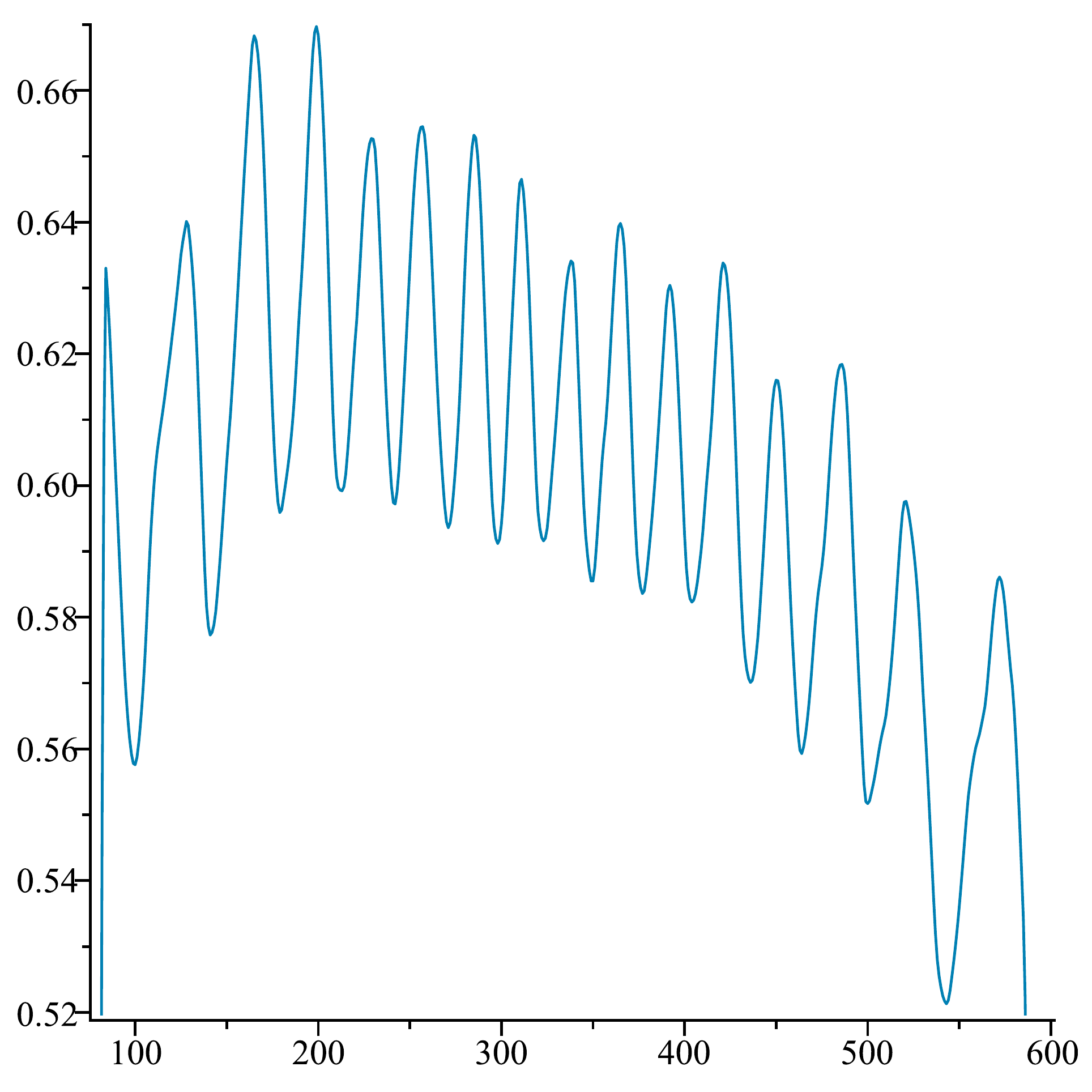}
 \includegraphics[width=0.45\textwidth]{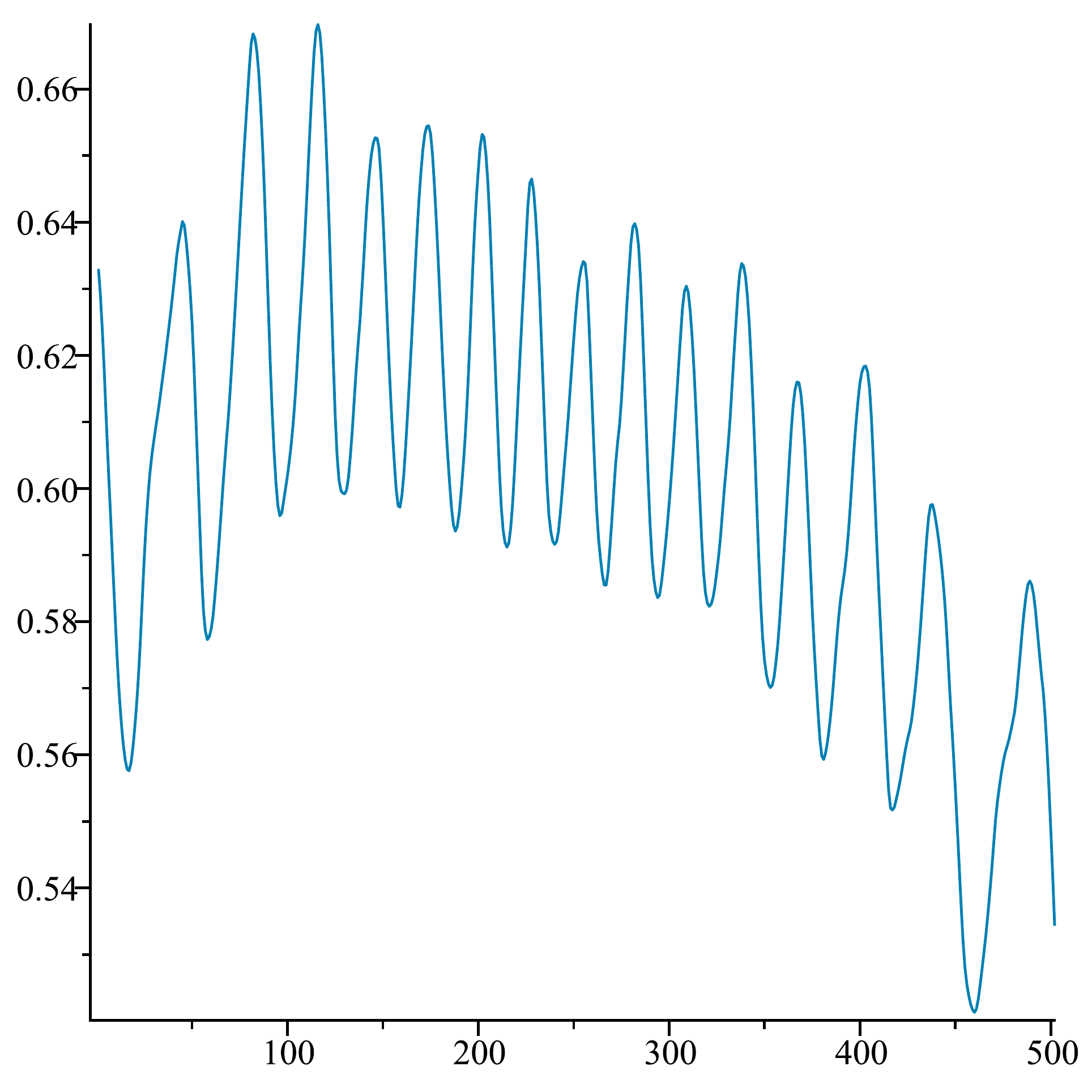}
\end{figure}

\begin{figure}[H]\caption{{\bf Messreihe IV}, Gl\"attung mit EA, Gleitl\"ange $10$.}\label{503}
 \includegraphics[width=0.45\textwidth]{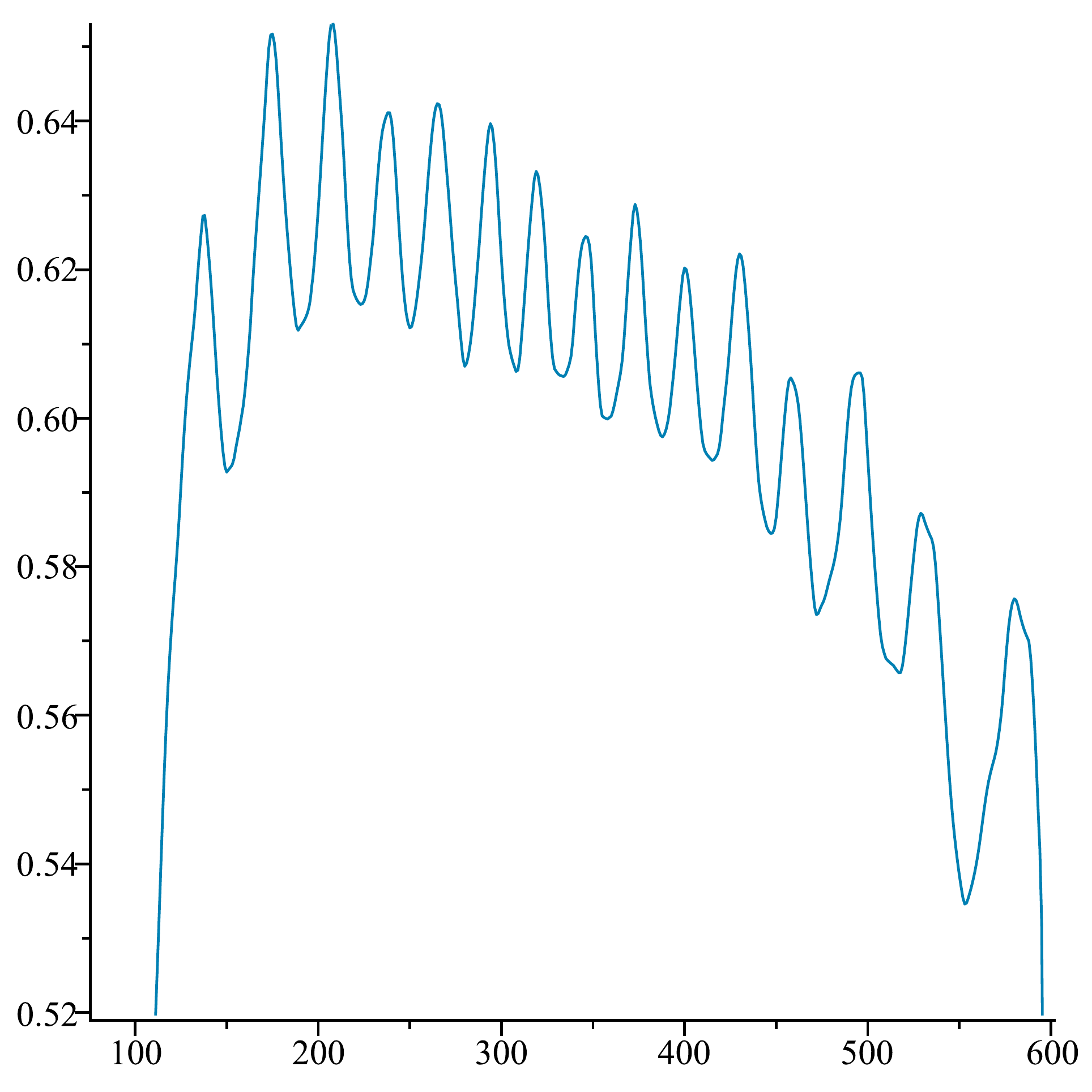}
 \includegraphics[width=0.45\textwidth]{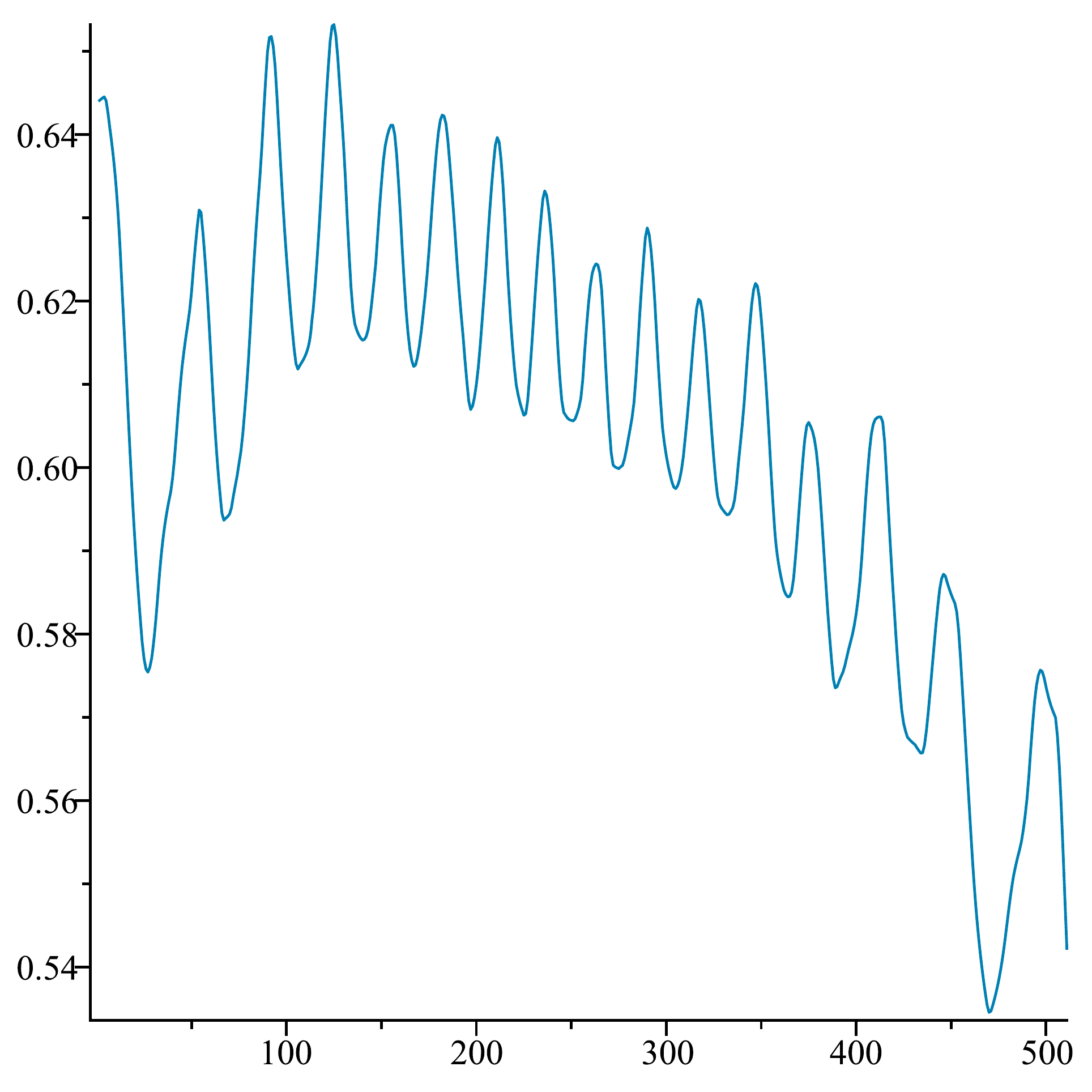}
\end{figure}

\begin{figure}[H]\caption{{\bf Messreihe IV}, Gl\"attung mit rEA, Gleitl\"ange $10$.}\label{504}
 \includegraphics[width=0.45\textwidth]{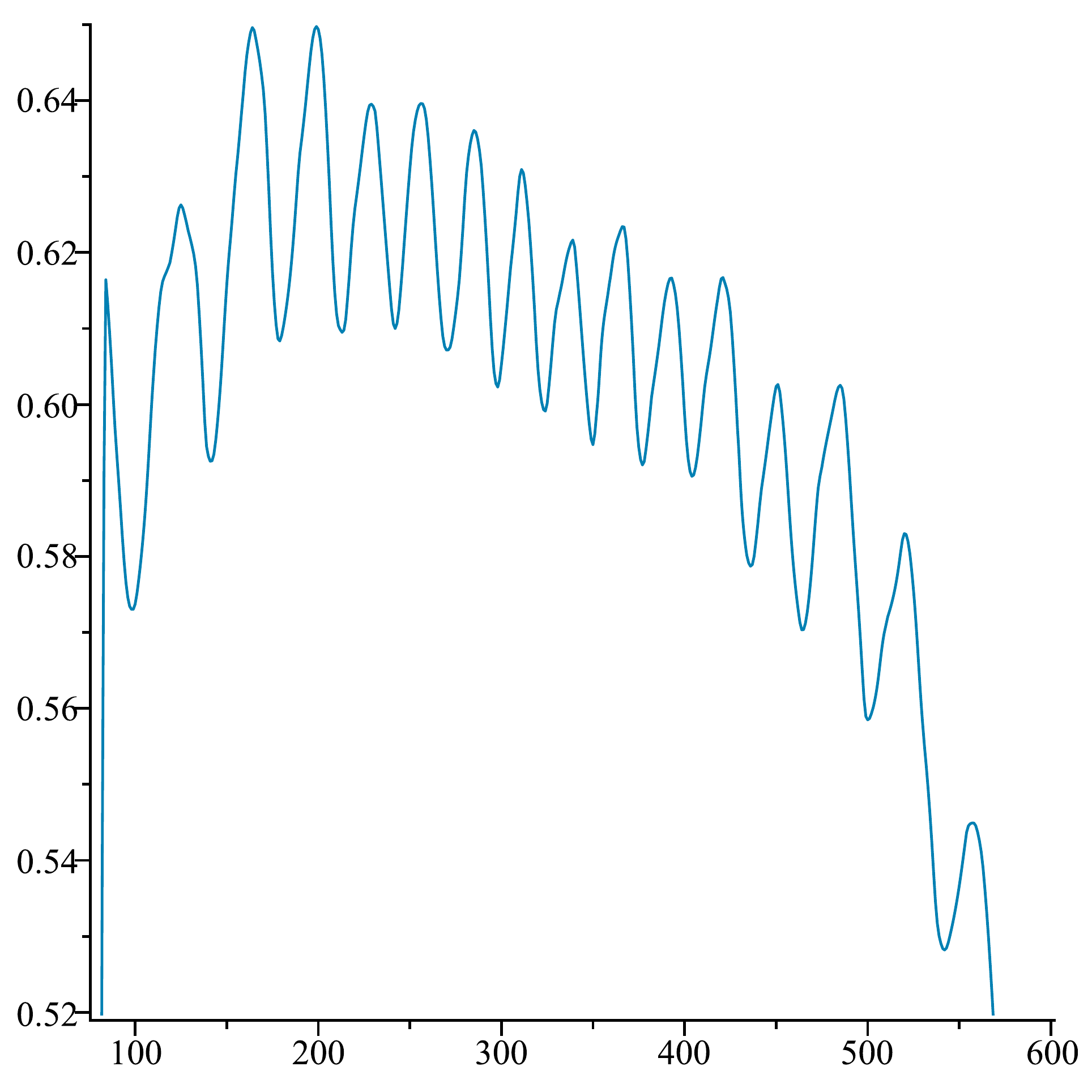}
 \includegraphics[width=0.45\textwidth]{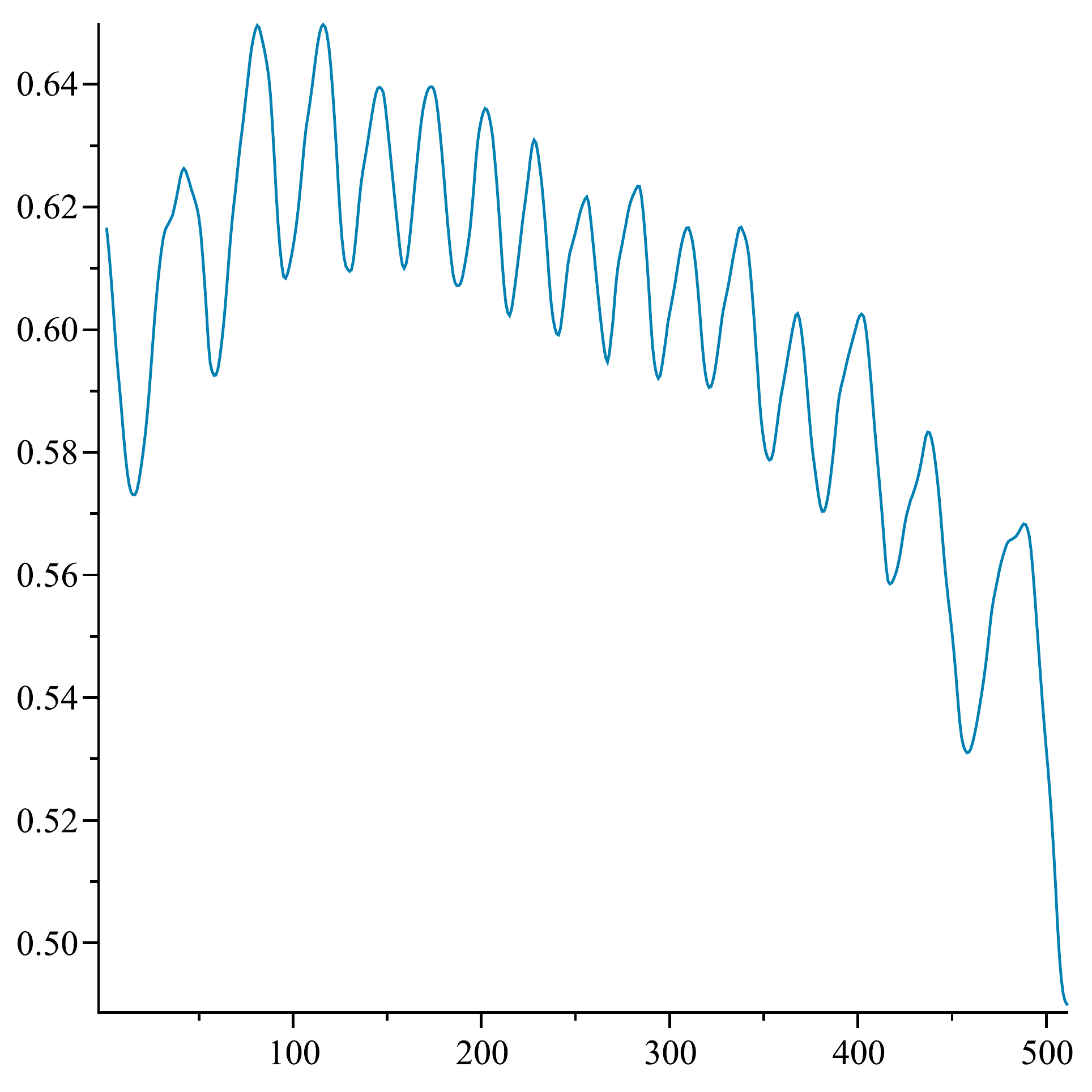}
\end{figure}

\begin{figure}[H]\caption{{\bf Messreihe IV}, Gl\"attung mit SEA, Gleitl\"ange $10$.}\label{505}
 \includegraphics[width=0.45\textwidth]{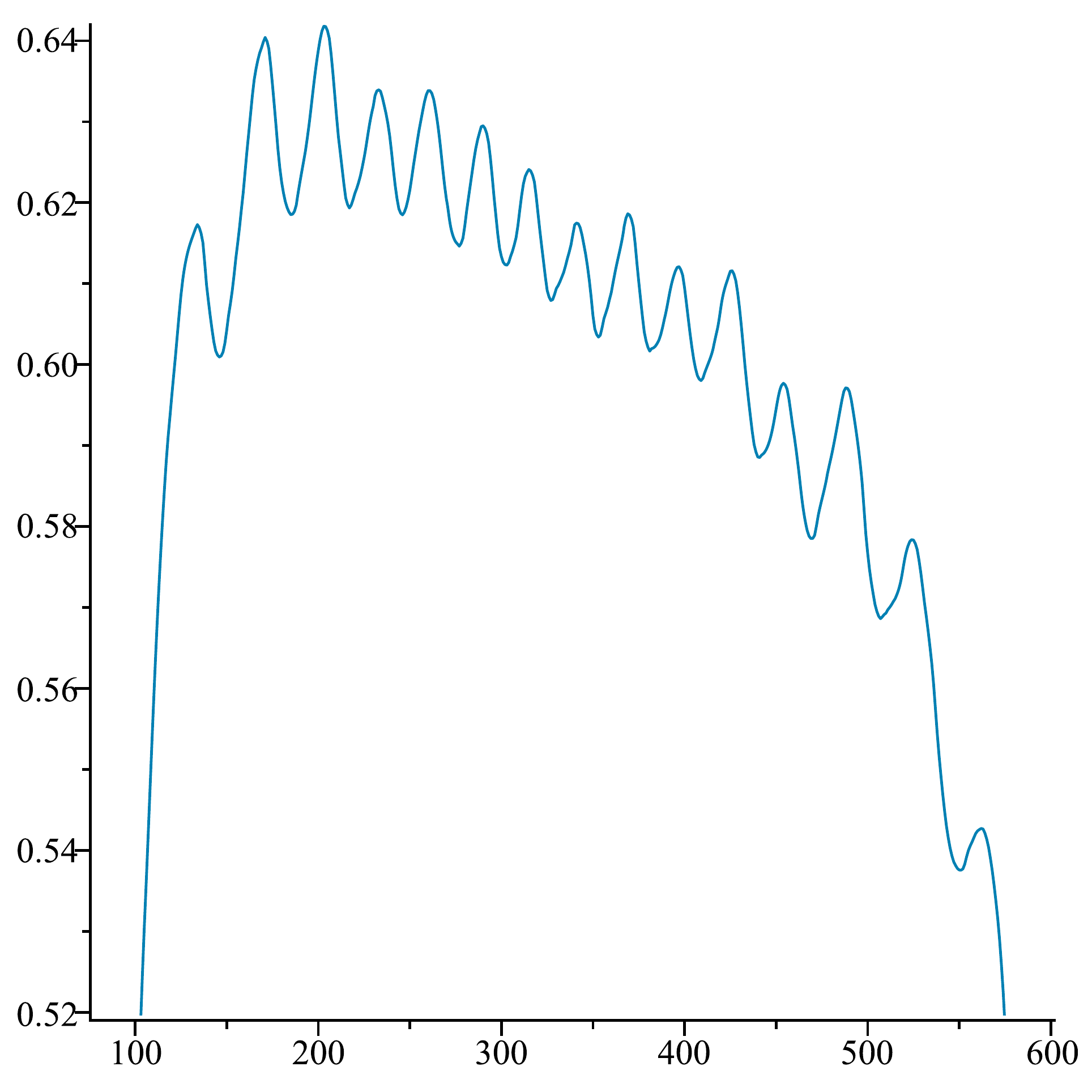}
 \includegraphics[width=0.45\textwidth]{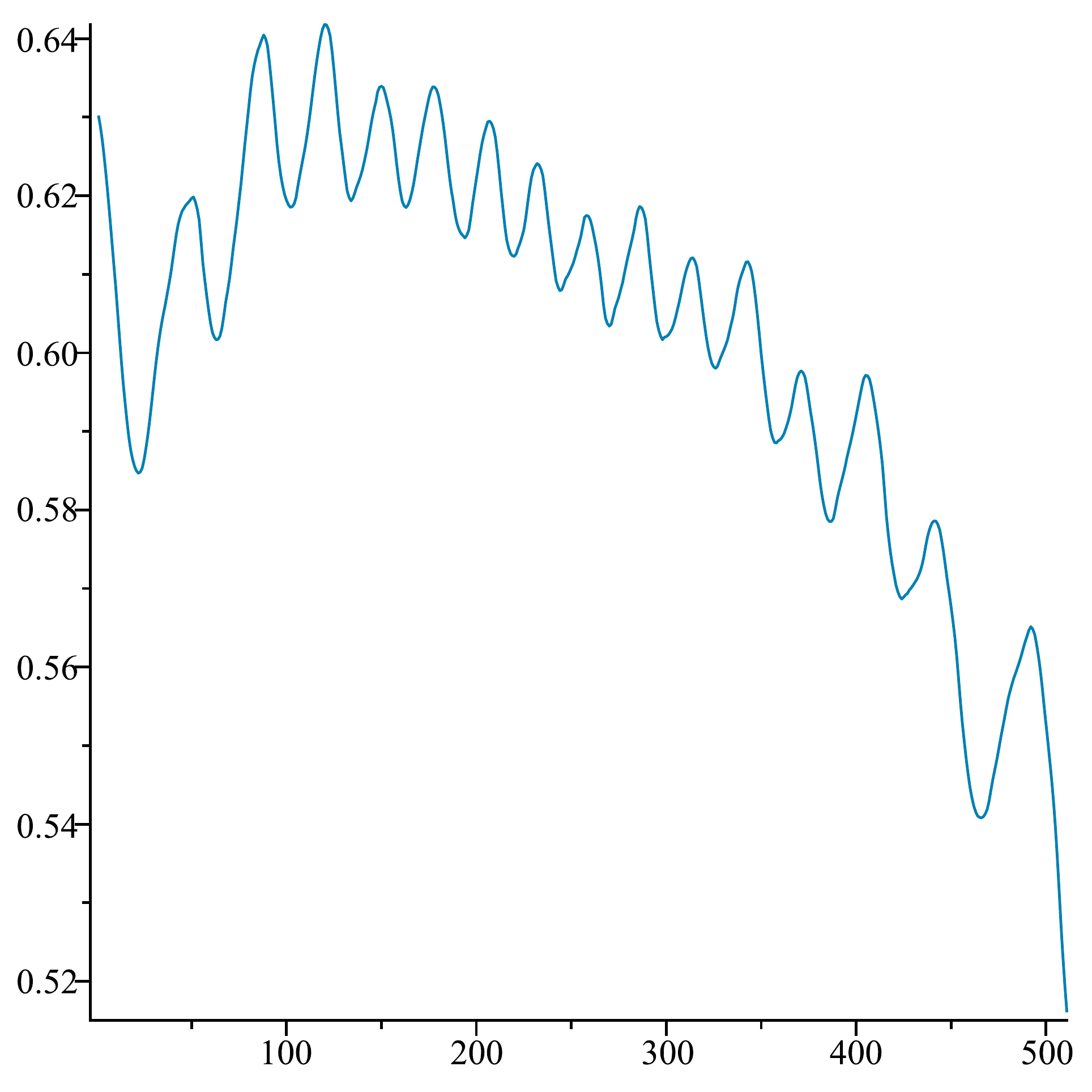}
\end{figure}

\begin{itemize}[leftmargin=1.5em]
\item Zun\"achst erkennt man, dass das Abtrennen der Nulll\"aufe von der Messreihe nur bedingt Auswirkung auf die Gl\"attung des relevanten Messbereiches hat.
\item Die im vorigen Punkt gemachte Einschr\"ankung ist aus folgenden Gr\"unden notwendig.
\begin{itemize}[leftmargin=1.5em]
\item  Bei einem hohen Nulllauf mit im Vergleich dazu schmalen Extrema werden die Randextrema gegebenenfalls weggegl\"attet, siehe Messreihe {\bf IV}, Abbildungen \ref{401}-\ref{405}. Da aber die Differenz zwischen der Anzahl der Z\"ahne auf der groben Seite der Scheibe einerseits und der feinen Seite andererseits relativ gross ist, w\"are der Verlust der Randextrema f\"ur die hier betrachtete praktische Fragestellung nicht relevant.
\item Auschlaggebend f\"ur die Entscheidung zur Abtrennung des Nulllaufs ist jedoch der folgende Punkt: In den Praxisdaten ist der Nulllauf nicht konstant, sondern er weist eine hochfrequente Grundschwingung auf, die auf eine Vibration des gesamten Versuchsaufbaus zur\"uckzuf\"uhren sein k\"onnte, vergleiche Abbildungen \ref{101}, \ref{201}, \ref{401} und \ref{501}. Diese Grundschwingung kann erst durch eine Erh\"ohung der Gleitl\"ange der Gl\"attung aufgefangen werden, und kann somit die Z\"ahlung beeinflussen. Vergleiche dazu auch die Diskussion in Abschnitt \ref{sub3}. 
\end{itemize}
\item Die symmetrisierte exponentielle Gl\"attung (SEA) gem\"a\ss~(\ref{SEA}) ist mit Abstand die effektivste Gl\"attungsmethode. Sie liefert schon bei kleiner Gleit\-l\"ange sehr gute Ergebnisse und auch bei einer hochfrequenten St\"orschwingung reicht eine moderate Gleitl\"ange zur Gl\"attung aus.
\item An Messreihe {\bf I} lassen sich die im Text beschriebenen Probleme der Gl\"at\-tungsverfahren im Zusammenhang mit stark zerkl\"ufteten Extrema erkennen, siehe Abbildungen \ref{101}-\ref{105} und dort speziell das vierte Maximum.
\item Der L\"osungsvorschlag der Vorgl\"attung mittels eines gleitenden Medians liefert hier die L\"osung. Vergleiche dazu in Abbildung \ref{106} die gezoomte Originalmessreihe (links) mit der Medianbildung (rechts), die dann als Ausgangsreihe herangezogen wird.
\end{itemize}

\subsection{Abbildungen zu den Praxisdaten}\label{Abbsec}\ 








\section{Fazit}\label{fazit}

Zur L\"osung der in der Einleitung beschriebenen technischen Fragestellung schlagen wir als Gl\"attungsverfahren die in Abschnitt \ref{subSEA} definierte symmetrisierte exponentielle Gl\"attung (SEA) vor. Um Ausrei\ss er abzufangen, sollte die Gl\"attung in Verbindung mit einer vorher durchgef\"uhrten Anwendung des gleitenden Medians durchgef\"uhrt werden. F\"ur beide sollte eine moderate Gleitl\"ange gew\"ahlt werden.

Eine a priori Bestimmung der Gleitl\"ange ist aufgrund der technischen Umsetzung und der damit willk\"urlichen Schwankungen unterliegenden Messreihen nicht m\"oglich. Vielmehr beruht seine Wahl auf Erfahrungswerten aus der Auswertung einer hinreichend gro\ss en Anzahl von Messreihen. 
Eine praktische Anwendung dieser empirischen a priori Festlegung der Gleitl\"ange kann etwa die Aussortierung defekter Bauteile sein. Die G\"ute der verwendeten Bauteile kann dann wiederum mit Hilfe einer Anpassung der Gleitl\"ange beeinflu\ss t werden.

Einerseits weist der SEA schon bei geringen Gleit\-l\"angen ein sehr gutes Gl\"at\-tungs\-verhalten auf und andererseits respektiert er die Symmetrie des Datensatzes. Der reine EA ist aus Gr\"unden der Unsymmetrie nicht geeignet, da zu fr\"uhe Artefakte nicht hinreichend gegl\"at\-tet werden k\"onnen, siehe Abschnitt \ref{sub2}, Abbildung \ref{fig23}. Die sehr guten Gl\"attungseigenschaften des EA kommen erst nach einer hinreichend langen Anfangsphase zum Tragen.
Trotz der guten Symmetrieeigenschaften ist der MA nicht in dem gleichen Ma\ss e geeignet wie der SEA, da man f\"ur vergleichbare Gl\"attung eine h\"ohere Gleitl\"ange ben\"otigt und somit eine starke Verk\"urzung des Datensatzes in Kauf nehmen muss.

Das Vorschalten des gleitenden Medians f\"angt zum einen Ausrei\ss er ein, zum anderen schafft er eine Vorabgl\"attung eventuell vorhandener zerkl\"ufteter Extrema. Mit dieser Form der Entartung haben alle Gl\"attungsverfahren -- wenn auch in unterschiedlichem Ma\ss e -- Probleme bei niedrigen Gleitl\"angen, vergleiche Abschnitt \ref{sub2} und ebenfalls die Messreihe II, Abbildungen \ref{101}-\ref{106}.

Zur Vorbereitung der Praxisdaten ist eine Abtrennung eines gegebenenfalls vorhandenen Nulllaufs angezeigt, da dieser die zur Gl\"attung notwendige Gleitl\"ange unn\"otig erh\"ohen kann, vergleiche Abschnitt \ref{sub7}.

\end{document}